  \documentclass[12pt]{article}
  \usepackage{amsfonts}
 \usepackage{amssymb,amsmath}

   \def\R{\mathbb{R}}
   \def\N{\mathbb{N}}
   \def\Z{\mathbb{Z}}

   \def\a{{\alpha}}
   \def\e{{\varepsilon}}
   \def\d{{\delta}}
   \def\D{{\nabla}}
   \def\l{{\lambda}}
   \def\vi{{\varphi}}
   \def\vt{{\theta}}

   \def\cA{{\cal A}}
   
   \def\cC{{\cal C}}
   \def\cD{{\cal D}}

   \def\cH{{\cal H}}
   \def\cI{{\cal I}}
   \def\cJ{{\cal J}}
   \def\cK{{\cal K}}
   \def\cL{{\cal L}}
   
   \def\cN{{\cal N}}
   
   \def\cR{{\cal R}}
   
   \def\cS{{\cal S}}

   \def\supp{\mathop{\rm supp}\nolimits}
   
   \def\dist{\mathop{\rm dist}\nolimits}
   \def\lo{\mathop{\longrightarrow}}

   \def\loc{\mathop{\rm loc}\nolimits}

   \def\no{\noindent}
   \def\qed{\hfill{\em q.e.d.}\\ \vspace{1mm}}
   \def\proof{\noindent{\underline {\sf Proof} \hspace{2mm}}}
   \def\tton{\stackrel{\mbox{\tiny{$n \to \infty$}}}{\longrightarrow}}
   \def\ttox{\stackrel{\mbox{\tiny{$|x| \to \infty$}}}{\longrightarrow}}

   \def\meas{\mathop{\rm meas}\nolimits}
   \def\intr{\mathop{\int_{\R^N}}}

   \newcommand{\beq}{\begin{equation}}
   \newcommand{\eeq}{\end{equation}}

\newtheorem{df}{Definition}[section]
\newtheorem{prop}[df]{Proposition}
\newtheorem{lemma}[df]{Lemma}
\newtheorem{teo}[df]{Theorem}
\newtheorem{rem}[df]{Remark}

\newtheorem{cor}[df]{Corollary}

 \newcommand{\sezione}[1]{\section{#1}\setcounter{equation}{0}}

\begin{document}


   \title{Multiplicity of positive and nodal solutions for scalar
     field equations} 
   
   \author{Giovanna Cerami\thanks{ Dipartimento di Matematica,
        Politecnico di Bari,  Via Amendola 126/B -  70126
     	Bari, Italia.  
		E-mail: {\sf  cerami@poliba.it}}
	\and
		Riccardo Molle\thanks{Dipartimento di
		Matematica, Universit\`a di Roma ``Tor Vergata", Via
                della Ricerca 
		Scientifica n$^o$ 1 - 00133 Roma. E-mail: {\sf
                  molle@mat.uniroma2.it}} 
	\and
		Donato Passaseo\thanks{Dipartimento di Matematica
                  ``E. De Giorgi'', 
	  	Universit\`a di Lecce P.O. Box 193 - 73100 Lecce,
                Italia.}		 
		}

 \date{}

 \maketitle

  

{\small {\sc \noindent \ \ Abstract.} - 
In this paper the question of finding infinitely many solutions to  
the problem $-\Delta u+a(x)u=|u|^{p-2}u$,  in $\R^N$,  $u \in  
H^1(\R^N)$, is considered when $N\geq 2$,   $p \in (2, 2N/(N-2))$,  
and  the potential $a(x)$ is a positive function which is not  
required to enjoy symmetry properties. 
Assuming that $a(x)$ satisfies a suitable ``slow decay at infinity''
condition and, moreover, that its graph has some ``dips'', we prove that
the problem admits either infinitely many nodal solutions either
infinitely many constant sign solutions. 
The proof method  is purely variational and allows to describe the
shape of the solutions. 
}

\vspace{2mm}


 \sezione{Introduction}


This paper deals with the question of the existence of multiple solutions for equations of the type

(E) \hspace{4cm}
 $-\Delta u+a(x) u=|u|^{p-1}u \ \ \ \ \ \ \  \mbox{in } \R^N,
    $

\vspace{2ex}    
\no where $N\geq2,\ p>1,$  $p< 2^* -1 = \frac{N+2}{N-2}$ if $N \geq
3,$ and the potential $a: \R^N \to \R $ is a positive smooth
function. 

Equation $(E)$ has strongly attracted the researchers attention and has
been extensively studied because it appears naturally in the study of
problems coming from applied Sciences. Actually, Euclidean Scalar
Field equations like $(E)$ arise in several contexts: the most known
is probably the search of certain kind of solitary waves in nonlinear
equations of the Klein-Gordon or Schr\"odinger type, but also questions
in nonlinear optics, condensed matter physics, laser propagation,
constructive field theory lead to look for solutions of equations like
$(E)$ (see f.i. \cite{BL}, \cite{Cmilan}). 

However, it is worth pointing out that a not least reason of interest
for $(E)$ rests on some of its mathematical features that give rise to
challenging difficulties. Equation $(E)$ is variational, so it is
reasonable to search its solutions as critical points of the
corresponding ``action'' functional defined in $H^1(\R^N)$ by  

$$I(u) = \frac{1}{2} \int_{\R^N}(|\D u|^2+a (x) u^2)dx-{1\over
  p+1}\int_{\R^N}|u|^{p+1}dx, $$
but a lack of compactness, due to the unboundedness of the domain
$\R^N,$ prevents a straight application of the usual variational
methods.  Whatever $p$ is, the embedding $j: H^1(\R^N)\rightarrow L^p
(\R^N)$ is not compact, hence arguments like those based on the
Palais-Smale condition can fail.  

 When the potential $a(x)$ enjoys radial symmetry the difficulty in
 treating $(E)$ can be overcome working in the subspace $H_r (\R^N)$
 of $H^1(\R^N)$ consisting of radially symmetric functions. $H_r
 (\R^N)$ embeds compactly in  $L^p (\R^N),\ p\in(2,2^*),$ so  standard
 variational techniques work and, as for bounded domains, one can show
 (\cite{S}, \cite{BL}, \cite{DN}) that $(E)$ possesses a positive
 radial solution and infinitely many changing sign radial
 solutions. Furthermore, the radial symmetry of the potential and, of
 course, of $\R^N$  plays a basic role also when one is
 looking for non-radially symmetric solutions (see f.i. \cite{BW},
 \cite{WY}).  
 
On the contrary, when $a(x)$ has no symmetry properties there is no
way to avoid the compactness question. Moreover, the difficulty is not
only a technical fact,  as one can understand considering $(E)$ with a
potential  increasing along a direction, then, as shown in \cite{CM},
$(E)$ has only the trivial solution.  Nevertheless, some considerable
progress has been performed also for non-symmetric potentials, even if
the situation is still not completely clear and the picture is
fragmentary. Most of researches have been concerned with the case in
wich 
$$
\lim_{|x|\to \infty} a(x) =  a_\infty > 0.
$$  

Under the above assumption various existence and multiplicity results
for $(E)$ have been stated, but it must be remarked at once that the
topological situation appears quite different according that $a(x)$
approaches $a_\infty$ from below or from above. In the first case the
existence of a least energy positive solution can be  obtained by
minimization, using concentration-compactness type arguments (see
\cite{PLL}, \cite{R}).  In the latter, the functional  $I$ may not
attain its infimum on the natural constraint $I'(u)[u] = 0$  and
critical points of $I$ have to be searched working at higher energy
levels and using more subtle variational and topological tools. By
this kind of arguments, in \cite{BaL} and \cite{BLi}, the existence of
a positive solution to $(E)$ has been proved under a suitable ``fast
decay'' assumption on $a(x)$. 

About the existence of multiple solutions, let us recall that a lot of
work has been done on the so called semi classical equation: 
$-\e^2 \Delta u+a(x) u=|u|^{p-1}u$.
Indeed, this equation, that by a change of variables can be written as
$(E)$ with potential $a(\e x),$ has been shown to have, when $\e$ is
small, a set of nontrivial positive solutions whose richness can be
related either to the topology of the sublevel sets of $a(x)$ either
to the number and/or the type of critical points of $a(x).$ But,
working with this kind of approach, it is not difficult to realize
that, as the number of solutions searched increases, the parameter
$\e$ goes to zero.  Making a list of all the interesting contributions
in this direction without forgetting something it is not an easy
matter, so we prefer to refer the interested reader to some  recent
papers on this subject \cite{BT}, \cite{LL}, and to \cite{DWY}
and to the rich list of references therein.

With respect to the question of getting, as in the radial case, infinitely many
solutions, again  different features appear
according the way the potential goes to its limit at infinity. 
 
When $a(x)$ goes to $a_\infty$ from below,  $(E)$  has been shown to
possess  infinitely many changing sign solutions in \cite{CDS}, by a
natural approximation method, under suitable ``slow'' decay 
assumptions on $a(x)$. 
Unlike, when $a(x)\ttox a_\infty$ from above, as far as we know, no
results of existence of infinitely many changing sign solutions are
until now known: namely, this is the main question we want to consider in
this paper. 
  
On the other hand, some  results of existence of infinitely many
positive ``multibump'' solutions hold in this situation. The first,
contained in \cite{CZR}, has been obtained under periodicity
assumptions on the potential, while, in the already mentioned paper
\cite{WY}, the attainment has been proved assuming $a(x)$ radially
symmetric and decaying at infinity with a prescribed polynomial
rate. The case in which $a(x)$ has no symmetry properties has been
successfully explored for the first time in \cite{CPS1}, where, by
using a purely variational method introduced in \cite{J1}, 
the following contribution  has been given: 

 \begin{teo}
\label{T1CPS}Let assumptions 

$$
\begin{array}{ccl}
(h_1) &\qquad&\qquad
a(x)\lo a_\infty > 0 \ \mbox{as} \  |x|\rightarrow \infty,\\
(h_2) &\qquad&\qquad
 a(x) \geq a_0 > 0 \   \forall x\in \R^N,\\
(h_3) &\qquad&\qquad
a \in L^{N/2}_ {loc}(\R^N),\\
(h_4) &\qquad&\qquad
\exists \ \bar{\eta} \in (0, \sqrt{a_\infty}) :\ \lim_{|x|\to \infty}
(a(x)-a_\infty) e^{\bar{\eta} 
 |x|} =  \infty 
\end{array}
$$ be satisfied.  

Then there exists a positive constant, $ \cA = \cA(N,\bar{\eta}, a_0, a_\infty) \in \R, $ such that, when 
$ \hspace{1cm}  (*) \hspace{2,5cm} |a(x)-a_\infty|_{{N/2},\,{loc}} := \sup_{y\in \R^N} |a(x) - a_\infty|_ {L^ {N/2}(B_1(y))}< \cA,$

\noindent equation (E) has infinitely many positive solutions belonging to $H^1(\R^N)$.
\end{teo}
The above result (see also the subsequent \cite{AW} for a different
proof and more general nonlinearities) is the starting point of our
work; some comments and questions come naturally  looking at its
statement. Indeed, assumptions $(h_1),-,(h_3)$ are standard and very
mild, moreover,  the slow decay condition $(h_4)$ is basic and it is
the deep motivation for the success in  obtaining ``multibump''
solutions.  
The solutions are found by a max-min argument on the
action functional $I$ restricted to special classes of multibump
functions and, in the maximization procedure, the attractive effect of
$a(x)$ on the bumps is dominating, just because of $(h_4),$ on the
repulsive disposition of the positive massess with respect to each
other and prevents the ``bumps'' to escape at infinity. 

On the contrary, the ``small oscillation'' condition $(*)$ appears
less natural, also
in view of the above quoted results in symmetric cases. 
Hence reasonably one wonders whether it is necessary or not and, more,
if it can be dropped. This last question is still largely open,
nevertheless in a recent paper \cite{DWY} it has been removed in the
planar case, under stronger regularity assumptions and imposing a
polynomial decay of the potential. 

In this paper we are also concerned with the above question, when the
potential satisfies $(h_1),-,(h_4)$ and no restriction on the dimension
of $\R^N$ is required.  Nevertheless, as already observed, our main
porpose is to face the problem of obtaining, by a variational approach
of the type  displayed in \cite{CPS1}, infinitely many ``multibump''
changing sign solutions to $(E)$ when $a(x) \ttox a_\infty$  from
above.  

Our results concern equations  in which the potentials are not
required to verify condition $(*)$, but, besides verifying
assumptions $(h_1),-,(h_4)$, they must have one or more regions around
the origin where they sink. The conclusion of such researches are
stated in the following two theorems:  

\begin{teo}
\label{T}
 Let $a(x)$ satisfy $(h_1)$,--,$(h_4)$ of Theorem \ref{T1CPS}.
Let $A,B\subset\R^N$ be bounded open sets such that $\ 0\in A\subset\subset
B$. 
Let  $\cN_A$ and $\cN_B$ be bounded neighborhoods of $\partial A$ and
$\partial B$, respectively, such that $\cN_A\cap\cN_B=\emptyset$. 
Let $ b:\cN \to\R,$ $\ \cN=\cN_A\cup\cN_B,$ be a continuous
nonnegative function such that
$$
\inf_{\partial A\cup\partial B} b =b_0>0,\qquad\qquad
\sup_{\cN}b<a_0.
$$
 Set $a_\e(x)= a(x) - b_\e(x)$ where $ b_\e : \R^N \to \R$ is defined as $  b_\e (x) =  b(\e x)$ if $x\in \cN/\e$, 
$  b_\e (x) = 0$ if $x \notin \cN/\e.$

Then   $\bar \e>0$ exists such that, for all $\e\in(0,\bar \e),$ 
$$
(P_\e)
\qquad\left\{\begin{array}{l}
 -\Delta u+a_\e(x) u=|u|^{p-1}u\\
     u\in H^1(\R^N)
\end{array}\right .
$$
has infinitely many nodal solutions.
\end{teo}

\vspace{2ex}

\begin{teo}
\label{T2}
Let $a(x)$ satisfy $(h_1)$,--,$(h_4)$ of Theorem \ref{T1CPS}. Let
$B\subset\R^N $ be an open bounded set such that $0\in B$. Let $\cN_B$
be a bounded 
neighborhood of $\partial B$ and let $\tilde b:\cN_B\to\R$ be a
continuous nonnegative function such that
$$
\inf_{\partial  B} \tilde b =\tilde b_0>0,\qquad\qquad
\sup_{\cN_B}\tilde b<a_0.
$$
 Set $\tilde a_\e(x)=  a(x) -\tilde b_\e(x)$ where $\tilde b_\e : \R^N
 \to \R$ is defined as $ \tilde b_\e (x) = \tilde b(\e x)$ if $x\in
 \cN_B/\e$,  
$ \tilde b_\e (x) = 0$ if $x \notin \cN_B/\e.$

Then  $\tilde \e>0$ exists such that, for all $\e\in(0,\tilde \e),$ problem
$$
(\tilde P_\e)
\qquad\left\{\begin{array}{l}
 -\Delta u+\tilde a_\e(x) u=|u|^{p-1}u\\
     u\in H^1(\R^N)
\end{array}\right .
$$
has infinitely many  constant sign solutions.
\end{teo}

\vspace{2ex}
\begin{rem}
{\em
It is worth stressing that the equations in $(P_\e)$ and  $(\tilde P
_\e)$ are exactly of the type $(E)$ for any choice of $\e$ and that
the claim of the above Theorems, unlike the above mentioned results
for semi classical equations, gives \textit{for all suitably small
  $\e$ infinitely many  solutions to the problems.} 

We notice also that Theorem \ref{T} still holds  assuming, on  
the bounded open sets $A,B,$ that $0\in A$ and $A\cap B=\emptyset$
instead of  $\ 0\in A\subset\subset 
B$.}
\end{rem}

In this research as in \cite{CPS1} we use variational tools:  the
solutions are found by looking, for all $k$, for critical points of
the functional $I_\e:H^1(\R^N)\to\R$ defined as
$$
I_\e(u)={1\over 2}\int_{\R^N}(|\D u|^2+a_\e(x)u^2)dx-{1\over
  p+1}\int_{\R^N}|u|^{p+1}dx
$$
in special classes of functions having $k$ ``bumps''. 
Such  functions, that are introduced in Section 2, are called
``functions emerging around $k$ points'' and are functions having the
most important part of their mass split in $k$ parts supported in $k$
disjoint balls. 
Of course, when we are interested in constant sign solutions we consider
only positive masses, while when we look for changing sign solutions
we need to consider functions having positive and negative bumps; in
this last case  a new difficulty arises: the possibility that ``bumps''
of opposite sign  collapse, because masses of opposite sign exert a
mutual attractive effect. The role played by the
potential dips in some regions of $\R^N$ is essential to control this
phenomenon.

\begin{rem}
{\em We point out that in section 2, after defining the  classes of
  ``admissible'' functions where the solutions are searched, Theorems
  \ref{T} and \ref{T2} are restated in Theorems
  \ref{T1} and \ref{T3} respectively in a more precise way,
  highlighting how the solutions 
  are actually constituted. We also remark that  in Section 2 some
  slightly more general result about the existence of nodal solutions
  are stated in Theorems \ref{T4} and \ref{T?}}. 
\end{rem}

Finally we observe that, making a deep study of the behavior of the
solutions found in Theorems \ref{T} and \ref{T2},  as the number of
the ``bumps'' increases going to infinity, it would be possible to show,
by arguments similar to those developed in \cite{CPS2},  that, in both
cases, those solutions converge to a solution having infinitely many bumps.  

\vspace{2ex}

The paper is organized as follows: in section 2 some notation is
introduced, useful facts and preliminary results are stated; in
sections 3 and 4 a max-min procedure is developed with the aim of
finding, in the special classes of functions emerging around $k$
points, for all $k\in \N,$ good candidates  to be critical points;
section 5 contains the analysis of the asymptotic properties (with
respect to shape and position) of the emerging part of the expected
solutions as $\e \to 0;$ finally in section 6 the proof of Theorems
\ref{T} and \ref{T2} is completed.


\sezione{Notations and variational framework}


Throughout the paper we make use of the following notations:

{\small

\begin{itemize}
\item
$H^1(\R^N)$ denotes the closure of $\cC^\infty_0(\R^N)$ with respect
to the norm
$$\|u\|:=\left[\int_{\R^N}(|\D u|^2+ a_\infty u^2)dx\right]^{1/2}.$$

\item
$L^q(D)$, $1\leq q \leq +\infty$, denotes the Lebesgue space with the
norm $|\cdot|_{L^q(D)}$ ; the norm in $L^q(\R^N)$ is denoted by $|\cdot|_q.$
\item
$B_r(x), $ when $x$ belongs to a metric space $X$, denotes the open
ball of $X$ having radius $r$ and centered at $x$.
\item
$L^q_{loc}(\R^N),$ $1\leq q < +\infty,$ [respectively
$H^1_{loc}(\R^N)$] denotes the set of functions $u$ such that 
$u \in L^q(B_1(x))$ $\forall x \in \R^N$ [$u \in H^{1,2}(B_1(x))$].

\noindent In $L^q_{\loc}(\R^N)$ we consider the subspace of
functions such that $\sup_{x\in \R^N} |u|_ {L^q(B_1(x))}$ is finite, endowed with
the norm
$$ 
|u|_{q,{\loc}} := \sup_{x\in \R^N} |u|_ {L^q(B_1(x))}.
$$

\item
Given a Lebesgue-measurable subset $\cD$ of $\R^N$, $|\cD|$ denotes
its Lebesgue measure. 

\item
C, c, $\tilde{c},\ \hat{c},\ c_i$ denote various positive constants.
\end{itemize}
}

We also set $\alpha_\e(x)=a_\e(x)-a_\infty$, $\bar a=a_0-\sup_\cN b$.

In what follows we consider 
$I_\infty:H^1(\R^N)\to \R$ defined by 
$$
I_\infty(u)={1\over 2}\int_{\R^N}(|\D u|^2+a_\infty u^2)dx-{1\over
  p+1}\int_{\R^N}|u|^{p+1}dx.
$$

The following three lemmas contain some useful facts, for the first
two, which are well known, see f.i. 
\cite{Cmilan} and \cite{AM} respectively, the third one is Lemma 3.3
of \cite{CPS1}. 

\begin{lemma}
\label{N2}
The problem
$$
(P_\infty)
\qquad\left\{\begin{array}{l}
 -\Delta u+a_\infty u=|u|^{p-1}u\\
     u\in H^1(\R^N)
\end{array}\right .
$$
 has a positive ground state solution $w$, unique up to translation,
 radially symmetric, 
decreasing when the radial coordinate increases and such that
\begin{equation}
\label{N2*}
\lim_{|x|\to\infty}|D^j w (x)||x|^{\frac{N-1}{2}}e^{\sqrt{a_\infty}|x|}=d_j>0,
\hspace{3mm} d_j\in\R,\ j=0,1.
\end{equation}
\end{lemma}

Moreover, setting
$$m_{\infty} := I_{\infty} (w) ,$$
$m_{\infty}$ can be characterized as:
\beq
\label{2.1}
m_{\infty} = \min_{\gamma \in \Gamma}\max_{t\in [0,1]} I_{\infty}(\gamma (t)),
\eeq
where, denoting by $e\in H^1( \R^N)$ any point for which $I_{\infty}(e)<0$,
\beq
\label{Gamma}
 \Gamma =\left\{\gamma \in C([0,1], H^1( \R^N)) : \gamma(0)=0,\
   \gamma(1)=e\right\}.
\eeq

\begin{lemma}
\label{N2'}
Set, for all $y \in \R^N,$
$$
 w_y (x) = w(x-y), 
$$
and
$$
\Phi = \left\{ w_y : \ y\in\R^N\right\}.
$$

Then $\Phi$ is nondegenerate, namely the following properties are true:

i)\ $(I_{\infty})''(w_y)$ is an index zero Fredholm map, for all $y \in \R^N$;

ii) Ker $(I_{\infty})''(w_y)$ = span $\left\{ \frac{\partial w_y}{\partial x_j}: 1 \leq j \leq N \right\} = T_{w_y}\Phi,$

$T_{w_y}\Phi$ being the tangent space to $\Phi$ at $w_y.$
\end{lemma}

\begin{lemma}
\label{N*} Let $\cD \subset \R^N$ be closed, let  $\lambda, s \in \R^+\setminus \left\{0\right\}$ and let $u$  be a positive solution of
$$
 \hspace{1cm} \left\{
        \begin{array}{ll}
         -\Delta u+ \lambda u \leq 0         & \mbox{ in }\ \R^N \setminus \cD,\\
          u \leq s                     & \mbox{ in }\ \R^N \setminus \cD.
        \end{array}
\right.
$$
Then, for all $d \in (0, \sqrt{\lambda}), $ there exists a positive
constant $ C_d$ (depending on $\lambda , d, N$) such that 
$$
u(x)\leq C_d s e^{-d\, dist(x,\cD)} \ \ \ \ \ \ \ \ \ \ \  \forall x \in \R^N \setminus \cD.
$$
\end{lemma}

\vspace{2ex}

In what follows we denote by $\delta$ a suitably small, fixed real
number such that
\beq
\label{3.1}
0<\delta<\min\left\{1,\left(\bar a\over
    p\right)^{1\over p-1},{w(0)\over 3},(a_\infty-\bar\eta ^2)^{1\over
  p-1}\right\}
\eeq
and by $\rho>0$ a real positive number fixed so that 
$$
w(x)<\d\qquad \forall x\in \R^N\setminus
B_{\rho/2}(0).
$$ 
For all $u\in H^1(\R^N)$, we can write 
$$
u=u_\delta +(u^+)^\delta-(u^-)^\delta,
$$
where 
$$
(u^+)^\delta=(u-\delta)^+,\qquad\qquad
(u^-)^\delta=(-u-\delta)^+,\qquad\qquad u_\delta=(u\wedge\delta)\vee(-\delta)
$$
are called, respectively, the {\em emerging positive part}, the {\em emerging
 negative  part} and the {\em middle part} of $u$.
Fixing $\delta$ and $\rho$ as before indicated, $l,m\in\N$, and $l+m$
points in $\R^N,$ $x_1,\ldots,x_l,y_1,\ldots,y_m,$ 
having interdistances not less than $2 \rho,$ we say that
a function $u\in H^1(\R^N)$ is
positively emerging around $(x_1,\ldots,x_l)$ and
negatively emerging around $(y_1,\ldots,y_m)$ (in  balls of radius $\rho$) if 
$$
(u^+)^\delta=\sum_{i=1}^l(u_i^+)^\delta,\qquad
(u^-)^\delta=\sum_{i=1}^m(u_i^-)^\delta,
$$
where
$$
(u^+_i)^\delta\in H^1_0(B_\rho(x_i)),\ (u^+_i)^\delta\ge 0,\
(u^+_i)^\delta\not\equiv 0,\ i\in\{1,\ldots,l\},
$$
$$
(u^-_i)^\delta\in H^1_0(B_\rho(y_i)),\ (u^-_i)^\delta\ge 0,\
(u^-_i)^\delta\not\equiv 0, \ i\in\{1,\ldots,m\}.
$$

\vspace{2ex}
Clearly, if $u\geq 0,$ we say that it is positively emerging around
$l$ points   $(x_1,\ldots,x_l)$ (in  balls of radius $\rho$) if its
emerging (positive) part  is such that
$(u^+)^\delta=\sum_{i=1}^l(u_i^+)^\delta,$ and the functions
$(u_i^+)^\delta \in H^1_0(B_\rho(x_i))$ are as  above described;
analougously,  we say that $u \leq 0$ is negatively emerging around
$m$ points  when $(u^-)^\delta=\sum_{i=1}^m(u_i^-)^\delta$ with
$(u^-_i)^\delta\in H^1_0(B_\rho(y_i))$ satisfying the above
relations. 
 
\vspace{2ex}
\textit{We are now ready to introduce the classes of multi-bump
  functions where we look for solutions of $(P_\e)$.} 

First of all we fix a neighborhood $\cN_0\subset\subset\cN$ of $\partial
A\cup\partial B$ such that $\inf_{\cN_0}b>0$
and, for $l,m\in\N\setminus\{0\}$, $\e>0$,
we define 
\begin{eqnarray*}
\nonumber\cH_l^\e=\{(x_1,\ldots,x_l)\in(\R^N)^l 
& : & 
 x_i\in(B\setminus A)/\e, \ \dist(x_i,\cN_0/\e)\ge 3\rho,\\
& & \ |x_i-x_j|\ge 3\rho,\
i,j\in\{1,\ldots,l\},  i\neq j\};  
\\
\nonumber \cK_m^\e=\{(y_1,\ldots,y_m)\in(\R^N)^m 
& : &
y_i\in(\R^N\setminus B)/\e,\ \dist(y_i, \cN_0/\e)\ge 3\rho, \\
& & |y_i-y_j|\ge 3\rho,\  i,j\in\{1 ,\ldots,m\}, i\neq j\}. 
\end{eqnarray*}

Then, for all  $(x_1,\ldots,x_l)\in \cH_l^\e$ and $(y_1,\ldots,y_{m})\in
\cK_{m}^\e,$ we set
\begin{eqnarray*}
S^\e_{x_1,\ldots,x_l,y_1,\ldots,y_m}\hspace{-2mm} &=\{&\hspace{-2mm}
u\in H^1(\R^N)\  :\  u \mbox{ 
  positively emerging around }(x_1,\ldots,x_l)\in \cH_l^\e,\\
& & \mbox{negatively emerging around }(y_1,\ldots,y_m)\in \cK_m^\e,
\mbox{ such that }\\
& &
 I'_\e(u)[(u_i^+)^\delta]=0,\ i\in\{1,\ldots,l\},\
 I'_\e(u)[(u_j^-)^\delta]=0,\ j\in\{1,\ldots,m\},\\
& & \beta_{z}(u)=0, \ \forall z\in\{x_1,\ldots,x_l,y_1,\ldots,y_m\}\}
\end{eqnarray*}
where
$$ 
\beta_z(u)\ =\ \frac{1}{\int_{B_\rho(z)}((|u|^+)^\delta)^2} \int_{B_\rho(z)} (x-z)
 ((|u(x)|^+)^\delta)^2dx.
$$

\vspace{2ex}

If $u\ge 0$ is positively emerging around $(x_1,\ldots,x_l)\in
\cH_l^\e$ [$u\le 0$ is negatively emerging around $(y_1,\ldots,y_m)\in
\cK_{m}^\e$], we can define analogously the set
$S^{\e,+}_{x_1,\ldots,x_l}$ [ $S^{\e,-}_{y_1,\ldots,y_m}$].

\vspace{2ex}

Obviously, {\em for all $l \in \N \setminus \left\{0\right\}$  an  $\e_l>0$
exists such that $\e\in(0,\e_l)$ implies 
  $\cH_l^\e\neq \emptyset.$   When $l= 0,\ \e_0$ can be defined arbitrarily, so we set $\e_0=\infty$.
 We remark  that from
$\cH_{l}^\e\neq \emptyset,$  it follows $\cH_{i}^{\e}\neq \emptyset$
for all $i \in \left\{1,...,l-1\right\}$.}

\vspace{2ex}

The solutions, whose existence is claimed in Theorem \ref{T} and
\ref{T2}  are searched fixing $h \in \N\setminus\left\{0\right\}$ and
$h=0$ respectively, and looking for critical points of $I_\e$ among
functions which are positively emerging around at most $h$ points and
negatively emerging around points whose number becomes arbitrarily
large. Namely, Theorem \ref{T} can be restated as follows:  

\begin{teo}
\label{T1}
Let  assumptions of Theorem \ref{T} be satisfied. Then, for all
$h\in \N\setminus\{0\}$ there exists $\bar \e_h>0$ such that for all
$\e\in(0,\bar \e_h)$, for all $k\in \N\setminus\{0\}$ and for all
$l\in\N$, $1\le l\le\min 
(h,k)$, there exists a solution of $(P_\e)$ positively emerging around
$l$ points  and negatively emerging around $k-l$ points.
\end{teo}
 Theorem \ref{T2} is obtained proving the following statement:

\begin{teo}
\label{T3}
Let assumptions of Theorem \ref{T2} be satisfied. Then,   $\bar \e \in \R \setminus \left\{0\right\}$ exists such that for all
$\e\in(0,\bar \e)$, and for all $k\in \N\setminus\{0\}$ there exists
  a negative solution of $(\tilde{P} _\e),$ negatively emerging around
$k$ points.
\end{teo}

\begin{rem}
{\em It is worth pointing out that in Theorems \ref{T1} and \ref{T3}
  the role of the positive and negative bumps can be exchanged. So,
  the same kind of arguments, we are going to develope, for proving
  Theorems \ref{T1} and  \ref{T3} could bring to conclude respectively
  that $(P_\e)$ possess infinitely many solutions, which are
  negatively emerging around at most $h$ points and positively
  emerging around points whose number becomes arbitrarily large  and
  that $(\tilde{P} _\e)$ has infinitely many positive solutions.} 
\end{rem}
\begin{rem}
{\em When the connected components of  $B\setminus \bar A$ are more than one, 
  more general multiplicity results can be proved. 
Actually, we might set on each connected component a prescribed number
of positive or negative bumps obtaining for instance:}
\end{rem}

\begin{teo}
\label{T4}
Let assumptions of Theorem \ref{T} be satisfied. 

Let $B\setminus \bar
A=\cup_{i=1}^d C_i$, $d\in\N\setminus\{0\}$,
$C_i$ being open sets such that   
$C_i\setminus\bar \cN\neq\emptyset$, $\bar C_i\cap \bar C_j=\emptyset$
if $i\neq j, \ \forall i, j \in\{1,\ldots,d\}$. 
  
Then, for all $h\in \N$ there exists $\bar \e_h>0$ 
such that for all $\e\in(0,\bar \e_h),$ for all $\bar{d}
\in\{0,1,\ldots,d\}$   for all $k \in \N\setminus\{0\},$ 
 for all $h_i, m_j \in\N$
$i\in\{1,\ldots,\bar{d}\}, \ j\in\{\bar{d}+1,\ldots,d +1\}$ such that
$\sum_{i=1}^{\bar{d}}h_i =h$, $\sum_{j=\bar{d}+1}^{d+1}m_j =k-h,$
problem $(P_\e)$ has a solution  
 positively emerging around $h$ points  and  negatively emerging
 around $k-h$ points: $h_i$ points in $C_i$,  $i\in\{1,\ldots,\bar{d}\},$ 
 $m_j$ points in $C_j$, $j\in\{\bar{d}+1,\ldots,d \}$ and $m_{d+1}$ points
 in $\R^N\setminus \bar B$. 
\end {teo}

\begin{teo}
\label{T?}
Let $a(x)$ satisfy $(h_1)$,--,$(h_4)$ of Theorem \ref{T1CPS}. 
 Let $B_i$ be open sets such that $0\in
B_0$ and  $B_{i-1}\subset\subset B_{i}$,  $i\in\{1,\ldots,d\}$, $d\in
\N\setminus\{0\}$. 
Let $\cN$ be a bounded neighborhood of $\cup_{i=0}^d \partial B_i$
such that $B_i\setminus(\bar\cN\cup B_{i-1})\neq \emptyset$, for
$i=1,\ldots,d$. 
Let $ b:\cN\to\R$ be a continuous nonnegative function such that
$$
\inf_{\cup_{i=0}^d \partial B_i}  b = b_0>0,\qquad\qquad \sup_{\cN} b<a_0.
$$     
 
Then, for all $h\in \N$ there exists $\bar \e_h>0$ 
such that for all $\e\in(0,\bar \e_h),$ for all $\bar{d}
\in\{0,1,\ldots,d\}$   for all $k \in \N\setminus\{0\},$ 
 for all $h_i, m_j \in\N$
$i\in\{1,\ldots,\bar{d}\}, \ j\in\{\bar{d}+1,\ldots,d +1\}$ such that
$\sum_{i=1}^{\bar{d}}h_i =h$, $\sum_{j=\bar{d}+1}^{d+1}m_j =k-h,$
problem $(P_\e)$ has a solution  
positively emerging around $h$ points  and  negatively emerging
around $k-h$ points: $h_i$ points in $B_{i}\setminus \bar B_{i-1}$,
$i\in\{1,\ldots,\bar{d}\},$  
$m_j$ points in $B_{i}\setminus \bar B_{i-1}$,
$j\in\{\bar{d}+1,\ldots,d \}$ and $m_{d+1}$ points 
in $\R^N\setminus \bar B_d$. 

\end {teo}

 It is clear that, for proving Theorems \ref{T4} and \ref{T?},
  the definition of the sets $\cH^\e_l$ 
 and  $\cK^\e_m$ as well as some arguments we are going to develope
 must be suitably refined, but this would imply only some technical
 complications and  heavier notation.

\vspace{2mm}

For what follows it is useful to define, on the subset of $H^1(\R^N)$
consisting of functions having finite measure support, the functional 
\begin{eqnarray*}
J_\e (u) = \frac{1}{2}\int_{\R^N}(|\nabla u|^2 +a_\e(x)u^2)dx + 
\int_{ \R^N} a_\e(x) \delta \, u \, dx  \hspace{6ex} 
\nonumber \\
- \frac{1}{p+1}\int_{\supp\, u} (\delta + |u|)^{p+1} dx +
\frac{1}{p+1}\ \delta ^{p+1} |\supp\ u|.  
\end{eqnarray*}
 We can then write for all $u\in H^1(\R^N)$ 
\begin{eqnarray}
I_\e(u)&=&I_\e(u_\delta)+J_\e((u^+)^\delta)+J_\e((u^-)^\delta) 
\nonumber \\
& = & I_\e(u_\delta-(u^-)^\delta)+J_\e((u^+)^\delta)
=I_\e(u_\delta+(u^+)^\delta)+J_\e((u^-)^\delta).
\label{1215}
\end{eqnarray}
The aim of next part of this section is to show that the above defined
sets $S^\e_{x_1,\ldots,x_l,y_1,\ldots,y_m}$ are not empty, to describe
the nature of the ``natural non smooth constraints''
$I'_\e(u)$ $[(u_i^+)^\delta]=0$ and  $I'_\e(u)[(u_j^-)^\delta]=0$ as
well as some important feature of the functionals $I_\e$ and $J_\e$ on
these sets. 
\begin{lemma}
\label{L3.1}
Let assumptions $(h_1),(h_2)$, and $(h_3)$ be satisfied. 
For all $\e>0$ and for all $u \in H^1 (\R^N)$ such that $(u^+)^\delta\neq 0$ 
 the function $\cJ^+_{\e,u}: [0,+\infty) \lo\R$ 
[for all $u \in H^1 (\R^N)$ such that $(u^-)^\delta \neq 0$ the
function $\cJ^-_{\e,u}: [0,+\infty) \lo\R$] defined as 
$$
\cJ^+_{\e,u}(t) = I_\e (u_\delta + t (u^+)^\delta-(u^-)^\delta)
\qquad\qquad
[\cJ^-_{\e,u}(t) = I_\e (u_\delta + (u^+)^\delta-t(u^-)^\delta)],
$$
has in $(0,+\infty)$ a unique critical point, which is a maximum. 
\end{lemma}

\proof
 Setting $\cS := \supp  (u^+)^\delta$ and using (\ref{1215}) we have
\begin{eqnarray}
(\cJ^+_{\e,u})' (t)  & = & \frac{d}{dt}\; J_\e (t(u^+)^\d) 
\nonumber \\ 
& = &t \left[ \int_{\cS}(|\nabla (u^+)^\delta|^2 +a_\e(x)((u^+)^\delta)^2)dx \right]
+ \delta \int_{ \cS} a_\e(x)(u^+)^\delta dx 
\nonumber \\ 
& & - \int_{ \cS}(\delta + t (u^+)^\delta)^{p} (u^+)^\delta dx.
\qquad \label{4.1}
\end{eqnarray}

Now
$$ 
\lim_{t \rightarrow +\infty} (\cJ^+_{\e,u})'(t)  = - \infty
$$
and,  by the choice of $\d$ in (\ref{3.1}),
\beq
\label{1836}
(\cJ^+_{\e,u})' (0) \geq  \delta \ [\bar a - \delta ^{p-1}]\ \int_{ \cS}
(u^+)^\delta dx  >  0, 
\eeq
moreover, it is easy to see that $(\cJ^+_{\e,u})'''(t)\le 0$, 
for all $t \in (0,+\infty),$ so 
$(\cJ^+_{\e,u})'(t)$ is concave and the claim for $\cJ^+_{\e,u}$
follows.
The same argument clearly works for $\cJ^-_{\e,u}$. 

\qed

 \textit{The maximum of $\cJ^+_{\e,u}$ in $(0,+\infty)$ is denoted by $\theta_\e^+(u)$, the maximum of $\cJ^-_{\e,u}$ in $(0,+\infty)$ is denoted by $\theta_\e^-(u)$ }.

\begin{cor}
\label{C5.1}
Let assumptions $(h_1),(h_2)$,and $(h_3)$ be satisfied.
 Let $l,  m\in\N$, 
and let $u\in
H^1(\R^N)$ be a function  positively emerging around
$(x_1,\ldots,x_l)\in \cH_l^\e$ and negatively emerging around
$(y_1,\ldots,y_m)\in \cK_m^\e$. 

Let consider, when $l>0$, for every $s\in\{1,\ldots,l\},$ the function
$$
t\mapsto I_\e\left(u_\delta+t(u^+_s)^\delta+\sum_{i\neq s}(u_i^+)^\delta-
\sum_{j=1}^m(u_j^-)^\delta\right),
$$
and, when $m>0$, for every $ s\in\{1,\ldots,m\},$ the function
$$
t\mapsto I_\e\left(u_\delta+\sum_{i=1}^l(u_i^+)^\delta-
t(u_s^-)^\delta-\sum_{j\neq s}(u_j^-)^\delta\right).
$$
Then each one of these functions has in $(0,+\infty)$ a  unique critical point, 
 which is a maximum point. Such points are denoted by  $\vt^{\e,+}_s(u)$ and $\vt^{\e,-}_s(u)$ respectively.
\end{cor}

\proof
The argument is quite analogous to that of Lemma \ref{L3.1} once one
writes 
$$
I_\e\left(\hspace{-1mm} u_\delta+t(u^+_s)^\delta+\sum_{i\neq s}(u_i^+)^\delta-
\sum_{j=1}^m(u_j^-)^\delta\hspace{-1mm}\right)=
I_\e\left(\hspace{-1mm}u_\delta+\sum_{i\neq s}(u_i^+)^\delta-
\sum_{j=1}^m(u_j^-)^\delta\hspace{-1mm}\right)+J_\e (t(u^+_s)^\d)
$$
and, respectively,
$$
I_\e\left(\hspace{-1mm}u_\delta
+\hspace{-1mm}\sum_{i=1}^l(u_i^+)^\delta\hspace{-1mm}- 
t(u_s^-)^\delta\hspace{-1mm}-
\hspace{-1mm}\sum_{j\neq s}(u_j^-)^\delta\hspace{-1mm}\right)=
I_\e\left(\hspace{-1mm}u_\delta
+\hspace{-1mm}\sum_{i=1}^l(u_i^+)^\delta\hspace{-1mm}
-\hspace{-1mm}\sum_{j\neq s}(u_j^-)^\delta\hspace{-1mm}\right)
+J_\e (t(u_s^-)^\delta ).
$$

\qed

\begin{df}
\label{D1307}
Let $u\in H^1(\R^N)$ be such that $(u^+)^\delta\neq 0,$ then  $u_\delta+\theta_\e^+(u)(u^+)^\delta-(u^-)^\delta$ is called the \emph{ projection} of $u$ on the set $\{u \in H^1 (\R^N) : I_\e'(u)[(u^+)^\delta] = 0\}.$ 

Let $u\in H^1(\R^N)$ be such that $(u^-)^\delta\neq 0,$ then  $u_\delta+(u^+)^\delta-\theta_\e^-(u)(u^-)^\delta$ is called the
 \emph{ projection}  of $u$ on the set $\{u \in H^1 (\R^N) : I_\e'(u) [(u^-)^\delta] =
0\}$. 

Let $u\in H^1(\R^N)$ be positively emerging around
$(x_1,\ldots,x_l)\in \cH^\e_l $ and negatively emerging around 
$ (y_1,\ldots,y_m)\in \cK^\e_m$. 

Let be $s \in \left\{1,\ldots,l\right\},$ the function
$u_\delta+\vt^{\e,+}_s(u)(u^+_s)^\delta+\sum_{i\neq s}(u_i^{+})^\delta- 
\sum_{j=1}^m(u_j^-)^\delta$ is called the \emph{projection} of $u$ on the set
$\{u\in H^1(\R^N)\ :\ I_\e'(u)[(u^+_s)^\delta]=0\}.$ 

Let be $s \in \left\{1,\ldots,m\right\}$,
$u_\delta+\sum_{i=1}^l(u_i^{+})^\delta- \vt^{\e,-}_s(u)(u^-_s)^\delta- 
 \sum_{j\neq s}(u_j^-)^\delta$ is called the \emph{projection} of $u$ on the set
$\{u\in H^1(\R^N)\ :\ I_\e'(u)[(u^-_s)^\delta]=0\}$. 
\end{df}

\begin{rem}
\label{R2.5}
{\em
We point out that, if $u\in S^\e_{x_1,\ldots,x_l,y_1,\ldots,y_m}$, then
the projection of $u$ on the sets $\{u\in H^1(\R^N)\ :\
I'_\e(u)[(u^+_s)^\delta]=0\}$ and $\{u\in H^1(\R^N)\ :\
I'_\e(u)[(u^-_s)^\delta]=0\}$ is nothing but $u$ itself.
Furthermore, we stress that $u\in S^\e_{x_1,\ldots,x_l,y_1,\ldots,y_m}$ implies
$$
\max_{t\geq 0}I_\e\left(\hspace{-1mm} u_\delta+t(u^+_s)^\delta+\sum_{i\neq s}(u_i^+)^\delta-
\sum_{j=1}^m(u_j^-)^\delta\hspace{-1mm}\right)=I(u)
$$
\beq
\label{punto} = 
\max_{t\geq 0} I_\e\left(\hspace{-1mm}u_\delta+\sum_{i=1}^l(u_i^+)^\delta-
t(u_s^-)^\delta-\sum_{j\neq s}(u_j^-)^\delta\hspace{-1mm}\right).
\eeq
}
\end{rem}

\begin{lemma}
\label{L5.2}
Assume $(h_1),(h_2)$, and $(h_3)$ hold. Then for all
$l,m\in\N\setminus\{0\}$ and for all
 $(x_1,\ldots,x_l)\in\cH^\e_l$, $(y_1,\ldots,y_m)\in\cK_m^\e$ 
\beq
\label{*}
\begin{array}{c}
\vspace{2mm} 
(i)\ S^{\e}_{x_1,\ldots,x_l,y_1,\ldots,y_m}\neq\emptyset,\quad 
(ii)\ S^{\e,+}_{x_1,\ldots,x_l}\neq\emptyset,\qquad \forall \e\in (0,\e_l)\\ 
(iii)\ S^{\e,-}_{y_1,\ldots,y_m}\neq\emptyset,\qquad \forall  \e>0.
\end{array}
\eeq
\end{lemma}
\proof
Let $v \in \cC_0^\infty(B_\rho(0))$ be a positive, radially
symmetric (around the origin) function such that $v (x)> \d$ on a
positive measure subset of $B_\rho(0)$. 
Let be $l\neq 0$, $m\neq 0$ and set, for all $i \in
\left\{1,...,l \right\}$ and $j\in\{1,\ldots,m\}$,  $v_{x_i}(x)=
v(x-x_i)$ and $v_{y_j}(x)= v(x-y_j)$.
Then denote by
$ \tilde{v}_{x_i}$ the projection of $v_{x_i}$ on the set
$\{u \in H^1 (\R^N) : I_\e'(u)
[(u^+)^\delta] = 0\}$ 
and by $ \tilde{v}_{y_j}$  the projection of $-v_{y_j}$ on $\{u \in
H^1 (\R^N) : I_\e'(u) [(u^-)^\delta] = 0\}$. 

Then the function
\[
u(x) = \ \ \left\{
        \begin{array}{ll}
         0          &  \forall\  x \in \ \R^N \setminus
         [(\bigcup^l_{i=1} \supp v_{x_i})\cup
(\bigcup^m_{j=1} \supp v_{y_j})] \\
          \tilde{v}_{x_i}(x) &  \forall\  x \in \ \supp v_{x_i}, \
          i\in\{1,\ldots,l\} \\
\tilde{v}_{y_j}(x) &  \forall\  x \in \ \supp v_{y_j}, \
          j\in\{1,\ldots,m\} 
         \end{array}
\right.
\]
is positively emerging around $(x_1,\ldots,x_l)$ and negatively
emerging around $(y_1,\ldots,y_m)$; moreover, by the symmetry of $v$, 
$\beta_z(u)=0$, $\forall z\in\{x_1,\ldots,x_l,y_1,\ldots,y_m\}$.
The conditions $I'(u)[(u^+_i)^\delta]=0$ and $I'(u)[(u^-_j)^\delta=0]$ are
obviously satisfied by the definition of $\tilde{v}^+_{x_i}$ and
$\tilde{v}_{y_j}^-$.  
So relation $(i)$ of (\ref{*}) is shown true.
The same argument clearly works either when $l=0,$ either when $m=0,$ and allows
to obtain relations $(ii)$ and $(iii)$ of (\ref{*}).

\qed

\begin{lemma}
\label{L6.1}
Let assumptions of Lemma \ref{L5.2} hold, and  let
$ l,\ m,\ \e,\ (x_1,\ldots,x_l)$, $(y_1,$ $\ldots,y_m)$ be as in Lemma \ref{L5.2}. Then
\beq
\label{e6.1}
u\in S^\e_{x_1,\ldots,x_l,y_1,\ldots,y_m} \quad\Rightarrow\quad
\left\{
\begin{array}{llc}
a)&I_\e(u)>0,&\\
b)&J_\e((u_i^+)^\delta)>0 &1\le i\le l,\\
c)&J_\e((u_j^-)^\delta)>0 &1\le j\le m.
\end{array}
\right.
\eeq
\end{lemma}

\proof
Let us consider $ u\in S^\e_{x_1,\ldots,x_l,y_1,\ldots,y_m}$. 
Being, by (\ref{3.1}), $|u_\d (x)|\le \d<\min\{\bar a^{1\over p-1},1\}$ for all
$x\in \R^N$, we can deduce  
\begin{eqnarray}
&I_\e(u_\d) =\frac{1}{2}\int_{\R^N}(|\nabla u_\d|^2 +a_\e(x)u_\d^2)dx
- \frac{1}{p+1}\int_{\R^N} |u_\d|^{p+1} dx 
\geq \nonumber \\
&\hspace{0,5cm} \label{6.1} 
\geq \frac{1}{2}\int_{\R^N}|\nabla u_\d|^2 dx+ \bar a \left( \frac{1}{2} -
\frac{1}{p+1}\right)\int_{\R^N} u_\d^2 dx >0.\qquad
\end{eqnarray}
Hence using (\ref{punto})  and (\ref{6.1}) we obtain
\beq
\label{2.7 b}
I_\e(u) = \max_{t\ge 0} I_\e(u_\d + t (u^+)^\d- (u^-)^\d) > 
\max_{t\ge 0} I_\e(u_\d - t(u^-)^\d)>I_\e(u_\d) > 0.
\eeq

Moreover, considering (\ref{1836}), we get
\beq
\label{2.7 c}
J_\e((u_i^+)^\d) = \max_{t\ge 0} J_\e( t (u_i^+)^\d) > J_\e(0) = 0,
\eeq
$$
J_\e((u_j^-)^\d) = \max_{t\ge 0} J_\e( t (u_j^-)^\d) > J_\e(0) = 0.
$$

\qed

\begin{rem}{\em
\label{CC}
It is worth observing that relation (\ref{6.1}) can be proved without
any change for any function belonging to 
$
\cC_\delta=\{u\in H^1(\R^N)\ :\
|u|\le \delta\}.
$
So, for all $\e>0$, $I_\e$ is coercive on $C_\delta$.
Moreover condition $\delta<\left(\bar a\over p\right)^{1\over
  p-1}$ allows to state that on the same set $I_\e$ is also convex
and, hence, weakly lower semicontinuous.
}
\end{rem}

\begin{rem}
\label{Rinfty}
{\em We stress the fact that, when $a_\e(x)= a_\infty$ and we work
  with the functionals $I_\infty$ and  $J_\infty,$ everything we have
  done considering $I_\e$ and $J_\e$ can be repeated.  

The definition of the sets of multi-bump functions satisfying  local
natural constraints and  barycenter type constraints still makes
sense, but in this case we use the notation
$S^\infty_{x_1,\ldots,x_l,y_1,\ldots,y_m},\ \
S^{\infty,+}_{x_1,\ldots,x_l},$ and $ S^{\infty,-}_{y_1,\ldots,y_m}.$  

The conclusions of Lemmas \ref{L3.1},\ \ref{L5.2},\ \ref{L6.1},  and
 Corollary \ref{C5.1} (stated with $I_\infty$ and
$J_\infty$ instead of $I_\e$ and $J_\e$) hold true. With respect to
Lemmas  \ref{L3.1} and Corollary \ref{C5.1} we remark that, in this
case, for the numbers whose existence is there claimed the notation
used is, respectively, $\theta_\infty^+(u), \ \theta_\infty^-(u), \
\vt^{\infty,+}_s(u)$ and $\vt^{\infty,-}_s(u)$. 

}\end{rem}

\begin{lemma}
\label{L2.10}
The sets $\{\theta_\e^+(w_z)\ :\ z\in\R^N\}$ and $\{\theta_\e^-(-w_z)$
: $z\in\R^N\}$ are bounded, uniformly with respect to $\e$.
\end{lemma}

\proof
By construction, for all $z \in \R^N$, $\theta_\e^+ (w_z) > 0$.
On the other hand, by the choice of $\rho$, $\supp  (w^+_z)^\d \subset
B_\rho(z)$, so from (\ref{4.1}) we deduce 
\begin{eqnarray*}
0 & = & \theta_\e^+(w_z) \left[ \int_{B_\rho(z)} (|\nabla (w^+_z)^\d|^2 +
a_\e (x)((w^+_z)^\d)^2)dx \right] + \d
\int_{B_\rho(z)}a_\e(x) (w_z^+)^\d dx  
\\
& & -\int_{B_\rho(z)}(\d + \theta_\e^+ (w_z)(w_z^+)^\d)^p (w_z^+)^\d
dx 
\\
& \leq & \theta_\e^+(w_z)\left[  \|(w_z^+)^\d\|^2 + 
\int_{B_\rho(z)} |\a_\e(x) ((w_z^+)^\d)^2)|dx \right] 
+  \d \int_{B_\rho(z)}a_\e (x)(w_z^+)^\d dx 
\\
& & -[\theta_\e^+(w_z)]^p \int_{B_\rho(z)}((w_z^+)^\d)^{p+1} dx
\end{eqnarray*}
from which, when $\theta_\e^+(w_z) > 1$, we get

$$
\hspace{-5cm} [\theta_\e^+(w_z)]^{p-1} 
$$
$$\leq \frac{ \|(w^+)^\d\|^2 + C_\rho
  |\a_\e|_{\frac{N}{2}, \loc}|(w^+)^\d|^2_ {2^*}+ \d [a_\infty
  |(w^+)^\d|_1+ C_\rho |\a_\e |_{\frac{N}{2},
    \loc}|(w^+)^\d|_{N\over(N-2)}]}{|(w^+)^\d|_{p+1}^{p+1}},
$$
where $C_\rho$ denotes a constant depending on $\rho$. 
Then the statement follows because $\sup_{\e>0}|\a_\e
|_{\frac{N}{2},{\loc}}<\infty$. 

 Same computations can be done considering $\theta_\e^-(-w_z)$ instead of  $\theta_\e^+(w_z).$
 
\qed

 

 \sezione{Minimization among functions emerging around a fixed set of points}


Purpose of  this section is proving that for all $l,\ m \in \N$ such
that $l+m \geq 1,$ and for all fixed $(l+m)$-tuple
$(x_1,\ldots,x_l,y_1,\ldots,y_{m})$, the functional $I_\e$  achieves
its infimum on the nonsmooth constraint $S^{\e}_{x_1,\ldots,x_l,y_1,\ldots,y_{m}}.$
Afterwards, we study the asymptotic decay of the minimizers and we
show that the middle part of a minimizer is a solution of $(P_\e)$ on
$\R^N$ except the support of the emerging parts. Finally, considering
that a minimizer $\hat{u}$ in principle could not be a critical point
of $I_\e,$ we deduce the relation $I'_\e(\hat{u})$ must satisfy.

\begin{prop}
\label{P6.2}
Assume $(h_1),(h_2)$, and $(h_3)$ hold. 
Then for all $l,m\in\N\setminus\{0\}$,  for all
$\e\in(0,\e_l)$,  for all   $(x_1,\ldots,x_l)\in\cH_l^\e$ and
$(y_1,\ldots,y_{m})\in\cK_{m}^\e$,
$$
\exists\ \hat u_\e\in S^{\e}_{x_1,\ldots,x_l,y_1,\ldots,y_{m}}
\quad \mbox{ such that } \quad
I_\e(\hat u_\e)=\inf_{S^{\e}_{x_1,\ldots,x_l,y_1,\ldots,y_{m}} }I_\e>0,
$$
\beq
\label{*0}
\exists\ \bar u_\e\in S^{\e,+}_{x_1,\ldots,x_l}
\quad \mbox{ such that } \quad
I_\e(\bar u_\e)=\inf_{S^{\e,+}_{x_1,\ldots,x_l} }I_\e>0,
\eeq
\beq
\label{*1}
\exists\ \tilde u_\e\in S^{\e,-}_{y_1,\ldots,y_m}
\quad \mbox{ such that } \quad
I_\e(\tilde u_\e)=\inf_{S^{\e,-}_{y_1,\ldots,y_m} }I_\e>0.
\eeq
\end{prop}

\proof
Let $l,m\in\N\setminus\{0\}$, $(x_1,\ldots,x_l)\in\cH_l^\e$,
$(y_1,\ldots,y_m)\in\cK_m^\e$ be fixed and let $ (u_n)_n$ be a
minimizing sequence in $S^{\e}_{x_1,\ldots,x_l,y_1,\ldots,y_{m}}$.
The sequence $(I_\e(u_n))_n$ is bounded from above and, considering also 
(\ref{6.1}) and (\ref{2.7 b}), we deduce
\beq
\label{m.1}
c_1 \|(u_n)_\d\|^2 \leq I_\e((u_n)_\d) \leq I_\e(u_n) \leq c_2. 
\eeq
Moreover, for all $i\in \left\{1,\ldots,l\right\},$ and
$j\in\{1,\ldots,m\}$, the sequences
$\left(\|((u_n)_i^+)^\d\|/\right.$ $\left.|((u_n)^+_i)^\d|_{p+1} \right)_n $ and 
$\left(\|((u_n)_j^-)^\d\|/|((u_n)^-_j)^\d|_{p+1} \right)_n $ are
bounded. 
Indeed, assume that for some $i\in
\left\{1,\ldots,l\right\},$ up to a subsequence,
$\lim_{n\to
  \infty}{\|((u_n)_i^+)^\d\|/ |((u_n)^+_i)^\d|_{p+1}}  = +
\infty$.
Then we get
\begin{multline}
\label{m.2}
\quad \lim_{n\to \infty} J_\e
\left( 
{((u_n)_i^+)^\d\over |((u_n)^+_i)^\d|_{p+1}} 
\right)  \\
\geq \lim_{n\to \infty}
\left[c_3\ 
{\|((u_n)_i^+)^\d\|^2\over |((u_n)^+_i)^\d|^2_{p+1}} 
-\frac{1}{p+1}\int_{B_\rho (x_i)} \left(\d +
{((u_n)_i^+)^\d\over |((u_n)^+_i)^\d|_{p+1}} 
\right)^{p+1} dx \right] 
=+ \infty.
\end{multline}

On the other hand, in view of Corollary \ref{C5.1}, Remark \ref{R2.5}
and Lemma \ref{L6.1}, we have 
\begin{eqnarray*}
I_\e(u_n) & =& \max_{t>0} I_\e((u_n)_\delta+ \sum_{h\neq i}(u_n^+)^\delta_h +
t (u_n^+)^\delta_i-\sum_{j=1}^m((u_n)_j^-)^\delta)\\
 & = & I_\e((u_n)_\d) + \sum_{h\neq i} J_\e ((u_n^+)^\delta_h)
 + \max_{t>0} J_\e \left(t
\frac{((u_n)_i^+)^\delta}{|((u_n)_i^+)^\d|_{p+1}}\right) 
+\sum_{j=i}^mJ_\e((u_n^-)^\delta_j)
\\
& \geq & \ J_\e \left(
{((u_n)_i^+)^\d\over |((u_n)^+_i)^\d|_{p+1}} 
\right),
\end{eqnarray*}
that, if (\ref{m.2}) is true, contradicts (\ref{m.1}).
The same argument brings to a contradiction assuming 
$\lim_{n\to
  \infty}{\|((u_n)_j^-)^\d\|/ |((u_n)^-_j)^\d|_{p+1}}  = +
\infty$ for some $j\in\{1,\ldots,m\}$.

Now, setting $\check{u}_n = (u_n)_\d + 
\sum_{i=1}^l\frac{((u_n)_i^+)^\delta}{|((u_n)_i^+)^\d|_{p+1}}-
\sum_{j=1}^m\frac{((u_n)_j^-)^\delta}{|((u_n)_j^-)^\d|_{p+1}}$, 
we can assert that $\check{u} \in H^1(\R^N)$ exists so that, up to a subsequence,
\beq
\label{13}
\begin{array}{cl}
(a) & \check{u}_n \rightharpoonup \check{u} \hspace{2cm} \mbox{in}\
H^1(\R^N)\ \mbox{and}\ L^{p+1}(\R^N)\\ 
(b)& \check{u}_n \rightarrow \check{u} \hspace{2cm} \mbox{in} \
L^{p+1}_{loc}(\R^N)\\ 
(c)&\check{u}_n (x) \rightarrow \check{u}(x) \hspace{0,85 cm} a.e. \
\mbox{in}\ \R^N. 
\end{array}
\eeq

Therefore, for all $z\in\{x_1,\ldots,x_l,y_1,\ldots,y_m\}$, we have 
$\check{u}_n \tton \check{u}$ in $L^{p+1}(B_\rho(z))$
and so $|(\check{u}_i^+)^\d|_{p+1} = 1=|(\check{u}_j^-)^\d|_{p+1} $
$\forall i\in\{1,\ldots,l\}$, $\forall j\in\{1,\ldots,m\}$;
moreover $\beta_z (\check{u})=0$, by the continuity of the map $\beta_z$,
and, by (\ref{13})(c), $|\check{u}(x)| \leq \d$ in 
$\R^N \setminus (\cup_{i=1}^l B_\rho(x_i))\cup(\cup_{j=1}^mB_\rho(y_j))$. 
Furthermore, taking into account Remark \ref{CC}, the inequality
\beq
\label{num}
I_\e(\check u_\delta)\le \liminf_{n\to \infty}I_\e((u_n)_\delta)
\eeq
 follows.

Let us, finally, define 
$\hat{u} = \check{u}_\d + \sum_{i=1}^l
\theta_i^{\e,+}(\check{u})(\check{u}_i^+)^\d
-\sum_{j=1}^m\theta_j^{\e,-}(\check{u})(\check{u}_j^-)^\d$,
then $\hat{u} \in S^{\e}_{x_1,\ldots,x_l,y_1,\ldots,y_{m}}$ 
and it is the desired minimizer because, by the strong convergence of
the emerging parts and (\ref{num}),
\begin{eqnarray*}
0 < I_\e(\hat{u}) & =  &I_\e \left(\check{u}_\d +
\sum_{i=1}^l \theta_i^{\e,+}(\check{u})(\check{u}_i^+)^\d
-\sum_{j=1}^m\theta_j^{\e,-}(\check{u})(\check{u}_j^-)^\d\right)\\
&\le &
 \liminf_{n \to \infty}I_\e \left((u_n)_\d + 
\sum_{i=1}^l \theta_i^{\e,+}(\check{u})(\check{u}_n^+)_i^\d
-\sum_{j=1}^m\theta_j^{\e,-}(\check{u})(\check{u}_n^-)_j^\d\right) \\
&\le &
 \lim_{n \to \infty} I_\e(u_n) = 
\inf_{S^{\e}_{x_1,\ldots,x_l,y_1,\ldots,y_{m}}} I_\e.
\end{eqnarray*}

The same argument, setting either $m=0$ or $l=0$, proves (\ref{*0}) and
(\ref{*1}), respectively.

\qed

Let be $(x_1,\ldots,x_l)\in\cH^\e_l$ and $(y_1,\ldots,y_m)\in\cK_m^\e$,  define
$$
\mu^\e(x_1,\ldots,x_l,y_1,\ldots,y_m):=
\min_{S^{\e}_{x_1,\ldots,x_l,y_1,\ldots,y_{m}}
}I_\e,
$$
$$
\mu^{\e}(x_1,\ldots,x_l):=\min_{S^{\e,+}_{x_1,\ldots,x_l}}I_\e,
$$
$$
\mu^{\e}(y_1,\ldots,y_m):=\min_{S^{\e,-}_{y_1,\ldots,y_m}}I_\e,
$$
and  set
$$
M^\e_{x_1,\ldots,x_l,y_1,\ldots,y_m}:=
\{u\in S^{\e}_{x_1,\ldots,x_l,y_1,\ldots,y_{m}}\ :\ I_\e(u)=
\mu^\e(x_1,\ldots,x_l,y_1,\ldots,y_m)\},
$$
$$
M^{\e,+}_{x_1,\ldots,x_l}:=
\{u\in S^{\e,+}_{x_1,\ldots,x_l}\ :\ I_\e(u)=
\mu^{\e}(x_1,\ldots,x_l)\},
$$
$$
M^{\e,-}_{y_1,\ldots,y_m}:=
\{u\in S^{\e,-}_{y_1,\ldots,y_m}\ :\ I_\e(u)=
\mu^{\e}(y_1,\ldots,y_m)\}.
$$

\vspace{2ex}
Aim of next lemmas is to describe some important feature of the middle
part of the functions belonging to the sets
$M^\e_{x_1,\ldots,x_l,y_1,\ldots,y_{m}}$ and $M^{\e,-}_{y_1,\ldots,y_{m}}$. 

The first Lemma is concerned with the asymptotic decay of a function
belonging to $M^{\e,-}_{y_1,\ldots,y_{m}}$ and, moreover, states that
any such  function solves $(\tilde{P}_\e)$ in the whole space $\R^N$
except  on the support of the emerging parts.  
Its proof can be done arguing as in Lemma 3.4 of \cite{CPS1}. 

\begin{lemma}
\label{L1207}
Let assumptions $(h_1)$, $(h_2)$, $(h_3)$, and $(h_4)$ be satisfied,
let $\tilde u_\e$ be as in Proposition \ref{P6.2}.
Then $(\tilde u_\e)_\d$ satisfies 
$$
\left\{
\begin{array}{cl}
-\Delta u+a_\e(x)u=|u|^{p-1}u &\mbox{ in }\R^N\setminus\supp((\tilde
u^-_\e)^\d)
\\
u=-\delta &\mbox{ on } \supp(\tilde u^-_\e)^\delta.
\end{array}\right.
$$
Moreover, for all $\eta\in(\bar\eta,\sqrt{(a_\infty-\d^{p-1})}\,)$ a
constant $c_\eta$ depending  on $(a_\infty-\d^{p-1})$, $\eta$, $N$
exists such that  the function  $\tilde u_\d$ verifies
$$
-c_\eta\delta e^{-\eta d_{1,\e}(x)}<(\tilde u_\e)_\d(x)<0 \qquad\forall x\in\R^N
$$
where $d_{1,\e}(x)=\dist (x,\supp(\tilde u_\e)^\delta\cup\supp \alpha_\e^-)$.
\end{lemma}
To describe the behaviour of the middle part of a function belonging
to $M^\e_{x_1,\ldots,x_l,y_1,\ldots,y_{m}},$ we need to state
beforehand next Lemma, whose proof can be obtained combining
Proposition 3.1 and Lemma 3.4 of \cite{CPS1}. 

\begin{lemma}
\label{L1214}
Let assumptions $(h_1)$, $(h_2)$, $(h_3)$, and $(h_4)$ be satisfied,
let $\hat u_\e$ be as in Proposition \ref{P6.2}.
Then 
\beq
\label{P+}
\left\{
\begin{array}{cl}
-\Delta u+a_\e(x)u=u^p &\mbox{ in }\R^N\setminus\supp(\hat
u^+_\e)^\delta 
\\
u=\delta &\mbox{ on } \supp(\hat u^+_\e)^\delta
\\
0\le u\le\d & \mbox{ in }\R^N
\end{array}\right.
\eeq
has a unique solution $v$, moreover, for all
$\eta\in(\bar\eta,\sqrt{(a_\infty-\d^{p-1})}\,)$ a constant $c_\eta$
depending  on $(a_\infty-\d^{p-1})$, $\eta$, $N$ exists  such that the
function $v$ verifies 
$$
0<v(x)<c_\eta\delta e^{-\eta d_{1,\e}(x)}\qquad\forall x\in\R^N
$$
where $d_{1,\e}(x)=\dist (x,\supp(\hat u^+_\e)^\delta\cup\supp \alpha_\e^-)$.
\end{lemma}

We are now ready to prove:

\begin{lemma}
\label{L7.1}
Let $(h_1)$, $(h_2)$, $(h_3)$, and $(h_4)$ hold,
let $l,m\in\N\setminus\{0\}$,  
$\e\in(0,\e_l)$,    $(x_1,\ldots,x_l)\in\cH_l^\e,$ 
$(y_1,\ldots,y_{m})\in\cK_{m}^\e,$  and $\hat u_\e\in
M^\e_{x_1,\ldots,x_l,y_1,\ldots,y_m}$.  Then  $(\hat u_\e)_\delta$  solves 
 \beq
\label{1224}
\left\{
        \begin{array}{cl}
         -\Delta u+a_\e(x)u=|u|^{p-1}u          & \mbox{ in } 
\R^N \setminus (\supp (\hat {u}^+_\e)^\delta\cup\supp (\hat {u}^-_\e)^\delta)\\
        u = \delta    & \mbox{ on } \supp (\hat {u}^+_\e)^\delta\\
        u = -\delta    & \mbox{ on } \supp (\hat {u}^-_\e)^\delta. 
        \end{array}
\right.
\eeq
Moreover, setting $d_\e(x)=\dist(x,\supp (\hat
{u}^+_\e)^\delta\cup\supp (\hat {u}^-_\e)^\delta\cup\supp\,
\alpha_\e^-)$,  for all $\eta\in(\bar\eta,$ $\sqrt{(a_\infty-\d^{p-1})}\,)$ a
constant $c_\eta$ depending  on $(a_\infty-\d^{p-1})$, $\eta$, $N$
exists such that 
\beq
\label{e7.1}
|(\hat u_\e)_\delta|\le c_\eta \,\delta e^{-\eta d_\e(x)}\qquad\forall
x\in\R^N.
\eeq
\end{lemma}

\proof
In view of Remark \ref{CC}, the functional $I_\e$ is coercive and
convex on the convex set 
$$
\cL_\e=\{u\in H^1(\R^N)\ :\ |u|\le \delta,\ u=\delta\mbox{ on }
  \supp(\hat u_\e^+)^\delta,\  u=-\delta\mbox{ on }
  \supp(\hat u_\e^-)^\delta\}
$$
 and $(\hat u_\e)_\delta$ is the only minimizer
for the minimization problem 
$\min\{I_\e(u)\ :\ u\in\cL_\e\}.$
So, $(\hat u_\e)_\delta$ satisfies 
\beq
\label{grad}
I'_\e((\hat u_\e)_\d)[\nu - (\hat u_\e)_\delta]\geq 0 \qquad \forall
\nu \in \cL_\e. 
\eeq
Now, arguing as in \cite{Dumi,DXL,MaP} we infer from (\ref{grad}) that
\beq
\label{fi}
I'_\e((\hat u_\e)_\d)[\varphi]\ge 0\qquad
\forall \varphi \in C^\infty_0(\R^N\setminus (\supp
(\hat{u}^+_\e)^\delta\cup\supp (\hat {u}^-_\e)^\delta)),
\eeq
which means $(\hat u_\e)_\d$ is a solution of (\ref{1224}). 
In fact, for all $\vi$ and for all $t\in[-1,1]$ we obtain by direct
computation 
\begin{eqnarray*}
I_\e((\hat u_\e)_\d+t\vi)
&
\ge 
& I_\e([(\hat u_\e)_\d+t\vi]_\d)\\
&  & +
I'_\e([(\hat u_\e)_\d+t\vi]_\d)[(\hat u_\e)_\d+
t\vi-[(\hat
u_\e)_\d+t\vi]_\d]-c_\vi t^2,
\end{eqnarray*}
\begin{eqnarray*}
I_\e([(\hat u_\e)_\d+t\vi]_\d)
&
\ge &
I_\e((\hat u_\e)_\d)\\
& & +I'_\e((\hat u_\e)_\d)[[(\hat
u_\e)_\d+t\vi]_\d-(\hat u_\e)_\d ]-c_\vi t^2,
\end{eqnarray*}
where $c_\vi$ is a suitable positive constant and
$$
\hspace{-3cm} I'_\e([(\hat u_\e)_\d+t\vi]_\d)[(\hat u_\e)_\d+
t\vi-[(\hat
u_\e)_\d+t\vi]_\d]
$$
$$
=
\int_{\R^N}[a_\e(x)-|[(\hat u_\e)_\d+t\vi]_\d|^{p-1}][(\hat
u_\e)_\d+t\vi]_\d
[(\hat u_\e)_\d+
t\vi-[(\hat
u_\e)_\d+t\vi]_\d]dx\ge 0\qquad\forall t\in\R
$$
because of the choice of $\d$.
Thus, taking into account (\ref{grad}) we get 
$$
I_\e((\hat u_\e)_\d+t\vi )\ge I_\e((\hat u_\e)_\d )-2c_\vi t^2
$$
which implies (\ref{fi}) (actually, one can prove furthermore that  $-\d<(\hat
u_\e)_\d<\d$ in $\R^N\setminus (\supp
(\hat{u}^+_\e)^\delta\cup\supp (\hat {u}^-_\e)^\delta)$).

In order to prove (\ref{e7.1}), we first observe that  $v$ solution of
(\ref{P+}) and $(\hat u_\e)_\d^+$ are both solutions  of problem 
\[
 \left\{
        \begin{array}{cl}
         -\Delta u+a_\e(x)u=u^{p}          & \mbox{ in } 
\R^N \setminus (\supp (\hat {u}^+_\e)^\delta\cup \supp (\hat {u}_\e)^-_\delta)\\
        u = \delta    & \mbox{ on }  \partial \supp (\hat {u}^+_\e)^\delta\\
        0\le u\le\d& \mbox{ in } \R^N \setminus (\supp (\hat
        {u}^+_\e)^\delta\cup \supp (\hat {u}_\e)^-_\delta),
        \end{array}
\right.
\]
and $(\hat u_\e)_\d^+\le v$ on $\partial (\R^N \setminus (\supp
(\hat {u}^+_\e)^\delta\cup \supp (\hat {u}_\e)^-_\delta))$.
Hence:
$$
0\le (\hat u_\e)_\d^+\le v(x)< C\d e^{- \eta d_{1,\e}(x)} \quad
\forall x\in \supp (\hat u_\e)_\d^+\setminus \supp (\hat
{u}^+_\e)^\delta,
$$
that gives (\ref{e7.1}) for $(\hat u_\e)_\d^+$ because
$d_{1,\e}(x)\ge d_\e(x)$.

Since a quite analogous argument can be done for $(\hat u_\e)_\d^-$,
(\ref{e7.1}) follows. 

\qed

We close this section by a proposition that states the equations the gradient of $I_\e,$ constrained on $S^{\e}_{x_1,\ldots,x_l,y_1,\ldots,y_{m}}$ and evaluated at a minimizer, has to satisfy.

We remark that $S^{\e}_{x_1,\ldots,x_l,y_1,\ldots,y_{m}}$
is not a smooth manifold in $H^1(\R^N)$ because the conditions
$I'_\e(u)[(u^+)_i^\delta]=0$ and $I'_\e(u)[(u^-)_i^\delta]=0$ give
rise to  natural but nonsmooth constraints, nevertheless we can argue on the smooth constraints coming from the
condition on the barycenters.

\begin{prop}
\label{CPS3.6}
Let $(h_1)$, $(h_2)$, and $(h_3)$ hold, 
let be $l,m\in \N\setminus\{0\}$, $\e\in(0,\e_l)$, $(x_1,\ldots$ $,x_l)\in\cH^\e_l$,
 $(y_1,\ldots,y_m)\in\cK^\e_m$, and $\hat u\in
M^\e_{x_1,\ldots,x_l,y_1,\ldots,y_m}.$ 
Then, for all $z\in\{x_1,\ldots,x_l,$ $y_1,\ldots,y_{m}\},$ there exists $\lambda_z\in\R^N$ such that
\beq
\label{CPSe3.9}
I_\e'(\hat u)[\psi]=\int_{B_\rho(z)}(|\hat u|^+)^\delta
\psi(x)(\lambda_z\cdot (x-z))dx\qquad\forall\psi\in
H_0^1(B_\rho(z)).
\eeq
\end{prop}

\proof 
First of all we point out that relation (\ref{CPSe3.9}) means that,
for all $z\in\{x_1,\ldots,x_l,$ $y_1,\ldots,y_{m}\}, \ \hat u$  is a
critical point for $I_\e$ constrained on the set $\Theta_z=\{u\in 
H^1(\R^N)\ :\ u=\hat u$ in $\R^N\setminus B_\rho(z)$, $(|u|^+)^\d\neq
0$ in $B_\rho(z)$, $\beta_{z}(u)=0 \}.$ 

To prove  (\ref{CPSe3.9}) we argue by contradiction and we assume
that there exists $ 
z\in\{x_1,\ldots,x_l, y_1,\ldots,y_{m}\}$  such that
$\hat u$ is not a critical point for $I_\e$ on $\Theta_{z}$.
Let us fix $z=x_i\in\{x_1,\ldots,x_l\}$ and set $\Theta=\Theta_{x_i}$. 
 An application of standard deformation techniques yields the existence of 
$\hat r>0$ and of a continuous map $\zeta:[0,1]\times\Theta\to\Theta$ such that 
\beq
\label{1234'}
\begin{array}{lll}
a)&\zeta(0,u)=u & \forall u\in \Theta\\ 
b)&\zeta(\tau,u)=u  &  \forall\tau\in[0,1],\ 
\forall u\in  \Theta\setminus B_{\hat r}(\hat u)\\
c)&I_\e(\zeta(\tau,u))<I_\e(u) & \forall\tau\in(0,1],\ 
\forall u\in  B_{\hat r}(\hat u).
\end{array}
\eeq
We can also choose $\hat r$ so small that $(u^+)^\delta\not\equiv
0$ in $B_\rho(x_i)$,
$\forall u\in\bar B_{\hat r}(\hat u)$.
Then set
$$
t^-=\inf\{t<0\ :\ \hat u + t (\hat u_i^+)^\delta\in B_{\hat r}(\hat
u)\},
\quad 
t^+=\sup\{t>0\ :\ \hat u + t (\hat u_i^+)^\delta\in B_{\hat r}(\hat
u)\}.
$$
and remark that $t^->-1$.

Now, let us consider the continuous curve in $\Theta$ defined for
$t\in[t^-,t^+]$ by
$
t\mapsto \zeta(1,\hat u + t (\hat u_i^+)^\delta).
$
Using  (\ref{1234'}) $(b)$ we deduce that $\zeta(1,\hat u + t (\hat u_i^+)^\delta)=\hat u+t(\hat
u_i^+)^\delta$ for $t=t^+$ and $t=t^-$, so, in view of the fact that $0$ is the only critical point (a maximum point)
in $(-1,+\infty)$ of $I_\e(\hat u + t (\hat u_i^+)^\delta)$, we obtain
 
$$
I_\e'(\zeta(1,\hat u + t (\hat u_i^+)^\delta))[((\zeta(1,\hat u + t (\hat u_i^+)^\delta))_i^+)^\d]\ \left\{
\begin{array}{ll}
> 0          & \mbox{ if } t = t^-\\
< 0          & \mbox{ if } t = t^+ .
\end{array}
\right.
$$
 
As a consequence, we can assert that $\hat t\in (t^-,t^+)$ exists such that 
$$
I_\e'(\zeta(1,\hat u + \hat t (\hat
u_i^+)^\delta))[((\zeta(1,\hat u + \hat t (\hat u_i^+)^\delta))_i^+)^\d]=0. 
$$
Now, by construction, 
$$
\begin{array}{c}
((\zeta(1,\hat u + \hat t (\hat
u_i^+)^\delta))_i^+)^\d\neq 0\\
\zeta(1,\hat u + \hat t (\hat
u_i^+)^\delta)=\hat u\mbox{ in
}\R^N\setminus B_\rho(x_i), \qquad
\beta_{x_i}(\zeta(1,\hat u + \hat t (\hat
u_i^+)^\delta))=0;
\end{array}
$$
so $\zeta(1,\hat u + \hat t (\hat
u_i^+)^\delta)\in
S^{\e}_{x_1,\ldots,x_l,y_1,\ldots,y_{m}}$ and,
 by $(c)$ in (\ref{1234'}), 
$$
I_\e(\zeta(1,\hat u + \hat t (\hat
u_i^+)^\delta))<
I_\e(\hat u + \hat t (\hat u_i^+)^\delta)\le I_\e(\hat u )
$$
which contradicts the minimality of $\hat u$ in
$S^{\e}_{x_1,\ldots,x_l,y_1,\ldots,y_{m}}$. 

The same argument works when $z=y_j\in\{y_1,\ldots,y_m\}$ completing the proof.

\qed

\begin{cor}
\label{gvinc}
Let $(h_1)$, $(h_2)$, and $(h_3)$ hold, let be $l,m\in
\N\setminus\{0\}$, $(x_1,\ldots,x_l)\in\cH^\e_l$ for $\e\in(0,\e_l)$,
$(y_1,\ldots,y_m)\in\cK^\e_m$ for $\e>0$,  $\bar{u} \in
M^{\e,+}_{x_1,\ldots,x_l},$ and $\tilde{u} \in
M^{\e_,-}_{y_1,\ldots,y_m}.$  
Then, for all $z\in\{x_1,\ldots,x_l\}$  there exists
$\lambda_z\in\R^N$ such that 
$$
I_\e'(\bar{u} )[\psi]=\int_{B_\rho(z)}(|\bar{u} |^+)^\delta
\psi(x)(\lambda_z\cdot (x-z))dx\qquad\forall\psi\in
H_0^1(B_\rho(z)),
$$
and for all $z \in\{ y_1,\ldots,y_{m}\},$ there exists
$\lambda_z\in\R^N$ such that 
$$
I_\e'(\tilde{u} )[\psi]=\int_{B_\rho(z)}(|\tilde{u} |^+)^\delta
\psi(x)(\lambda_z\cdot (x-z))dx\qquad\forall\psi\in
H_0^1(B_\rho(z)).
$$

\end {cor}
 

 \sezione{ Maximization on the set of minimizers obtained on varying the (l+m)-tuples }


Aim of this section is to show that when  the $(l+m)$-tuples 
vary on $\cH_l^\e\times\cK_m^\e$ the supremum of the set of  minima
$\mu^\e(x_1,\ldots,x_l,y_1,\ldots,y_m)$   is, actually, a maximum. We
show also that an analogous result holds  for the sets consisting of
the minima $\mu^\e(x_1,\ldots,x_l),$ as $(x_1,\ldots,x_l)\in \cH_l,$
and $\mu^\e(y_1,\ldots,y_m)$  as $(y_1,\ldots,y_m) \in \cK_m^\e.$  
 
For all $l,m \in \N \setminus \left\{0\right\}$ we set
$$
\mu_{l,m}^\e=\sup_{\cH_l^\e\times
  \cK^\e_m}\min_{S^\e_{x_1,\ldots,x_l,y_1,\ldots,y_m}}I_\e, 
$$
moreover we put 
$$
\mu_{l,0}^\e=\sup_{\cH_l^\e}\min_{S^{\e,+}_{x_1,\ldots,x_l}}I_\e\ ; \ \ \ \ \ 
\mu_{0,m}^\e=\sup_{\cK_m^\e}\min_{S^{\e,-}_{y_1,\ldots,y_m}}I_\e.
$$

\begin{lemma}
\label{N16}
The relation
\beq
\label{16.1}
\mu^\infty :=\sup_{x\in\R^N}\min_{S_x^{\infty,+}}I_\infty=
\sup_{y\in\R^N}\min_{S_y^{\infty,-}}I_\infty= m_\infty 
\eeq
holds and, for all $ z \in \R^N,$
\beq
\label{16.2}
M^{\infty,+}_z:=\{u\in S^{\infty,+}_z\ :\ I_\infty(u)= m_\infty\} =
\left\{ w(\cdot -z)\right\}, 
\eeq
\beq
\label{meno}
M^{\infty,-}_z:=\{u\in S^{\infty,-}_z\ :\ I_\infty(u)= m_\infty\} =
\left\{- w(\cdot -z)\right\}. 
\eeq
\end{lemma}

\proof 
 First of all  we observe that, by the $\Z_2$-symmetry of $I_\infty$,
 the set $M^{\infty,-}_z$ consists just of functions $-u$ with $u \in
 M^{\infty,+}_z$.  We point out, also, that the same argument of
 Proposition \ref{P6.2} holds when one works with the functional
 $I_\infty$ on the set $S_z^{\infty,+},$  so we can assert that
 $M^{\infty,+}_z \neq \emptyset$, $M^{\infty,-}_z \neq \emptyset.$ 
 
 By the invariance of $\R^N$ and $I_\infty$ under the action
of the translations group, it is easily understood that
$\mu^\infty=\min_{S_z^{\infty,+}}I_\infty$ for all $z \in \R^N$.
Thus, to prove (\ref{16.1}), we just show that $m_\infty =
\mu^\infty(0):=\min_{{S}_0^{\infty,+}}I_\infty$.

Clearly, $w \in {S}_0^{\infty,+}$ hence
\beq
\label{16.3}
\mu^\infty(0) \leq I_\infty(w) = m_\infty.
\eeq
To show the reverse inequality let consider $u \in M_0^{\infty,+}$,
in view of Lemma \ref{L3.1} and Remark \ref{Rinfty}, $ 1$ is
the only maximum point in $(0,+\infty)$ of $\ \cI_u^\infty(t) :=
I_\infty (u_\d +t (u^+)^\d).$ Moreover, $\lim_{t \rightarrow +\infty}
\cI_u^\infty(t) = -\infty,$ so $\hat{t} > 1$ exists so that $\
\cI_u^\infty(\hat{t}) = I_\infty (u_\d +\hat{t} (u^+)^\d) < 0 .$ Now,
defining $\gamma : [0,1] \to H^1(\R^N)$ as 
\[
\gamma(t) =  \left\{
        \begin{array}{ll}
         2t u_\d         & t \in [0,1/2]\\
          u_\d +2(t - (1/2))\hat{t}(u^+)^\d & t\in[1/2,1],
        \end{array}
\right.
\]
$\gamma$ turns out to belong to the family $\Gamma$ defined in
(\ref{Gamma}) and, taking into account Remark \ref{CC} and
(\ref{2.1}), we obtain 
\beq
\label{m<mu}
\mu^\infty(0) = I_\infty(u)= \mbox{$\max$}_{t \in [0,1]}
I_\infty(\gamma(t))\geq \mbox{$\min$}_{\gamma \in
  \Gamma}\mbox{$\max$}_{t \in [0,1]}\ I_\infty(\gamma(t)) = m_\infty , 
\eeq
that, together with (\ref{16.3}), gives (\ref{16.1}). The relation
(\ref{16.2}) and, then (\ref{meno}), is  consequence of Lemma \ref{N2}. 

\qed

\begin{lemma}
\label{L1844}
Let $(h_1)$, $(h_2)$, $(h_3)$ hold,   let $l,m\in\N$ and $\e\in
(0,\e_l)$.  Then the
function $\mu^\e$ is upper semi-continuous on $\cH^\e_l$,  $\cK^\e_m$
and  $\cH^\e_l\times \cK^\e_m$.
\end{lemma}
\proof
We show
that $\mu^\e$ is upper semi-continuous on $\cH^\e_l;$ the proof of the upper semi-continuity on $\cK^\e_m$
and on $\cH^\e_l\times \cK^\e_m$ can be obtained exactly by the same argument. 

Let $(\bar x_1,\ldots,\bar x_l)\in\cH_l^\e$ be fixed, let be $\bar
u =\bar u_\d+\sum_{i=1}^l(\bar u_i^+)^\d \in M^{\e,+}_{ \bar x_1,\ldots,\bar x_l}.
$
and $((x_1^n,\ldots,$ $x_l^n))_n \in \cH^\e_l$ any sequence  such that
$(x_1^n,\ldots,x_l^n)\to (\bar x_1,\ldots,\bar x_l)$, as $n\to\infty$.
Let us consider 
$$
u_n=(u_n)_\d+\sum_{i=1}^l\theta_i^{\e,+}(u_n) (\bar u_n^+)_i^\d
$$
where
$$
(\bar u_n^+)_i^\d=(\bar u_i^+)^\d(x+\bar x_i-x_i^n)\qquad\mbox{ for
}i\in\{1,\ldots,l\} 
$$
and $(u_n)_\d$ is the unique function realizing 
$$
\min \{I_\e(u)\ :\ u\in H^1(\R^N),\ |u|\le \delta,\ u=\delta\mbox{ on
}\supp( u_n^+)^\delta\}.
$$
Then $u_n \in S^{\e,+}_{x_1,\ldots,x_l},$ $u_n \to\bar u$ in $H^1(\R^N)$ and we get
$$
\limsup_{n\to\infty}\mu^\e(x_1^n,\ldots,x_l^n)\le
\limsup_{n\to\infty}I_\e(u_n)=I_\e(\bar u)=
\mu^\e(\bar x_1,\ldots,\bar x_l),
$$
as desired.

\qed

\begin{lemma}
\label{L10.1}
Let $(h_1)$, $(h_2)$, $(h_3)$, and $(h_4)$ hold. Then, for all $l\in
\N\setminus\{0\}$ and for all $\e\in(0,\e_l),$ 
there exists $(\bar 
x_{1,\e},\ldots,\bar x_{l,\e})\in \cH^\e_l$ such that
\beq
\label{e10.1}
\mu^\e(\bar x_{1,\e},\ldots,\bar x_{l,\e})=\mu^\e_{l,0}.
\eeq
Moreover, for all neighborood $\tilde\cN$ of
$\partial A\cup\partial B$ such that $\cN_0\subset\subset\tilde\cN
\subset \cN$ and 
$\inf_{\tilde\cN}b>0$, the relation 
\beq
\label{9.59}
\lim_{\e\to 0}
\min_{i\in\{1,\ldots,l\}} \dist(\bar x_{i,\e},\tilde \cN/\e)=\infty
\eeq
holds.
\end{lemma}

\proof
 Relation (\ref{e10.1}) is a straight consequence of the compactness of $\cH^\e_l$ and of the
 upper semi-continuity of $\mu^\e$ on it.

\vspace{1mm}

The proof of (\ref{9.59}) is carried out in two steps: first we
consider $l=1$, then $l>1$.

{\bf Step 1: l=1.} We start showing that
\beq
\label{10.28}
\liminf_{\e\to 0} \mu_{1,0}^\e\ge m_\infty.
\eeq
Let us consider $\bar x\in B\setminus(A\cup \bar \cN)$ and functions $u_\e\in
M^{\e,+}_{\bar x/\e}$,  $\tilde
u_\e=(u_\e)_\delta+\theta_\infty^+(u_\e)(u_\e^+)^\delta$.
We have 
\beq
\label{4.13}
\mu_{1,0}^\e\ge I_\e(u_\e)\ge I_\e(\tilde u_\e)=I_\infty (\tilde u_\e)+
{1\over 2}\int_{\R^N}\alpha_\e(x)\tilde u_\e^2dx\ge m_\infty
+
{1\over 2}\int_{\R^N}\alpha_\e(x)\tilde u_\e^2dx,
\eeq
which gives (\ref{10.28}) if
\beq
\label{1851}
\left|\int_{\R^N}\alpha_\e(x)\tilde u_\e^2dx\right|=o(1)
\eeq
holds true.

Now, fixing arbitrarily  $\zeta>0$ real numbers $R>0$ and $\e>0$ can be found such that $|a(x)-a_\infty|<\zeta$
$\forall x\in\R^N\setminus B_R(0)$
and $\tilde u_\e^2<\zeta/|a-a_\infty|_{L^1(B_R(0))}$ in $B_R(0)$. 
Then, taking into account Lemma \ref{L7.1}, we get
\begin{eqnarray*}
\left|\int_{\R^N}\alpha_\e(x)\tilde u_\e^2dx\right|&\le&
\int_{\R^N}|a(x)-a_\infty|\tilde u_\e^2dx+
\int_{\R^N}|b(\e x)| \tilde u_\e^2dx\\
&\le & 
\int_{\R^N\setminus B_R(0)}|a(x)-a_\infty|\tilde u_\e^2dx+
\int_{B_R(0)}|a(x)-a_\infty|\tilde u_\e^2dx\\
& & +
\int_{\cN/\e}|b(\e x)| \tilde u_\e^2dx\\
&\le &
\zeta|\tilde u_\e|_2^2+\zeta+\max_{\cN}|b|\,{|\cN|\over
  \e^N} (c\,e^{-\tilde\eta \dist(\bar x,\cN)/\e}),
\end{eqnarray*}
where $\tilde\eta\in(\bar\eta,\sqrt {(a_\infty-\d^{p-1})})$.
Hence, (\ref{1851}) follows once $|\tilde u_\e|_2 < C$ for all $\e>0$
is proved. 
In order to do this
we set $\tilde w_{\bar x/\e}=(w_{\bar  x/\e})_\delta+\theta^{+}_\e(w_{\bar x/\e})(w_{\bar
  x/\e}^+)^\delta$ and, reminding that,  by Lemma \ref{L2.10},
$\theta^{+}_\e(w_{\bar x/\e})$ is bounded uniformly with respect to
$\e,$ we get 
$$
I_\e(u_\e)=\min_{S^{\e,+}_{\bar x/\e}}I_\e\le I_\e(\tilde w_{\bar
    x/\e})\le c\ I_\infty (w_{\bar x/\e})=c\ m_\infty.
$$
Moreover,
$$
\int_{\R^N}\alpha_\e^-(x)\tilde u_\e^2dx=\int_{B_{\dist(\bar
    x,\cN)/2\e}(\bar x/\e)}\alpha_\e^-(x)\tilde u_\e^2dx+
\int_{\R^N\setminus B_{\dist(\bar
    x,\cN)/2\e}(\bar x/\e)}\alpha_\e^-(x)\tilde u_\e^2dx
$$
which, by $(h_1)$ and the asymptotic decay of $u_\e$ and of
$\tilde{u}_\e$, yields  
$$
\liminf_{\e\to 0}\int_{\R^N}\alpha_\e(x)\tilde
u_\e^2dx\ge 0.
$$
Therefore, using (\ref{4.13}), (\ref{6.1}), and observing that, in
view of Remark \ref{Rinfty}, relations analogous to (\ref{e6.1})$(a)$
and (\ref{2.7 c}) hold, we infer    
\begin{eqnarray}
C & > & J_{\infty}((\tilde u_\e^+)^\delta)=\max_{t>0}J_{\infty}(t(\tilde u_\e^+)^\delta) 
\nonumber \\ & =& 
\max_{t>0}J_{\infty}\left(t\,{( \tilde u_\e^+)^\delta\over | ( \tilde u_\e^+)^\delta|_{p+1}}\right)\ge J_{\infty}\left({(\tilde  u_\e^+)^\delta\over | ( 
   \tilde u_\e^+)^\delta|_{p+1}}\right)
\nonumber\\ &\ge & 
c_1{\|( \tilde u_\e^+)^\delta\|^2\over | ( 
   \tilde u_\e^+)^\delta|^2_{p+1}}-c_2 \qquad \forall\e>0,
\qquad  \label{1309}
\end{eqnarray}
from which $\|(\tilde  u_\e^+)^\delta\|\le c |( \tilde
u_\e^+)^\delta|_{p+1}$ follows. 
Now, if  $\|(\tilde u_{\e}^+)^\delta\|$ where unbounded, we could
deduce the existence of a sequence of positive numbers $(\e_n)_n$,
$\e_n\to 0$ such that 
$$
J_{\infty}(( \tilde
u_{\e_n}^+)^\delta)\le\tilde c_1 \|( \tilde
u_{\e_n}^+)^\delta\|^2
-\tilde c_2  |( \tilde
u_{\e_n}^+)^\delta|^{p+1}_{p+1}+\tilde c_3\tton -\infty
$$
contradicting (\ref{e6.1})$(b)$. 
So, $\|(\tilde u_\e^+)^\delta\|\le c$ for all $\e>0$ and, using again
the asymptotic decay of $\tilde u_\e$ we get the boundedness of
$|\tilde u_\e|_2 $  and (\ref{1851}).

\vspace{1mm}

We are now in position to prove (\ref{9.59}). We argue by
contradiction and we assume the existence of a number 
$\bar r\ge0$ and of sequences $(\e_n)_n$, $(x_{\e_n})_n$ such that
$0<\e_n\tton 0,$ $x_{\e_n}\in\cH^{\e_n}_1$,
$\mu^{\e_n}(x_{\e_n})=\mu_{1,0}^{\e_n}$ and $\lim_{n\to\infty}
\dist(x_{\e_n},\tilde \cN/\e_n)=\bar r$. 
 
Setting $\tilde
w_{x_{\e_n}}=(w_{x_{\e_n}})_\delta+\theta_{\e_n}^+(w_{x_{\e_n}})(w_{x_{\e_n}}^+)^\delta$,
in view of Lemma \ref{N16}, we infer
$M^{\infty,+}_{x_{\e_n}}=\{w_{x_{\e_n}}\}$ so
$$
\mu_{1,0}^{\e_n}=\mu^{\e_n}(x_{\e_n}) \le I_{\e_n} (\tilde w_{x_{\e_n}})
=I_\infty(\tilde w_{x_{\e_n}})+ 
{1\over 2}\int_{\R^N}\alpha_{\e_n}(x)\tilde w_{x_{\e_n}}^2dx
$$
\beq
\label{11.58}
\le
I_\infty( w_{x_{\e_n}})+ 
{1\over 2}\int_{\R^N}\alpha_{\e_n}(x)\tilde w_{x_{\e_n}}^2dx
= m_\infty+
{1\over 2}\int_{\R^N}\alpha_{\e_n}(x)\tilde w_{x_{\e_n}}^2dx.
\eeq
Now, setting $\nu:=\inf_{\tilde\cN}b>0$, it is not difficult to see
that $\bar C:=\liminf_{n\to \infty}\meas\{x\in B_{\bar R}(x_{\e_n})$ :
$\alpha_{\e_n}^-(x)>\nu/2\}>0$,  where $\bar R=\max\{\rho,2\bar r\}$,
furthermore assumption $(h_1)$ implies  $\lim_{n\to
 \infty}\sup_{B_{ R}(x_{\e_n})}\alpha_{\e_n}^+=0$, $\forall R>0$, so 
\begin{eqnarray}
\int_{\R^N}\alpha_{\e_n}(x)\tilde w_{x_{\e_n}}^2dx & = & 
\int_{\R^N\setminus B_{\bar R}(x_{\e_n})}\alpha_{\e_n}(x)\tilde w_{x_{\e_n}}^2dx+ 
\int_{B_{\bar R}(x_{\e_n})}\alpha_{\e_n}(x)\tilde w_{x_{\e_n}}^2dx
\nonumber\\
& \le & \int_{{\R^N}\setminus B_{\bar R}(x_{\e_n})}\hspace{-2mm}
\alpha_{\e_n}^+(x) \tilde w_{x_{\e_n}}^2dx
+\int_{B_{\bar R}(x_{\e_n})} \hspace{-2mm}
(\alpha_{\e_n}^+(x)-\alpha_{\e_n}^-(x))\tilde
w_{x_{\e_n}}^2dx
\nonumber\\
& \le &  
-{\nu\over 2}\, \bar C \inf_{B_{\bar R}(0)}w^2+o(1)
\label{12.06}
\end{eqnarray}
that inserted in (\ref{11.58}) gives $\mu^{\e_n}_{1,0}\le
m_\infty-c+o(1)$ bringing to a contradiction with (\ref{10.28}).

{\bf Step 2: $\mathbf{l>1}$.} We claim that
\beq
\label{16.35}
\mu_{l,0}^\e\ge\mu_{l-1,0}^\e+m_\infty+o(1).
\eeq
Actually, once proved the above inequality, it is not difficult to
prove (\ref{9.59})  arguing by contradiction. 
Indeed, let us assume that a number $\bar r\ge 0$ and  sequences
$\e_n$ and $(\bar x_{1,\e_n},\ldots,\bar x_{l,\e_n})$ $\in\cH_l^{\e_n}$ exist so that $0<\e_n \to 0$,
$\mu^{\e_n}(\bar x_{1,\e_n},\ldots,\bar 
x_{l,\e_n})=\mu_{l,0}^{\e_n}$ and, for some $i\in\{1,\ldots,l\}$,
$\lim_{n\to \infty}\dist(\bar x_{i,\e_n},\tilde \cN/\e_n)=\bar r$. 
Assuming without any loss of generality $i=1$, we  choose
$z_{\e_n}\in  M^{\e_n,+}_{\bar x_{1,\e_n}}$, $s_{\e_n}\in M^{\e_n,+}_{\bar
  x_{2,\e_n},\ldots,\bar x_{l,\e_n}}$ and we observe that
$z_{\e_n}\vee s_{\e_n}\in S^{\e_n,+}_{\bar x_{1,\e_n},\ldots,\bar
x_{l,\e_n}}$, $I_{\e_n}(z_{\e_n}\wedge s_{\e_n})>0,$ because
$|z_{\e_n}\wedge s_{\e_n}|\le \delta$.
Moreover, the same argument displayed to prove (\ref{12.06}) gives
$$
I_{\e_n}(z_{\e_n})\le m_\infty-c+o(1).
$$
Therefore, we obtain
\begin{eqnarray*}
\mu_{l,0}^{\e_n} &\le & I_{\e_n}(z_{\e_n}\vee
s_{\e_n})=I_{\e_n}(z_{\e_n})+I_{\e_n}(s_{\e_n})-I_{\e_n}(z_{\e_n}\wedge s_{\e_n})
\nonumber \\
& \le &
I_{\e_n}(z_{\e_n})+I_{\e_n}(s_{\e_n})\le m_\infty-c+o(1)+\mu^{\e_n}(\bar
x_{2,\e_n},\ldots,\bar x_{l,\e_n})
\nonumber\\
& \le & 
\mu^{\e_n}_{l-1,0}+m_\infty-c+o(1),
\qquad 
\end{eqnarray*}
that contradicts (\ref{16.35}).

Let us now prove the claim.
Let $(\bar x_{1,\e},\ldots,\bar x_{l-1,\e})\in\cH^\e_{l-1}$ be such
that $\mu_\e(\bar x_{1,\e},\ldots,$ $\bar x_{l-1,\e})=\mu^\e_{l-1,0}$.
Passing eventually to a sequence, we can assume $\e\bar x_{i,\e}\to \xi_i\in\bar
B\setminus A$, for all $i\in\{1,\ldots,l-1\}$ and we consider
 $v_\e\in M^{\e,+}_{\bar x_{1,\e},\ldots,\bar x_{l-1,\e},\bar\xi/\e}$,
with  $\bar\xi\in B\setminus\bar\cN$  such that
 $\bar\xi\neq\xi_i$ $\forall i\in\{1,\ldots,l-1\}$. 
The function $v_\e$ can be written as
$$ 
v_\e = (v_\e)_\d + \sum_{i=1}^{l-1} ((v_\e)_i^+)^\d +
((v_\e)_{l}^+)^\d ,
$$
$((v_\e)_{l}^+)^\d $ denoting the part of $v_\e$ emerging around
$\bar\xi/\e$, and
$$
 \beta_{\bar x_{i,\e}}(v_\e) = 0 \ , \ \ 
\supp  ((v_\e)_i^+)^\d \subset B_\rho(\bar{x}_{i,\e}), \ \ \ \ 
\forall i \in \left\{ 1,\ldots,l-1\right\}; 
$$
$$
\beta_{\bar\xi/\e}(v_\e) = 0 \ , \ \supp((v_\e)_{l}^+)^\d   
\subset B_\rho(\bar\xi/\e).
$$
\no Set, now, $\hat v_\e(x)=\chi_\e(x)  v_\e(x)$, where $\chi_\e\in
C^\infty(\R^N,[0,1])$, $\chi_\e(x)= \chi(({\e/ R})\,$ $
  \left|x-{\bar\xi/\e}\right|)$, being  $R>0$
so that $B_R(\bar\xi)\subset B\setminus\bar\cN$ and $\chi\in
C^\infty([0,+\infty),[0,1])$  such
that $\chi(t)=1$ $\forall t\in [0,1/4]\cup[1,+\infty)$, $\chi(t)=0$
$\forall t\in[1/2,3/4]$.

Let us evaluate $I_\e(\hat v_\e)$. 
Taking into account the asymptotic behaviour stated
in (\ref{e7.1}) and setting $D_\e=B_{R/\e}(\bar\xi/\e)\setminus
B_{R/4\e}(\bar\xi/\e)$, we obtain 
\begin{eqnarray}
I_\e(\hat v_\e) & = & I_\e(\chi_\e v_\e)
\nonumber \\
 & \leq &
 I_\e(v_\e) + \frac{1}{2} \int_{D_\e}\left[ |\D\chi_\e|^2 
- \frac{1}{2} \Delta \chi_\e^2\right] v_\e^2 
+ \frac{1}{p+1} \int_{D_\e}(1- \chi_\e^{p+1})(v_\e)^{p+1} 
\nonumber \\ 
& \leq &
\mu^\e_{l,0} + o(1). \label{4.19}
\end{eqnarray}

On the other hand we can write $\hat v_\e(x)=\hat v_\e^I(x)+\hat
v_\e^{II}(x)$ with $\hat v_\e^I\in S^{\e,+}_{\bar \xi/\e}$, $v_\e^{II}\in
S^{\e,+}_{\bar x_{1,\e},\ldots,\bar x_{l-1,\e} }$, $\supp v_\e^I \cap \supp
v_\e^{II}=\emptyset$ and
\begin{eqnarray}
I_\e(\hat v_\e) & = & I_\e(\hat v_\e^I+\hat v_\e^{II})=
I_\e(\hat v_\e^I)+I_\e(\hat v_\e^{II})
\nonumber \\
&\ge  & \mu^\e(\bar\xi/\e)+\mu^\e(\bar x_{1,\e},\ldots,\bar x_{l-1,\e} )
\nonumber \\
& \ge & m_\infty+o(1)+\mu_{l-1,0}^\e
\label{16.30}
\end{eqnarray}
where the last inequality follows arguing as for proving (\ref{4.13})
and (\ref{1851}). 

Combining (\ref{16.30}) and (\ref{4.19}) we get (\ref{16.35}).

\qed
\begin{rem}
{\em It is worth pointing out that, arguing as before  for proving
  (\ref{16.35}), for all $m \in \N$ the relation  
\beq
\label{17.24}
\mu_{l,m}^{\e_n}\ge\mu_{l-1,m}^{\e_n}+m_\infty+o(1),
\eeq
 can be obtained.
}
\end{rem}

\begin{prop}
\label{P10.2}
Let $(h_1)$,--,$(h_4)$ hold, let be $l,m\in\N$, $m\ge 1$, and
$\e\in(0,\e_l)$, then  
\beq
\label{10.3}
\mu_{l,m}^\e>\mu_{l,m-1}^\e+ m_\infty.
\eeq
Moreover,   $(\tilde x_{1,\e},\ldots,\tilde
x_{l,\e})\in\cH_l^\e$, $(\tilde y_{1,\e},\ldots,\tilde
y_{m,\e})\in\cK^\e_m$ exist such that 
\beq
\label{10.4}
\mu_{l,m}^\e=\mu^\e(\tilde x_{1,\e},\ldots,\tilde
x_{l,\e},\tilde y_{1,\e},\ldots,\tilde
y_{m,\e}).
\eeq
\end{prop}

\proof
When $l=0$ the statement is nothing but  Proposition 4.4 in
\cite{CPS1}. 
So we assume $l\ge 1$.

The proof is carried out by an inductive argument on $m$.

{\bf Step 1: m=1.}\hspace{2mm}   
Let  $(\bar x_1,\ldots,\bar x_l)\in\cH_l^\e$ be such that $\mu_\e(\bar
x_1,\ldots,\bar x_l)=\mu_{l,0}^\e$ and let
$y_n=\sigma_n\tau$ be a sequence in $\R^N$, with $\tau\in\R^N$,
$|\tau|=1$, $\sigma_n\in \R$, and $\sigma_n\tton
+\infty$.
Let us define
$$
\Sigma_n=\left\{x\in\R^N\ :\ {\sigma_n\over2}-1<(x\cdot\tau)<{\sigma_n\over2}+1\right\}
$$
and consider  $u_n\in M^\e_{\bar x_1,\ldots,\bar x_l,y_n}$ where $n$ is so large
 that $y_n\in \cK^\e_1$.
We can write $u_n$ as
$$
u_n(x)=(u_n)_\delta(x)+\sum_{i=1}^l((u_n)_i^+)^\delta(x)-(u_n^-)^\delta(x)
$$
where $((u_n)_i^+)^\delta$ are the positively emerging parts
around $\bar x_i$, $i\in\{1,\ldots,l\}$, of $u_n$, and $(u_n^-)^\delta$ is the
negatively emerging part around $y_n$.
Clearly, for large $n$,
$$
\bigcup_{i=1}^l B_\rho(\bar x_i)\subset\{x\in\R^N\ :\
(x\cdot\tau)<\sigma_n/2-1\}
$$
$$
B_\rho(y_n)\subset\{x\in\R^N\ :\
(x\cdot\tau)>\sigma_n/2+1\}.
$$
Set now $v_n(x)=\tilde\chi_n(x) u_n(x)$, with
\beq
\label{chi}
\begin{array}{c}
\vspace{1mm}
\tilde\chi_n(x)
=\tilde\chi\left(\left|(x\cdot\tau)-\frac{\sigma_n}{2}\right|\right), 
\quad \mbox{ where } \tilde\chi\in C^\infty(\R^+,[0,1]) \mbox{ is such
  that }\\
\tilde\chi(t)=0\mbox{ if }t\le 1/2, \ \tilde\chi(t)=1\ \mbox{if}\ t\ge
1. 
\end{array}
\eeq
Arguing as in (\ref{4.19}), we infer 
\begin{eqnarray}
I_\e(v_n) & = &I_\e(\tilde\chi_n(x)u_n(x))
\nonumber \\ 
 & \le &
\mu_{l,1}^\e+c_1\int_{\Sigma_n}[(u_n)_\delta]^2dx+c_2\int_{\Sigma_n}[(u_n)_\delta]^{p+1}dx
\nonumber \\ & \le &
\mu^\e_{l,1}+o(e^{-\bar\eta\sigma_n}),
\label{12.1}
\end{eqnarray}
where last inequality follows from (\ref{e7.1}), taking into account
that, being $\supp \alpha_\e^-$ bounded,
$\dist(\supp \alpha_\e^-,\Sigma_n)\ge {\sigma_n\over 2}-c$.

On the other hand,
$$
v_n=v_n^I+v_n^{II}
$$
with
$v_n^I\in S^{\e,+}_{\bar x_1,\ldots, \bar x_l}$
defined by
$$
v_n^I(x)=\quad\left\{\begin{array}{lc}
0 &\mbox{if }(x\cdot\tau)\ge{\sigma_n\over 2}-{1\over 2}\\
\sum_{i=1}^l ((u_n)_i^+)^\delta (x)+\tilde\chi_n(x)(u_n)_\delta(x)
&\mbox{ if } (x\cdot\tau)<{\sigma_n\over 2}-{1\over 2}
\end{array}
\right.
$$
and $v_n^{II}\in S_{y_n}^{\e,-}$ defined by
$$
v_n^{II}(x)=\quad\left\{\begin{array}{lc}
0 &\mbox{if }(x\cdot\tau)\le {\sigma_n\over 2}+{1\over 2}
\\
-(u_n^-)^\delta (x)+\tilde\chi_n(x)(u_n)_\delta(x) 
&\mbox{ if } (x\cdot\tau)>{\sigma_n\over 2}+{1\over 2}.
\end{array}
\right.
$$
The choice of $(\bar x_1,\ldots,\bar x_l)$ and of $v_n^I$ implies
\beq
\label{12.3}
I_\e(v_n^I)\ge \mu^\e(\bar x_1,\ldots,\bar x_l)=\mu_{l,0}^\e.
\eeq
Considering $\tilde v_n\in S^{\infty,-}_{y_n}$ defined by
\beq
\label{10.44}
\tilde
v_n=(v_n^{II})_\delta-\theta^-_\infty(v_n^{II})(u_n^-)^\delta
\eeq
we deduce, using  (\ref{16.2}),
\begin{eqnarray*}
I_\e(\tilde v_n) & = & I_\infty (\tilde v_n)+{1\over
  2}\int_{\R^N}\alpha_\e(x)\tilde v_n^2dx
\\ & \ge &
m_\infty+ 
{1\over 2}\left[\int_{B_\rho(y_n)}\alpha_\e(x)\tilde v_n^2dx
-\sup_{\supp \alpha_\e^-}\tilde
v_n^2 \ \left(\int_{\supp\alpha_\e^-}|\alpha_\e(x)|dx\right)\right]
\end{eqnarray*}
and, since for large $n$ $\supp \alpha_\e^-\subset\{x\in\R^N\ :\ (x\cdot\nu)\le
\sigma_n/2-1\}$, we obtain 
\beq
I_\e(v^{II}_n)\ge I_\e(\tilde v_n)\ge m_\infty +{1\over
  2}\int_{B_\rho(y_n)}\alpha_\e(x)\tilde v_n^2dx.
\label{13.1}
\eeq
Therefore, combining (\ref{12.1}), (\ref{12.3}) and (\ref{13.1}), we infer
\begin{eqnarray*}
\mu_{l,1}^\e & \ge & I_\e(v_n)-o(e^{-\bar\eta\sigma_n})
\\ & \ge &\mu_{l,0}^\e+m_\infty+{1\over
  2}\int_{B_\rho(y_n)}\alpha_\e(x)\tilde
  v_n^2dx-o(e^{-\bar\eta\sigma_n}),
\end{eqnarray*}
that gives (\ref{10.3}), thanks to the assumption $(h_4)$ of
 slow asymptotic decay of $\alpha_\e$.

\vspace{2mm}

To  prove (\ref{10.4}), we just need to show that a sequence
$(x_1^n,\ldots,x_l^n,\bar y^n)$ such that 
$$
\lim_{n\to \infty}\mu^\e(x_1^n,\ldots,x_l^n,\bar y^n)=\mu^\e_{l,1}
$$
is bounded. Then (\ref{10.4}) follows by the upper semi-continuity of $\mu^\e$ on
$\cH^\e_l\times\cK^\e_1$.

If $(x_1^n,\ldots,x_l^n,\bar y^n)$ is not  bounded, then $|\bar y^n|\tton
\infty$ must be true. 
Setting
\beq
\label{11.22}
w_n(x)=(w_{\bar y^n})_\delta+\theta^+_\e(w_{\bar y^n}) (w_{\bar
  y^n}^+)^\delta,
\eeq
 $-w_n$ turns out to be negatively emerging around $\bar y^n$ and 
\begin{eqnarray*}
I_\e(-w_n)=I_\e(w_n)& = & I_\infty(w_n)
+{1\over  2}\int_{\R^N}\alpha_\e(x)w_n^2dx
\\ & \le &
I_\infty(w)
+{1\over  2}\int_{\R^N}|\alpha_\e(x)w_n^2|dx
\end{eqnarray*}
and, since $\alpha_\e\ttox 0$ and $|\bar y^n|\tton \infty$, writing 
$$
\int_{\R^N}|\alpha_\e(x)w_n^2(x)|dx=
\int_{\R^N\setminus B_{\tilde r}(0)}|\alpha_\e(x)w^2_n(x)|dx
+\int_{B_{\tilde r}(0)}|\alpha_\e(x)w^2_n(x)|dx
$$
for $\tilde r>0$, we easily obtain
$$
\lim_{n\to \infty} \int_{\R^N} |\alpha_\e(x)w^2_n(x)|dx=0.
$$
Therefore
\beq
\label{14.1}
I_\e(-w_n)\le m_\infty+o(1).
\eeq
Now, let $s_n=(s_n)_\delta+\sum_{i=1}^l(s_n^+)^\delta_i \in M^{\e,+}_{x_1^n,\ldots,x_l^n}$, 
and set
$
z_n=(z_n)_\delta+\sum_{i=1}^l(s_n^+)^\delta_i-(w_n^+)^\delta
$
with $(z_n)_\delta=\xi_n(x)(s_n)_\delta-\hat \xi_n(x)(w_{\bar
  y^n})_\delta$, where 
$\xi_n,\hat\xi_n\in C^\infty(\R^N,[0,1])$ are the cut-off functions
$$
\xi_n(x)=\tilde\chi\left({2\over |\bar y_n|}|x-\bar y_n|\right),
\qquad
\hat \xi_n(x)=
1-\tilde\chi\left({2\over |\bar y_n|}|x-\bar y_n|\right)
$$
and $\tilde{\chi}$ is defined as in (\ref{chi}).

A direct computation shows that
\begin{eqnarray}
I_\e(z_n) & = &
I_\e((z_n)_\delta)+\sum_{i=1}^lJ_\e(s_n^+)_i^\delta+J_\e((w_n^+)^\delta)
\nonumber \\ & = &
I_\e((z_n)_\delta)+I_\e(s_n)-I_\e((s_n)_\delta)+I_\e(-w_n)-I_\e((-w_n)_\delta)
\label{15.1}
\end{eqnarray}
and that, arguing as in (\ref{4.19}) and using (\ref{N2*}) and (\ref{e7.1}),
\beq
\label{15.2}
I_\e((z_n)_\delta)-I_\e((s_n)_\delta)-I_\e((-w_n)_\delta)=o(1).
\eeq
Thus, combining  (\ref{14.1}), (\ref{15.1}), (\ref{15.2}) and taking
into account $I_\e(s_n)\le\mu_{l,0}^\e$, we get
\beq
\label{15.3}
I_\e(z_n)\le \mu_{l,0}^\e+m_\infty+o(1).
\eeq
On the other hand, fixing  $u_n\in M^\e_{x_1^n,\ldots, x_l^n,\bar y^n}$
and considering $z_n\in S^\e_{x_1^n,\ldots,x_l^\e,\bar y^n}$, we have
\beq
\label{15.4}
\mu^\e(x_1^n,\ldots,x_l^\e,\bar y^n)=I_\e(u_n)\le I_\e(z_n).
\eeq
Hence, combining (\ref{15.3}) and (\ref{15.4}) and letting
$n\lo\infty$, we obtain
$$
\mu^\e_{l,1}\le\mu^\e_{l}+m_\infty
$$
that contradicts (\ref{10.3}) and allows us to conclude that 
$((x_1^n,\ldots,x_l^n,\bar y^n))_n$ is bounded.

\vspace{2mm}

{\bf Step 2:  $\mathbf{m>1}$.}\hspace{2mm}  
Assume now that (\ref{10.3}) and (\ref{10.4}) hold true for $m-1$
and let us prove them for $m$.

Let $(\tilde x_{1,\e},\ldots,\tilde x_{l,\e},\tilde
y_{1,\e},\ldots,\tilde y_{m-1,\e})$ be such that
\beq
\label{17.3}
\mu^{\e}_{l,m-1}=\mu^\e(\tilde x_{1,\e},\ldots,\tilde x_{l,\e},\tilde
y_{1,\e},\ldots,\tilde y_{m-1,\e})
\eeq
and let $y_n=\sigma_n\tau$ and $\Sigma_n$  be defined as in Step 1.
Let us consider a function
$u_n\in M^\e_{\tilde x_{1,\e},\ldots,\tilde x_{l,\e},\tilde
y_{1,\e},\ldots,\tilde y_{m-1,\e},y_n}$
and observe that it can be written as
$$
u_n=(u_n)_\delta+\sum_{i=1}^l(u_n^+)_i^\delta-\sum_{j=1}^{m-1}
(u_n^-)_j^\delta-(\hat u_n^-)^\delta
$$
where $(\hat u^-_n)^\delta$ is the part of $u_n$ negatively emerging around
$y_n$.
Notice that for large $n$ we have
$$
\left(\cup_{i=1}^l B_\rho(\tilde x_{i,\e})\right)\bigcup\left(
\cup_{j=1}^{m-1} B_\rho(\tilde y_{j,\e})\right)\subset\{x\in\R^N\ :\
(x\cdot\tau)<(\sigma_n/2)-1\}
$$
$$
B_\rho(y_n)\subset\{x\in\R^N\ :\ (x\cdot\tau)>(\sigma_n/2)+1\}.
$$
Now we argue as for proving (\ref{12.1}): we set $v_n(x)=\tilde
\chi_n(x)u_n(x),$  $\tilde{ \chi}_n$ being as in (\ref{chi}), and 
we  show that
\beq
\label{18.1}
I_\e(v_n)\le \mu^\e_{l,m}+o(e^{-\bar \eta\sigma_n}).
\eeq
On the other hand, we write $v_n=v_n^I+v_n^{II}$ with
\beq
\label{18.2}
v_n^I\in S^\e_{\tilde x_{1,\e},\ldots,\tilde x_{l,\e},\tilde
y_{1,\e},\ldots,\tilde y_{m-1,\e}},\qquad v^{II}_n\in S^{\e,-}_{y_n}
\eeq
defined by
$$
v_n^I(x)=\left\{
\begin{array}{cc}
0 & \mbox{if }(x\cdot\tau)\ge {\sigma_n\over 2}-{1\over 2}\\
\tilde \chi_n(x)(u_n)_\delta(x)+\sum_{i=1}^l(u^+_n)^\delta_i(x)-\sum_{j=1}^{m-1}(u_n^-)^\delta_j(x)&
\mbox{if }(x\cdot\tau) < {\sigma_n\over 2}-{1\over 2},
\end{array}
\right.
$$
$$
v_n^{II}(x)=\left\{
\begin{array}{cc}
0 & \mbox{if }(x\cdot\tau)\le {\sigma_n\over 2}+{1\over 2}\\
-(\hat u_n^-)^\delta(x)+\tilde \chi_n(x)(u_n)_\delta(x) &
\mbox{if }(x\cdot\tau)> {\sigma_n\over2}+{1\over2}.
\end{array}
\right.
$$
Now, by (\ref{18.2}) and (\ref{17.3})
\beq
\label{19.1}
I_\e(v_n^I)\ge \mu^\e(\tilde x_{1,\e},\ldots,\tilde x_{l,\e},\tilde
y_{1,\e},\ldots,\tilde y_{m-1,\e})=\mu^\e_{l,m-1}
\eeq
follows, moreover defining $\tilde v_n$ as in
(\ref{10.44}), working as for proving  (\ref{13.1}), we deduce
\beq
\label{11.05}
I_\e(v_n^{II})\ge m_\infty+{1\over
  2}\int_{B_\rho(y_n)}\alpha_\e(x)(\tilde v_n(x))^2dx.
\eeq
Hence, combining (\ref{18.1}), (\ref{19.1}) and (\ref{11.05}) we get
\begin{eqnarray*}
\mu_{l,m}^\e  & \ge &
I_\e(v_n)-o(e^{-\bar\eta\sigma_n}) \\
& \ge &
\mu^\e_{l,m-1}+m_\infty+ {1\over
  2}\int_{B_\rho(y_n)}\alpha_\e(x)(\tilde
v_n(x))^2dx-o(e^{-\bar\eta\sigma_n}) ,
\end{eqnarray*}
that, in view of $(h_4)$, implies (\ref{10.3}).

\vspace{1mm}

In order to prove (\ref{10.4}), we consider a sequence
$(x_1^n,\ldots,x_l^n,y_1^n,\ldots,y_m^n)$, such that
\beq
\label{19.2}
\lim_{n\to \infty}\mu^\e(
x_1^n,\ldots,x_l^n,y_1^n,\ldots,y_m^n)=\mu^\e_{l,m}
\eeq
and we show it is bounded, this done then (\ref{10.4}) follows by the
upper semi-continuity of $\mu^\e$. 
Arguing by contradiction, we can assume without any loss of generality
that, up to a subsequence, there exists $t \in\{1,\ldots,m\}$ such
that $|y^n_i|\tton \infty$ for $i\in\{t ,\ldots,m\}$ and
$(y_i^n)_n$ is bounded otherwise.
For all $i\in\{t ,\ldots,m\}$ we consider  $w^i_n$ defined as in
(\ref{11.22}), with $y^n_i$ instead of $\bar y^n$, so,
 arguing as for proving (\ref{14.1}), we obtain
\beq
\label{1624}
I_\e(-w^i_n)=I_\e(w_n^i)\le m_\infty+o(1).
\eeq
Now,  let us choose $s_n\in
M^\e_{x_1^n,\ldots,x_l^n,y_1^n,\ldots,y_{t -1}^n}$, $s_n\in
M^{\e,+}_{x_1^n,\ldots,x_l^n}$ if $t =1$, be: 
$$
s_n=(s_n)_\delta
+\sum_{i=1}^l(s_n^+)_i^\delta-\sum_{j=1}^{t -1}(s_n^-)^\delta_j,
$$
and let us consider
$$
z_n=(z_n)_\delta+\sum_{i=1}^l(s_n^+)_i^\delta
-\sum_{j=1}^{t -1}(s_n^-)^\delta_j
-((w_n^t \vee\ldots\vee w_n^m)^+)^\delta
$$
where 
$$
(z_n)_\delta(x)=\left[1-\tilde \chi\left({2\over d_n}\,
    |x|\right)\right](s_n(x))_\delta
-\tilde \chi\left({2\over d_n}\, |x|\right) 
((w_n^t \vee\ldots\vee w_n^m)(x))_\delta,
$$
with $\tilde{\chi}$ defined as in (\ref{chi}) and  $d_n=\min\{|y_i^n|$
: $i\in\{t ,\ldots,m\}\}$. 

Since $ |w^{i_1}_n\wedge w_n^{i_2}|<\d,$  
considering Remark \ref{CC} and  using (\ref{1624}), we get
\begin{eqnarray*}
I_\e(-(w_n^t \vee\ldots\vee w_n^m)) &  = & I_\e(w_n^t )+
I_\e(w_n^{t +1}\vee\ldots\vee w_n^m)
-I_\e(w^t \wedge(w_n^{t +1}\vee\ldots\vee w_n^m ))
\\
& < &  I_\e(w_n^t )+ I_\e(w_n^{t +1}\vee\ldots\vee w_n^m).
\end{eqnarray*}
Now, arguing analogously, we have
$$
 I_\e(w_n^{t +1}\vee\ldots\vee w_n^m) <
I_\e(w_n^{t+1} )+I_\e(w_n^{t +2}\vee\ldots\vee w_n^m);
$$
so, repeating $m-t+1$ times the same kind of computation, we get
$$
I_\e(-(w_n^t \vee\ldots\vee w_n^m)) < \sum_{i=t}^m I_\e(w_n^{i} ) 
\le (m-t +1)\ m_\infty+o(1).
$$
Therefore, arguments analogous to those developed in (\ref{15.1}),--,
(\ref{15.3}) bring to 
\begin{eqnarray}
I_\e(z_n)& = &
I_\e(s_n)+I_\e(-(w_n^{t }\vee\ldots\vee
w_n^m))  
\nonumber \\
&  &
+I_\e((z_n)_\delta)-I_\e((s_n)_\delta)-I_\e(-(w_n^{t }\vee\ldots\vee
w_n^m)_\delta) 
\nonumber \\
& \le &
\mu_{l,t -1}^\e+(m-t +1)m_\infty+o(1).
\label{20.1}
\end{eqnarray}
On the other hand, being  $z_n\in S^\e_{x_1^n,\ldots,x_l^n,y_1^n,\ldots,y_{m}^n}$
$$
\mu^\e({x_1^n,\ldots,x_l^n,y_1^n,\ldots,y_{m}^n})\le
I_\e(z_n),
$$
holds true, so, letting $n\to\infty$ in the above inequality, and using (\ref{19.2}) and
(\ref{20.1}), we obtain
$$
\mu^\e_{l,m}\le \mu^\e_{l,t -1}+(m-t +1 )m_\infty,
$$
which contradicts  (\ref{10.3}) completing the proof. 

\qed



 \sezione{Behaviour of  max-min functions as $\e\to 0$}


The functions realizing the max min values of $I_\e$ are expected to
be solutions of $(P_\e);$ a basic ingredient for the attainment of
this conclusion turns out to be the study, to which this section is
devoted, of some asymptotic properties of these functions  as $\e\to
0$.  
The results of Proposition \ref{P22.1} and \ref{CPS5.1}
respectively, allow to assert  that the distances between the points
around which these functions emerge and the boundary of $A$ and $B$
increase  going to infinity, as $\e\to 0,$ and that, also, the
interdistances between the same points become larger and larger as
$\e\to 0$.  
In Proposition \ref{CPS5.4} we analyze the asymptotic shape of the
emerging parts of these ``candidate solutions'' and we conclude they
approach more and more, as $\e\to 0,$ the ground state solution of
$(P_\infty)$.    

\vspace{2ex}
In what follows  we denote by {\bf $u_{l,m,\e}$ a function in
  $M^\e_{\tilde x_{1,\e},\ldots,\tilde x_{l,\e},\tilde
    y_{1,\e},\ldots,\tilde y_{m,\e}}$}.

\begin{prop}
\label{P22.1}
Let $(h_1)$,--,$(h_4)$ hold and let $\tilde\cN$ be as in Lemma \ref{L10.1}. Let be $l\in\N, \ m\in\N\setminus\{0\}, \ \e\in(0,\e_l),$ and 
let   $(\tilde x_{1,\e},\ldots,\tilde x_{l,\e},\tilde
y_{1,\e},\ldots,\tilde  y_{m,\e})$  be  a (l+m)-tuple for which $\mu_{l,m}^\e=\mu^\e(\tilde x_{1,\e},\ldots,\tilde
x_{l,\e},\tilde y_{1,\e},\ldots,\tilde
y_{m,\e})$ holds. Then, for all $r>0$ an
 $\e_r\in(0,\e_l)$ exists such that for all $\e\in(0,\e_r)$ and for
 all $m\in\N\setminus\{0\}$ 
\beq
\label{e22.2}
 \min\{\dist(\tilde y_{i,\e}, \tilde\cN /\e)\ :\ 
i\in\{1,\ldots,m\}\}>r,
\eeq
\beq
\label{e22.1}
\min\{ \dist(\tilde x_{i,\e},  \tilde \cN/\e)\ :\ 
{i\in\{1,\ldots,l\}}\}>r.
\eeq
\end{prop}

\proof
We argue by contradiction and we assume that there exist $\bar r\ge0,$ and
sequences $(\e_n)_n$, $0<\e_n\tton 0$, $(m_n)_n$ in
$\N\setminus\{0\}$, $(\tilde
x_{1,\e_n},\ldots,\tilde x_{l,\e_n},\tilde 
y_{1,\e_n},\ldots,\tilde
y_{m_n,\e_n})\in\cH^{\e_n}_l\times\cK^{\e_n}_{m_n}$, with $\mu^{\e_n}
(\tilde x_{1,\e_n},\ldots,$ $\tilde x_{l,\e_n},\tilde
y_{1,\e_n},\ldots,\tilde
y_{m_n,\e_n})=\mu^{\e_n}_{l,m_n}$, such that either (\ref{e22.2}), or
(\ref{e22.1}), does not hold.

 Without any loss of generality we can suppose that the relation not
 true is (\ref{e22.2}), and that
$$
\dist(\tilde y_{1,\e_n}, \tilde \cN/\e)\le\bar r\qquad\forall
n\in\N.
$$
Now, let us choose $z_n\in M^{\e_n,-}_{\tilde y_{1,\e_n}}$, $s_n\in
M^{\e_n}_{\tilde x_{1,\e_n},\ldots,\tilde x_{l,\e_n},\tilde
  y_{2,\e_n},\ldots,\tilde y_{m_n,\e_n}}$, a neighborhood $\hat\cN$ of
$\partial B$, with $\hat\cN\subset\subset\cN_0$, and a cut-off function 
$$
\hat\chi\in C^\infty(\R^N,[0,1])\quad\mbox{ such that }\
\hat\chi_{|_{B\setminus\hat\cN}}\equiv 1,\quad 
\hat\chi_{|_{\R^N\setminus(B\cup \hat\cN)}}\equiv 0.
$$
and define the function  $v_n\in S^{\e_n}_{\tilde x_{1,\e_n},\ldots,\tilde x_{l,\e_n},\tilde
  y_{1,\e_n},\ldots,\tilde y_{m_n,\e_n}}$as $v_n=\hat\chi(\e_nx)s_n+(1-\hat\chi(\e_nx))(z_n\wedge
s_n)$.
 The arguments used in (\ref{11.58}) and
(\ref{12.06}) show that $I_{\e_n}(z_n)\le m_\infty-c+o(1)$. 
Hence, taking into account   Lemma
\ref{L7.1}
and the coercivity of $I_{\e_n}$ on $C_\delta$ ,
we get
$$
 \mu^{\e_n}_{l,m_n}\le I_{\e_n}(v_n)
$$
$$
=
\int_{(B\setminus \hat\cN)/\e_n}\left[{1\over 2}(|\D
  s_n|^2+a_{\e_n}(x)s_n^2)-
{1\over p+1}|s_n|^{p+1}\right]dx
$$
$$
+\int_{\hat\cN/\e_n}\left[{1\over 2}(|\D
  v_n|^2+a_{\e_n}(x)v_n^2)-
{1\over p+1}|v_n|^{p+1}\right]dx
$$
$$
+\int_{\R^N\setminus (B\cup \hat\cN)/\e_n}\left[{1\over 2}(|\D
 (z_n\wedge s_n)|^2+a_{\e_n}(x)(z_n\wedge s_n)^2)-
{1\over p+1}|z_n\wedge s_n|^{p+1}\right]dx
$$
$$
=I_{\e_n}(s_n)+I_{\e_n}(z_n)-I_{\e_n}((1-\hat\chi({\e_n}x))(s_n\vee
z_n))+o(1)
$$
\beq
\label{1710}
\le \mu^{\e_n}_{l,m_n-1}+m_\infty-c+o(1),
\eeq
which contradicts (\ref{10.3}).

If the relation not true were (\ref{e22.1}), we could argue in a similar way, by
using (\ref{17.24})

\qed

\begin{rem}
\label{og}
{\em
Since the conclusion of Proposition \ref{P22.1} holds when $\tilde\cN\subset\cN$ is any neighborood of $\partial A\cup\partial B$,  such that
$\inf_{\tilde\cN} b>0$, and $(h_1)$ holds, we can assert that
 $\forall R>0$
$$
\lim_{\e\to 0}\sup_{B_R(\tilde x_{i,\e})}|\alpha_\e|=
 \lim_{\e\to 0}\sup_{B_R(\tilde y_{j,\e})}|\alpha_\e|=0\quad
\forall i\in\{1,\ldots,l\},\ \forall
j\in\{1,\ldots,m\},
$$
and that the above limits are uniform with respect to $m\in\N$.
}
\end{rem}

\begin{prop}
\label{CPS5.1}

Let assumptions of Proposition \ref{P22.1} hold. Let $l,\ m, \ \e_r, \ (\tilde x_{1,\e},\ldots,\tilde x_{l,\e},$ $\tilde
y_{1,\e},\ldots,\tilde  y_{m,\e})$ be as in Proposition \ref{P22.1}.
Then, for all $r>0,$  $\bar\e_r\in(0,\e_r)$  exists such that for
all $\e\in(0,\bar\e_r)$ and for all $m\in\N\setminus\{0\}$
\beq
\label{e23.2}
\min\{ |\tilde y_{i,\e}-\tilde y_{j,\e}|\  :\
i,j\in\{1,\ldots,m\},\ i\neq j\}>r.
\eeq
Moreover, when $l\ge 1, $ relation 
\beq
\label{e23.5}
 \min\{ |\tilde x_{i,\e}-\tilde y_{j,\e}|\  :\
i\in\{1,\ldots,l\},\ j\in\{1,\ldots,m\}\}>r,
\eeq
holds, and, when $l\ge 2,$ also
\beq
\label{e23.1}
\min\{ |\tilde x_{i,\e}-\tilde x_{j,\e}|\  :\
i,j\in\{1,\ldots,l\},\ i\neq j\}>r
\eeq
is true.
\end{prop}

\proof
  We prove (\ref{e23.2}) by contradiction, hence we assume it false.
Then, a real number $\bar r\ge 3\rho$ and sequences $(\e_n)_n,\ \e_n \in \R,$
$(m_n)_n$ in $\N$, $(\tilde x_{1,\e_n},\ldots,\tilde x_{l,\e_n},\tilde
y_{1,\e_n},\ldots,\tilde
y_{m_n,\e_n})\in\cH^{\e_n}_l\times\cK^{\e_n}_{m_n}$ exist  such that  $0<\e_n\tton 0$, $\mu^{\e_n}
(\tilde x_{1,\e_n},\ldots,$ $\tilde x_{l,\e_n},\tilde
y_{1,\e_n},\ldots,\tilde
y_{m_n,\e_n})=\mu^{\e_n}_{l,m_n}$ and
$$
\min\{|\tilde y_{i,\e_n}-\tilde y_{j,\e_n}|\ :\ 
i,j\in\{1,\ldots,m_n\},\ i\neq j\}\le\bar r \qquad\forall n\in\N.
$$
Without any loss of generality we can assume that
\beq
\label{1850}
|\tilde y_{1,\e_n}-\tilde y_{2,\e_n}|\le \bar r\qquad\forall n\in\N.
\eeq
We claim that from (\ref{1850}) relation
$$
\hspace{-3cm} \limsup_{n\to\infty}|\mu^{\e_n}(\tilde x_{1,\e_n},
\ldots,\tilde x_{l,\e_n},\tilde y_{1,\e_n},\tilde
y_{2,\e_n}\,\ldots,\tilde y_{m_n,\e_n}) 
$$
\beq
\label{1857}
\hspace{3cm} -
\mu^{\e_n}(\tilde x_{1,\e_n},
\ldots,\tilde x_{l,\e_n},\tilde y_{2,\e_n},\ldots,\tilde
y_{m_n,\e_n})|<m_\infty
\eeq
follows.
Actually, once the claim is proved, we are done, because
$\mu^{\e_n}(\tilde x_{1,\e_n},
\ldots,\tilde x_{l,\e_n},$ $\tilde y_{1,\e_n},\tilde y_{2,\e_n},\ldots,\tilde
y_{m_n,\e_n})=\mu^{\e_n}_{l,m_n}$ and $\mu^{\e_n}(\tilde x_{1,\e_n},
\ldots,\tilde x_{l,\e_n},\tilde y_{2,\e_n}\ldots,\tilde
y_{m_n,\e_n})\le\mu^{\e_n}_{l,m_n-1}$, so, for large $n$, (\ref{1857})
contradicts (\ref{10.3}).

Let us consider the functions $\bar v_n(x)=\hat\chi(\e_nx)\bar u_n
(x)+(1-\hat\chi(\e_nx))(\bar u_n\wedge \tilde w_{\tilde y_{1,\e_n}})(x)$,
where $\bar u_n\in M^{\e_n}(\tilde x_{1,\e_n},
\ldots,\tilde x_{l,\e_n},\tilde y_{2,\e_n},\ldots,\tilde
y_{m_n,\e_n})$, $\tilde w_{\tilde y_{1,\e_n}}=(-w_{\tilde
  y_{1,\e_n}})_\delta-\theta^-_{\e_n}(-w_{\tilde y_{1,\e_n}})$ $(w_{\tilde
  y_{1,\e_n}}^+)^\delta\in S^{\e_n,-}_{\tilde y_{1,\e_n}},$ and the cut off
 $\hat\chi, $ and the related set $\hat\cN$ , are defined as in Proposition \ref{P22.1}. 

Since $\bar v_n\in S^{\e_n}_{\tilde x_{1,\e_n},\ldots,\tilde
  x_{l,\e_n},\tilde y_{1,\e_n}\ldots,\tilde y_{m_n,\e_n} }$, arguing
as in (\ref{1710}), we get
$$
 \mu^{\e_n}(\tilde x_{1,\e_n},
\ldots,\tilde x_{l,\e_n},\tilde y_{1,\e_n},\tilde
y_{2,\e_n},\ldots,\tilde y_{m_n,\e_n})\le 
I_{\e_n}(\bar v_n)
$$
$$
=I_{\e_n}(\tilde w_{\tilde y_{1,\e_n}})+I_{\e_n}(\bar u_n)
-I_{\e_n}((1-\hat\chi(\e_nx))(\bar u_n\vee
\tilde w_{\tilde y_{1,\e_n}}))+o(1),
$$
from which
$$
\mu^{\e_n}(\tilde x_{1,\e_n},
\ldots,\tilde x_{l,\e_n},\tilde y_{1,\e_n},\tilde
y_{2,\e_n},\ldots,\tilde y_{m_n,\e_n})- 
\mu^{\e_n}(\tilde x_{1,\e_n},
\ldots,\tilde x_{l,\e_n},\tilde y_{2,\e_n},\ldots,\tilde y_{m_n,\e_n})
$$
\beq
\label{1732}
\le
I_{\e_n} (\tilde w_{\tilde y_{1,\e_n}}
)-I_{\e_n}((1-\hat\chi(\e_nx))(\bar u_n\vee \tilde w_{\tilde
  x_{1,\e_n}}) )+o(1) .
\eeq

Now, let us show that 
\beq
\label{1547}
I_{\e_n}(\tilde w_{\tilde y_{1,\e_n}})\le m_\infty+o(1).
\eeq
Indeed, by Lemma \ref{N16}, $I_{\infty}(\tilde w_{\tilde y_{1,\e_n}})\le
m_\infty,$ and, by Lemma \ref{L2.10}, 
$\{\theta^-_{\e_n}(- w_{\tilde y_{1,\e_n}})$ : $n\in\N\}$ is bounded,
thus 
\begin{eqnarray}
I_{\e_n}(\tilde w_{\tilde y_{1,\e_n}})& = & 
I_{\infty}(\tilde w_{\tilde y_{1,\e_n}})+{1\over
  2}\int_{\R^N}\alpha_{\e_n}(x)
\tilde w_{\tilde y_{1,\e_n}}^2dx
\nonumber \\
&\le & m_\infty+c\int_{\R^N}|\alpha_{\e_n}(x)|\,
 w_{\tilde x_{1,\e_n}}^2dx 
\label{1654}
\end{eqnarray}
holds. 
Moreover, fixing arbitrarily a real number $\zeta>0$, a real number
$R_\zeta>0$ and a natural $n_\zeta$ can be found so that
$ |w_{\tilde y_{1,\e_n}}|_{L^{2^*}(\R^N\setminus B_{R_\zeta}(\tilde
  y_{1,\e_n}))}<\zeta$ and, by Remark \ref{og}, 
$|\alpha_{\e_n}|_{L^\infty(B_{R_\zeta}(\tilde y_{1,\e_n}))}$ $<\zeta$ for all
$ n>n_\zeta$, so
$$
\int_{\R^N}|\alpha_{\e_n}(x)|\, w^2_{\tilde y_{1,\e_n}}dx \le 
\zeta\int_{B_{R_\zeta}(\tilde y_{1,\e_n})}w^2_{\tilde y_{1,\e_n}}dx+
 \int_{\R^N\setminus B_{R_\zeta}(\tilde y_{1,\e_n})}
|\alpha_{\e_n}(x)|\, w^2_{\tilde y_{1,\e_n}}dx, 
$$
therefore, writing $\R^N=\cup_{i=1}^\infty Q_i$, where $Q_i$ are
$N$-dimensional disjoint unit hypercubes,  we deduce, for large $n$,
\begin{eqnarray*}
\int_{\R^N}|\alpha_{\e_n}(x)|\, w^2_{\tilde y_{1,\e_n}}dx &
 \le &\zeta |w|_2^2+\sum_{i=1}^\infty|\alpha_{\e_n}|_{L^{N/2}(Q_i)}|
w_{\tilde y_{1,\e_n}}|^2_{L^{2^*}(Q_i\cap(\R^N\setminus B_{R_\zeta}(\tilde
  y_{1,\e_n})))}
\\
& \le &
(|w|_2^2+\sup_{i\in \N}|\alpha_{\e_n}|_{L^{N/2}(Q_i)})\zeta\le c\zeta
\end{eqnarray*}
with $c>0$ independent of $n$, giving
\beq
\label{1835}
\lim_{n\to \infty} \int_{\R^N}|\alpha_{\e_n}(x)|\,
 w_{\tilde y_{1,\e_n}}^2dx= 0
\eeq
and then (\ref{1547}).

To estimate the second addend in (\ref{1732}), we start observing that 
$| (1-\hat\chi(\e_nx))(\bar u_n\vee
\tilde w_{\tilde y_{1,\e_n}})|\le\delta$ and 
$ (1-\hat\chi(\e_nx))(\bar u_n\vee
\tilde w_{\tilde y_{1,\e_n}})=(1-\hat\chi(\e_nx))(\bar u_n\vee
 (-w_{\tilde y_{1,\e_n}}))$,
so, in view of Remark \ref{CC}, we obtain
\begin{eqnarray*}
I_{\e_n}((1-\hat\chi(\e_nx))(\bar u_n\vee
\tilde w_{\tilde y_{1,\e_n}})) & = &
I_{\e_n}((1-\hat\chi(\e_nx))(\bar u_n\vee
 (-w_{\tilde y_{1,\e_n}})))\\
& \ge &
\tilde c\int_{\R^N}[(1-\hat\chi(\e_nx))|\bar u_n\vee
 (-w_{\tilde y_{1,\e_n}})|]^2dx.
\end{eqnarray*}
Now, considering  that $b_{\bar
  r}:=\inf_{B_{2\bar r}(0)}w>0$, and that, by (\ref{1850})  and
(\ref{e22.2}), $ \supp(\bar u_n^-)^\delta_2\subset B_\rho(\tilde
y_{2,\e_n})\subset  B_{2\bar r}(\tilde
y_{1,\e_n})\subset(\R^N\setminus(B\cup\hat\cN))/\e_n$   we get
$$
I_{\e_n}((1-\hat\chi(\e_nx))(\bar u_n\vee
\tilde w_{\tilde y_{1,\e_n}}))\ge \tilde c\int_{\supp (\bar
  u_n^-)^\delta_2}w^2_{\tilde y_{1,\e_n}}\ge\tilde c \, b_{\bar r}^2|\supp
(\bar  u_n^-)^\delta_2|.
$$

Hence, the proof of (\ref{1857}) is reduced to show the inequality:
\beq
\label{1751}
\liminf_{n\to \infty}
|\supp (\bar  u_n^-)^\delta_2 |>0.
\eeq
Again we argue by contradiction and we assume  (\ref{1751}) false, so
that, up to a subsequence,  
$ \lim_{n\to \infty} |\supp (\bar  u_n^-)^\delta_2 |=0$.
Thus, relation
\beq
\label{quozdiv}
\lim_{n\to \infty}
\|(\bar  u_n^-)^\delta_2 \|/|(\bar  u_n^-)^\delta_2 |_{p+1}=\infty
\eeq
must be true, up to a subsequence. Indeed, otherwise,  the sequence
$(\bar u_n^-)^\delta_2(x-\tilde y_{2,\e_n})$ $/|(\bar u_n^-)^\delta_2|_{p+1}$
would be bounded in $H_0^1(B_\rho(0))$ and then, up to a subsequence,  converging
to a function $\hat u\in H_0^1(B_\rho(0)),$ weakly in
$H_0^1(B_\rho(0))$ and strongly in $L^{p+1}(B_\rho(0))$ 
 so that $\hat u=0$  which is
impossible. 
 As (\ref{quozdiv}) holds, we deduce by a computation
quite analogous  to (\ref{m.2}), 
\begin{eqnarray}
\lim_{n\to \infty}J_{\e_n}((\bar  u_n^-)^\delta_2 )
& = & 
\lim_{n\to \infty}\max_{t>0}J_{\e_n}(t(\bar  u_n^-)^\delta_2  )
\nonumber \\
& = & 
\lim_{n\to \infty}\max_{t>0}J_{\e_n}\left(t\, {(\bar
      u_n^-)^\delta_2 \over |(\bar
      u_n^-)^\delta_2 |_{p+1}}\right)
\nonumber \\
&\ge &\lim_{n\to \infty} J_{\e_n}\left( {(\bar
      u_n^-)^\delta_2 \over |(\bar
      u_n^-)^\delta_2 |_{p+1}}\right)=+\infty.
\label{1543}
\end{eqnarray}

Now, we define $\tilde u_n\in S^{\e_n}_{\tilde
  x_{1,\e_n},\ldots,\tilde x_{l,\e_n}, \tilde y_{2,\e_n},\ldots,\tilde
  y_{m_n,\e_n}}$ as  
$$
\tilde u_n=\hat\chi(\e_nx)\bar u_n
+(1-\hat\chi(\e_nx))([(\bar u_n)_\delta +\sum_{j=3}^{m_n}
(\bar u_n^-)_j^\delta]\wedge \tilde w_{\tilde y_{2,\e_n}})
$$
with $\tilde w_{\tilde y_{2,\e_n}}=(-w_{\tilde y_{2,\e_n}})_\delta
-\theta^-_{\e_n}(-w_{\tilde y_{2,\e_n}} ) 
(w_{\tilde  y_{2,\e_n}}^+ )^\delta
\in S^{\e_n,-}_{\tilde y_{2,\e_n} }.$ 
By definition of $\bar u_n$
$$
I_{\e_n}(\bar u_n)\le I_{\e_n}(\tilde u_n)
$$
holds, on the other hand, for large $n$ we have 
\begin{eqnarray*}
I_{\e_n}(\tilde u_n)-I_{\e_n}(\bar u_n) & \le &
\hspace{-1mm} [I_{\e_n}((\bar u_n)_\delta)+\hspace{-1mm}\sum_{i=1}^lJ_{\e_n}((\bar
u_n^+)_i^\delta)+\hspace{-1mm}\sum_{j=3}^{m_n} 
J_{\e_n}((\bar u_n^-)_i^\delta)+I_{\e_n}(\tilde w_{\tilde y_{2,\e_n}})+o(1)]
\nonumber \\
&  &
-[I_{\e_n}((\bar u_n)_\delta)+\sum_{i=1}^lJ_{\e_n}((\bar
u_n^+)_i^\delta)+\sum_{j=2}^{m_n} 
J_{\e_n}((\bar u_n^-)_i^\delta)] < 0,
\end{eqnarray*}
that is obtained using (\ref{1543})  and computations
analogous to those for getting (\ref{1547}) and (\ref{1710}).
Therefore we are in contradiction, this means (\ref{1857}) holds true
and the proof of  (\ref{e23.2}) is completed. 

\vspace{2mm}

Inequality (\ref{e23.5}) follows straightly from definition of $\cH^\e_l$, $\cK^\e_m$. 

\vspace{2mm}

The proof of (\ref{e23.1}) can be carried out  analogously to that of (\ref{e23.2}), substituting (\ref{10.3}) by (\ref{17.24}).

\qed

\begin{cor}
\label{CPS5.2}
Let $(h_1)$,--,$(h_4)$ hold. Let $\tilde\cN$ and $l$ be as in
Proposition \ref{P22.1}. 
Let $(\e_n)_n$ and $(m_n)_n$ be sequences such that $0<\e_n\tton 0$,
$(m_n)_n \in \N\setminus\{0\}$. 
For all $n\in\N$, let $(\tilde x_{1,\e_n},\ldots, \tilde x_{l,\e_n}, \tilde
y_{1,\e_n},\ldots, \tilde y_{m_n,\e_n})\in\cH_l^{\e_n}\times\cK^{\e_n}_{m_n}$  be a
$(l+m_n)$-tuple for which $\mu^{\e_n}
(\tilde x_{1,\e_n},\ldots,$ $\tilde x_{l,\e_n},\tilde
y_{1,\e_n},\ldots,\tilde
y_{m_n,\e_n})=\mu^{\e_n}_{l,m_n}.$  Then
\beq
\label{1031}
\lim_{n\to\infty}  \min\{\dist(\tilde y_{i,\e_n}, \tilde\cN /\e_n)\ :\ 
i\in\{1,\ldots,m_n\}\}=\infty.
\eeq
Moreover, if $l\ge 1$
\beq
\label{1032}
\lim_{n\to\infty}  \min\{ \dist(\tilde x_{i,\e_n},  \tilde \cN/\e_n)\ :\ 
{i\in\{1,\ldots,l\}}\}=\infty,
\eeq
$$
 \lim_{n\to \infty} \min\{ |\tilde x_{i,\e_n}-\tilde y_{j,\e_n}|\  :\
i\in\{1,\ldots,l\},\ j\in\{1,\ldots,m_n\}\}=\infty,
$$
if $l\ge 2$
\beq
\label{1014}
\lim_{n\to \infty} \min\{ |\tilde x_{i,\e_n}-\tilde x_{j,\e_n}|\  :\
i,j\in\{1,\ldots,l\},\ i\neq j\}=\infty,
\eeq
if $m_n\ge 2$
\beq
\label{1016}
 \lim_{n\to \infty} \min\{ |\tilde y_{i,\e_n}-\tilde y_{j,\e_n}|\  :\
i,j\in\{1,\ldots,m_n\},\ i\neq j\}=\infty.
\eeq
\end{cor}

\begin{prop}
\label{CPS5.4}
Let $(h_1)$,--,$(h_4)$ hold. Let $l,\ (\e_n)_n
 \ (m_n)_n,$ and  $\tilde\cN$ be as in Corollary \ref{CPS5.2}.  Let $((\tilde
x_{1,\e_n},\ldots, \tilde x_{l,\e_n},\tilde y_{1,\e_n},\ldots,
\tilde y_{m_n,\e_n}))_n$ be a sequence of $l+m_n$-tuples, belonging to
$\cH_l^{\e_n}\times\cK^{\e_n}_{m_n},$  verifying (\ref{1031}),
(\ref{1032}) and, if $l,m_n>1,$ also (\ref{1014}), (\ref{1016}).
Let  $u_{l,m_n,\e_n}$ be in $M^{\e_n}_{\tilde
x_{1,\e_n},\ldots, \tilde x_{l,\e_n}, \tilde y_{1,\e_n},\ldots,
\tilde y_{m_n,\e_n}}$.  Then, for all $r>0$
\beq
\label{e23.4}
\lim_{n\to \infty}
(\sup\{|u_{l,m_n,\e_n}(x+\tilde y_{i,\e_n})+w(x)|\ :\ i\in\{1,\ldots,m_n\},\
|x|\le r\})=0,
\eeq
if $l\ge 1$
\beq
\label{e23.3}
\lim_{n\to \infty}
(\sup\{|u_{l,m_n,\e_n}(x+\tilde x_{i,\e_n})-w(x)|\ :\ i\in\{1,\ldots,l\},\
|x|\le r\})=0.
\eeq
\end{prop}

\proof
We prove (\ref{e23.3}); the argument to show (\ref{e23.4}) can be developed 
analogously. 
Without any loss of generality, from now on we put $i=1$.

We set 
$\zeta_n(x) = \zeta \left(|x-\tilde x_{1,\e_n}| - \frac{d_n}{2}\right)$, where
$d_n = \min\{\dist(\tilde x_{1,\e_n},\partial (B\setminus A)/\e_n)$, $|\tilde
x_{1,\e_n} - \tilde x_{i,\e_n}|$, $ i\in\{ 2,...,l\} \}$ 
and $\zeta \in C^\infty (\R , [0,1])$ is such that $\zeta (t) = 0$ 
when $|t| \leq 1/2,$ and $\zeta(t) = 1,$ when $|t| \geq 1$.
 Then, we consider $\zeta_n(x)\ u_{l,m_n,\e_n} (x)$.
By (\ref{1014}) and (\ref{1032}), $d_n \tton \infty$, and for
large $n$ we can write 
$$ 
\zeta_n(x)\ u_{l,m_n,\e_n} (x) = \bar{u}_n(x) + \hat{u}_n (x)
$$
with $\hat{u}_n (x) \in S^{\e_n}_{\tilde x_{2,\e_n},..., \tilde
  x_{l,\e_n},\tilde y_{1,\e_n},\ldots,\tilde y_{m_n,\e_n}}$ and
$\bar{u}_n(x) \in S^{\e_n,+}_{\tilde x_{1,\e_n}}$ 
 such that $\supp\ \hat{u}_n \subset \R^N \setminus$  $B_{\frac{d_n}{2}
   + \frac{1}{2}}(\tilde x_{1,\e_n})$,  $\supp\ \bar{u}_n\subset
 B_{\frac{d_n}{2} - \frac{1}{2}}(\tilde x_{1,\e_n})$.
Now, our reasoning is carried out by proving step by step the following points:

\vspace{2ex}

{\bf (a) } \hspace{1cm} $\lim_{n\to \infty}I_{\e_n}(\bar u_n)= m_\infty;$

{\bf (b) } \hspace{1cm} $ \bar u_n (\cdot + x_{1,\e_n}) \tton w, $  in $ H^1(\R^N)$;

{\bf (c) } \hspace{1cm} $ \bar u_n (\cdot + x_{1,\e_n}) \tton w, $  uniformly in $K, \ \ \forall K \subset \subset B_\rho(0)$ compact. 

\vspace{2ex}
\noindent Actually, once realized these points, we easily conclude. 
Indeed,  the choice of $\delta$ and $\rho$ that implies  $w(x)<\delta$ 
when $|x|>\rho/2$,  the exponential decay of $u_{l,m_n,\e_n}$, the
relations  $u_{l,m_n,\e_n}=\bar u_n$ in $B_{{d_n\over  2}-1}(\tilde
x_{1,\e_n}),$ and $d_n\tton \infty$ allow us to conclude that $\bar
u_n(\cdot +\tilde x_{1,\e_n})$ verifies, for all fixed $R>\rho$ and
for large $n$,  
$$
-\Delta u+a_{\e_n}(x)u=|u|^{p-2}u\quad\mbox{ and }\quad |u|\le c\,
e^{-\bar \eta |x|}\quad\mbox{ in }B_R(0)\setminus B_{\rho/2}(0),
$$
 with a constant $c>0$ independent of $R$.
Hence, taking into account point (b), Lemma \ref{N2}, and the fact that
$a_{\e_n}(\cdot+\tilde x_{1,\e_n})\tton a_\infty$ uniformly in $
B_R(0)\setminus B_{\rho/2}(0)$,  regularity arguments give   
$$
u_n(\cdot +\tilde x_{1,\e_n})\tton w\quad\mbox{ uniformly in }
 B_R(0)\setminus B_{\rho/2}(0)
$$
and this relation, together with point (c), yields  (\ref{e23.3}).

\emph{Proof of Point (a)}

Let us consider $\bar z_n\in M^{\e_n,+}_{\tilde x_{1,\e_n}}$, first of all we want to show that

\beq
\label{1622}
\lim_{n\to \infty}I_{\e_n}(\bar z_n)=\lim_{n\to \infty}I_{\e_n}(\bar
u_n).
\eeq
Arguing as in
(\ref{4.19}) and taking into account the asymptotic decay (\ref{e7.1}) of $u_{l,m_n,\e_n},$
 we obtain
\beq
\label{1301}
I_{\e_n}(\bar u_n)+I_{\e_n}(\hat u_n)=
I_{\e_n}(\bar u_n+\hat u_n)=I_{\e_n}(\zeta_n u_{l,m_n,\e_n})\le
I_{\e_n}(u_{l,m_n,\e_n})+o(1).
\eeq
On the other hand, defining $z_n(x)=\tilde \zeta_n(x)\cdot \bar z_n(x)$, where
$\tilde \zeta_n(x)=\tilde \zeta\left(|x-\tilde x_{1,\e_n}|-{d_n\over
    2}\right)$ with $\tilde \zeta\in C^\infty(\R,[0,1])$,
$\tilde\zeta(t)=1$ if $t\le {1\over 4}$, $\tilde \zeta (t)=0$ if $t\ge
{1\over 2},$ 
by construction we have that $z_n\in S^{\e_n,+}_{\tilde x_{1,\e_n}}$ and $z_n+\hat u_n\in S^{\e_n}_{\tilde x_{1,\e_n},\tilde
  x_{2,\e_n},..., \tilde x_{l,\e_n},\tilde y_{1,\e_n},\ldots,\tilde
  y_{m_n,\e_n}}$. Using (\ref{e7.1}) and  $d_n\tton
  \infty,$ it is not difficult to see that 
$I_{\e_n}(z_n)=I_{\e_n}(\bar z_n)+o(1),$
so we infer
\beq
\label{1300}
I_{\e_n}(u_{l,m_n,\e_n})
\le I_{\e_n}(z_n+\hat u_n)
=I_{\e_n}(z_n)+I_{\e_n}(\hat u_n)
=I_{\e_n}(\bar z_n)+I_{\e_n}(\hat u_n)+o(1).
\eeq
Combining (\ref{1301}) and (\ref{1300}) we get then
$$
I_{\e_n}(\bar u_n)
\le I_{\e_n}(u_{l,m_n,\e_n})+o(1)-I_{\e_n}(\hat u_n)
\le I_{\e_n}(\bar z_n)+o(1),
$$
which, together with $I_{\e_n}(\bar u_n)\ge I_{\e_n}(\bar z_n)$ gives (\ref{1622}).

Now, we have to evaluate the limit as $n\to \infty$ of $I_{\e_n}(\bar z_n)$ and to show it is equal to $m_\infty.$

Arguments analogous to those of
 (\ref{11.58}), (\ref{1654}), and  (\ref{1835}) give:
$$
I_{\e_n}(\bar z_n)
\le I_\infty(w_{\tilde x_{1,\e_n}})+c\int_{\R^N}|\alpha_{\e_n}(x)|\,
w^2_{\tilde x_{1,\e_n}}dx
=m_\infty+o(1).
$$
Therefore, if we prove the reverse inequality, we are done. To this
end, let us consider $\tilde z_n=(\bar 
z_n)_\delta+\theta^+_\infty(\bar z_n)(\bar z_n^+)^\delta$ and observe that 
$$
c\ge I_{\e_n}(\bar z_n)\ge I_{\e_n}(\tilde z_n)=I_\infty(\tilde
z_n)+{1\over 2}\int_{\R^N}\alpha_{\e_n}(x)\tilde z_{n}^2dx,
$$
then, we can argue as in step 1 of Lemma \ref{L10.1} to show that $\|(\tilde
z_n^+)^\delta\|$ is bounded, and, in turn, thanks to  Lemma \ref{L7.1}, that $\|\tilde z_n\|$ is bounded too. Now,  relation
\beq
\label{limalfa}
\lim_{n\to \infty}\int_{\R^N}|\alpha_{\e_n}(x)|\,\tilde z_n^2dx=0,
\eeq
follows by the same argument used in  (\ref{1835}).
Lastly, using (\ref{m<mu}) and (\ref{limalfa}), we get as desired
$$
I_{\e_n}(\bar z_n)\ge I_{\e_n}(\tilde z_n)\ge I_\infty(\tilde
z_n)-{1\over 2}\int_{\R^N}|\alpha_{\e_n}(x)|\, \tilde z_n^2dx\ge
m_\infty-o(1).
$$

\emph{Proof of Point (b)}

Setting $\tilde u_n:=(\bar u_n)_\delta +\theta^+_\infty(\bar
u_n)(\bar u_n^+)^\delta$,  we first show that
\beq
\label{1532}
\lim_{n\to \infty}I_\infty(\tilde u_n)=m_\infty.
\eeq
holds. 

Arguing as before for $\tilde
z_n$, one can prove that $\|\tilde u_n\|$ is bounded;
furthermore, being $\bar u_n\in S^{\e_n,+}_{\tilde x_{1,\e_n}}$,  using the previous point,
(\ref{m<mu}), and (\ref{e7.1}) we obtain
$$
m_\infty + o(1) = I_{\e_n}(\bar u_n)\ge I_{\e_n}(\tilde u_n)\ge I_\infty(\tilde
u_n)-{1\over 2}\int_{\R^N}|\alpha_{\e_n}(x)|\, \tilde u^2_ndx\ge
m_\infty-o(1),
$$
which implies $\lim_{n\to
  \infty}I_{\e_n}(\tilde u_n)=m_\infty$.
So, since 
$$
\lim_{n\to \infty}|I_{\e_n}(\tilde u_n)-I_\infty(\tilde u_n)|
=\lim_{n\to \infty}{1\over 2}\left| 
\int_{\R^N}\alpha_{\e_n}(x)\, \tilde u^2_ndx\right|=0
$$
(\ref{1532}) follows.

Set now $v_n:=\tilde u_n(\cdot +\tilde x_{1,\e_n}): \ (\|v_n\|)_n$ is bounded, so 
 up to a subsequence,   $\tilde v\in H^1(\R^N)$ exist such
that
\beq
\label{1105}
v_n\tton \tilde v \mbox{ in } L^{p+1}_{\loc}(\R^N), \mbox{ a.e. in
}\R^N\mbox{ and weakly in }H^1(\R^N)\mbox{ and }L^{p+1}(\R^N)
\eeq
furthermore, we remark that by the uniform exponential decay of $\bar
u_n$ (inherited from $u_{l, m_n, \e_n}$) more is true: 
\beq
\label{vienne}
 v_n\tton \tilde v \ \ \ \ \ \ \ \ \ \ \ \ \mbox{in}\ \  L^{p+1}(\R^N).
\eeq
Our next goal is to show that $(\tilde v^+)^\delta\not\equiv 0$.

The sequence $(\|(v_n^+)^\delta\|/|(v_n^+)^\delta|_{p+1})_n$ is
bounded, otherwise a direct computation as in (\ref{1309}) would imply
$J_\infty((v_n^+)^\delta)\tton +\infty$, contradicting
(\ref{1532}).
Then, setting $\hat
v_n=(v_n)_\delta+{(v_n^+)^\delta/|(v_n^+)^\delta|_{p+1}}$, 
  a  $\hat v\in H^1(\R^N)$ exists such that $\hat v_n\tton  \hat v$ in
  $L^{p+1}_{\loc}(\R^N)$,  a.e. in $\R^N$ and weakly in $H^1(\R^N)$
  and $L^{p+1}(\R^N)$. 
In particular, we remark that $|(\hat v^+)^\delta|_{p+1}=1$, $\beta_0(\hat v)=0$,
$|\hat v|=|\tilde v|\le\d$ in $\R^N\setminus B_\rho(0)$.
Thus, considering  $ \check v=(\hat v)_\delta+\theta^+_\infty(\hat v)(\hat
v^+)^\delta,$  projection of $|\hat v|$ on $S^{\infty,+}_0,$  
by  the minimality of $m_\infty,$ (\ref{1105}),  and the weakly lower semicontinuity of $I_\infty$ we obtain 
$$
m_\infty
\le I_\infty(\check v)
=I_\infty((\hat v)_\delta+\theta_\infty^+(\hat v)(\hat v^+)^\delta)
\le \liminf_{n\to \infty} I_\infty\left((v_n)_\delta+\theta_\infty^+(\hat
  v){(v_n^+)^\delta\over |(v_n^+)^\delta|_{p+1}}\right)
$$
\beq
\label{1550}\le\liminf_{n\to \infty}I_\infty((v_n)_\delta+(v_n^+)^\delta)
=\lim_{n\to \infty}I_\infty(v_n)=m_\infty.
\eeq
Then, by Lemma \ref{N16}, $\check v=w,$ so, in particular, $(\tilde
v)_\delta=w_\delta$. 
Now, if $(\tilde v^+)^\delta\not\equiv 0$ were not true, from
$v_n\tton \tilde v$ in $L^{p+1}(\R^N)$ and
$(\|(v_n^+)^\delta\|/|(v_n^+)^\delta|_{p+1})_n$ bounded we would infer
$(v_n^+)^\delta\tton 0$ in $H^1(\R^N)$ and then, by (\ref{1550}),
 $J_\infty((\check v^+)^\delta)=J_\infty((\tilde
v^+)^\delta)=0$, contrary to $\check v=w$. Therefore, $(\tilde v^+)^\delta\not\equiv 0$ and, analogously to
 (\ref{1550}), we have
$$
m_\infty
\le I_\infty(\tilde v_\delta+\theta^+_\infty(\tilde v)(\tilde
v^+)^\delta)
\le \liminf_{n\to \infty}I_\infty((v_n)_\delta+\theta^+_\infty(\tilde
v)(v^+_n)^\d)
$$
\beq
\le\liminf_{n\to \infty} I_\infty((v_n)_\d +(v^+_n)^\d)
= \lim_{n\to \infty} I_\infty (v_n)=m_\infty.
\label{1107}
\eeq
which, together with  (\ref{1105}) and (\ref{vienne}), implies $\lim_{n\to 
  \infty}\|v_n\|=\|\tilde v\|$.
  Hence
$v_n\to \tilde v$ in $H^1(B_\rho(0)),\ \theta^+_\infty(\tilde
v)=1$, and, by (\ref{1107}), $m_\infty = I_\infty (\tilde v):$ 
 these facts mean that $\tilde v\in M_0^{\infty,+}$, so, by Lemma
\ref{N16}, $\tilde v=w$, that is
$$
\tilde{u}_n (\cdot + \tilde{x}_{1,\e_n}) \tton w  \ \ \ \ \ \  \mbox{in}\ \  H^1(\R^N). 
$$
As a consequence of the above relation and of the definition of $\tilde{u}_n$,  $\theta^+_\infty(\bar
u_n(\cdot + x_{1,\e_n}))\tton 1,$ and we get as desired
$\bar u_n(\cdot +\tilde x_{1,\e_n})\tton w$ in $H^1(\R^N)$.

\vspace{2ex}
\emph{Proof of Point (c)}

Since $u_{l,m_n,\e_n}\in M^{\e_n}_{\tilde x_{1,\e_n},\ldots,\tilde
  x_{l,\e_n},\tilde y_{1,\e_n},\ldots,y_{m,\e_n}}$ and $d_n\to
\infty$, by Proposition \ref{CPS3.6}, for all $n\in\N$ there exists
$\lambda_n\in\R^N$ such that
\beq
\label{r.1}
I'_{\e_n}(\bar u_n)[\psi]=\int_{B_\rho(\tilde x_{1,\e_n})}
(\bar u_n^+)^\d(x)\psi(x)(\lambda_n\cdot (x-\tilde x_{1,\e_n}))dx\qquad \forall
\psi\in H^1_0(B_\rho(\tilde x_{1,\e_n}))
\eeq
then a standard bootstrap argument (see f.i. \cite{Ber}) allows us to
conclude that $\bar u_n(\cdot +\tilde x_{1,\e_n})\in C^{1,\sigma}(K)$
for all compact sets $K\subset B_\rho(0)$.

Moreover, 

\beq
\label{r.2}
\lim_{n\to \infty}\lambda_n=0.
\eeq
Indeed,  let us assume that  $|\lambda_n|\neq 0$, $\forall
n\in\N$. 
Testing (\ref{r.1}) with $\psi_n(x-\tilde
x_{1,\e_n}):=  ((\l_n/|\l_n|) \cdot (x-\tilde x_{1,\e_n} ) ) \phi(x-\tilde x_{1,\e_n}),$ 
where $\phi \in C_0^\infty (B_\rho(0))$ is radial  and $\phi(x)> 0$ in
$B_\rho(0)$, and considering that $\bar u_n(\cdot+\tilde x_{1,\e_n})\tton w$ in
$H^1(\R^N)$,  we get 
$$
\int_{B_\rho (0)} 
[(\nabla \bar u_n(x+\tilde x_{1,\e_n}) \cdot \nabla \psi_n (x)) 
+ a_{\e_n}(x+x_{1,\e_n})\bar u_n (x+x_{1,\e_n})\psi_n (x)]dx 
$$
$$
- \int_{B_\rho (0)}
|\bar u_n(x+\tilde x_{1,\e_n})|^{p-1}\bar u_n (x+\tilde x_{1,\e_n})\psi_n(x)dx
$$
$$
 = |\l_n| \int_{B_\rho (0)}(\bar u_n^+)^\d(x+\tilde
 x_{1,\e_n})\phi(x)\left(\frac{\l_n}{|\l_n|}\cdot x \right)^2 
 \, dx
$$
$$
=\left(\int_{B_\rho(0)}(w^+(x))^\d\phi(x)\,(e\cdot x)^2dx+o(1)\right)|\l_n|
$$
\beq
\label{5.17}
\ge c|\l_n|,
\eeq 
where $c> 0$ and $e\in\R^N$, $|e|=1$.
On the other hand the left hand side in (\ref{5.17}), as $n\to \infty$,
goes to $I'_\infty(w)[ (e\cdot x)\, \phi]=0$, so (\ref{r.2}) follows.

Hence, by (\ref{r.1}), we can assert that the sequence $\bar u_n(x+
\tilde x_{1,\e_n})$ is bounded in $C^{1,\sigma}(K)$, for all compact
sets $K \subset B_\rho(0)$, and, then, 
$$
\bar u_n(\cdot + \tilde x_{1,\e_n})\tton w \quad\mbox{ uniformly in }
K,\qquad \forall K\subset\subset B_\rho(0),
$$
completing the proof.

\qed

\begin{cor}
\label{Corag}
Let assumptions of Proposition \ref{CPS5.4} hold . Let $\ l,\ \e_n,
 \ m_n, \ (\tilde
x_{1,\e_n},\ldots,$ $\tilde x_{l,\e_n},\tilde y_{1,\e_n},\ldots,
\tilde y_{m_n,\e_n}),$ and $u_{l,m_n,\e_n}$ be as in  Proposition \ref{CPS5.4}. Then,  for large $n$, 
$$
\supp ((u_{l,m_n,\e_n})_i^-)^\delta\subset\subset B_\rho(\tilde y_{i,\e_n})
\qquad\forall i\in\{1,\ldots,m_n\},
$$
and, if $l\ge 1,$
$$
\supp ((u_{l,m_n,\e_n})_i^+)^\delta\subset\subset B_\rho(\tilde x_{i,\e_n})
\qquad\forall i\in\{1,\ldots,l\}.
$$
\end{cor}
\begin{cor}
\label{CPS5.5}

Let assumptions of Proposition \ref{CPS5.4} hold . Let $\ l,\ \e_n,
 \ m_n, \ (\tilde
x_{1,\e_n},\ldots,$ $\tilde x_{l,\e_n},\tilde y_{1,\e_n},\ldots,
\tilde y_{m_n,\e_n}),$ and $u_{l,m_n,\e_n}$ be as in  Proposition \ref{CPS5.4}.

Then, for large $n,$ $\lambda_{\tilde x_{i,\e_n}},\ 
 \lambda_{\tilde y_{j,\e_n}} \in \R^N$, $i\in\{1,\ldots,l\}, \ j\in\{1,\ldots,m_n\}$, exist for which the relation 
\beq
\label{CPS5.22}
-\Delta u_{l,m_n,\e_n} (x)+a_{\e_n}(x) u_{l,m_n,\e_n} (x)= |u_{l,m_n,\e_n} (x)
|^{p-1} u_{l,m_n,\e_n} (x) 
\eeq
$$
+ \sum_{i=1}^l (u_{l,m_n,\e_n}^+)_i^\delta (x)(\lambda_{\tilde x_{i,\e_n}}
\cdot (x-\tilde x_{i,\e_n}))
+ \sum_{j=1}^{m_n} (u_{l,m_n,\e_n}^-)_j^\delta (x)(\lambda_{\tilde y_{j,\e_n}}
\cdot (x-\tilde y_{j,\e_n})) \qquad x\in\R^N,
$$
holds true.
\end{cor}

\proof Relation (\ref{CPS5.22}) can be obtained combining the results of Proposition
\ref{CPS3.6}, Lemma \ref{L7.1}, and Corollary \ref{Corag}.


 \sezione{Proof of the  results }


The max-min method, dispayed in the previous sections, made us find
functions $u_{l,m,\e},$ that are good candidates to be critical
points. Actually, we can already assert they  are solutions of
$(P_\e)$ in the whole space $\R^N$ except the support of their
emerging parts. Furthermore, we know they satisfy  relation
(\ref{CPSe3.9}), so to get the desired conclusion what is needed is to
show that the Lagrange multipliers appearing in (\ref{CPSe3.9}) are
zero. This achievement will be obtained when $\e$ is suitably small. 

Clearly, the same reasoning can be used
to face problem $(\tilde{P}_\e)$ setting $l= 0$ and considering
constant sign functions. 

\vspace{0,5cm}

{\underline {\sf Proof of Theorem \ref{T1}} \hspace{2mm}}

 The strategy for finding infinitely many nodal solutions is the
 following: we fix arbitrarily $h \in \N \setminus\left\{0\right\} ,$
 and we show that for all  $l\in (0,h)$ it is possible to find a
 suitably small $\tilde{\e}_l$   such that { \bf for all} $\e \in (0,
 \tilde{\e}_l)$ and { \bf for all} $m\in \N\setminus
 \left\{0\right\},$  every function $u_{l,m,\e} \in M^\e_{ x
   _{1,\e},\ldots, x_{l,\e}, 
 y_{1,\e},\ldots, y_{m,\e} }$ such that $I_\e(u_{l,m,\e})=
\mu_{l,m}^\e$ is a solution of $(P_\e)$, that is $I'_\e(u_{l,m,\e})=
0.$ Then, of course, Theorem \ref{T1} follows choosing
$\bar\e_h=\min\{\tilde \e_l\ :\ 
l\in\{1,\ldots,h\}\}$.

We argue by contradiction, so, we assume that, for all $n,$ 
$(\e_n)_n \in \R^N,\  m_n \in \N\setminus\{0\}, ( x_{1,{\e_n}},\ldots, x_{l,{\e_n}},
y_{1,{\e_n}},\ldots, y_{m,{\e_n}})\in\cH^{\e_n}_l\times \cK^{\e_n},$
and $ u_{l,m_n,{\e_n}}\ 
$ exist such that
$$
\begin{array}{c}
0<\e_n\tton 0 , \  u_{l,m_n,{\e_n}}\in
M^{\e_n}_{ x_{1,{\e_n}},\ldots, x_{l,{\e_n}},
y_{1,{\e_n}},\ldots, y_{m,{\e_n}} },\\
  I_{\e_n}(u_{l,m_n,{\e_n}})=
\mu_{l,m_n}^{\e_n} , \ \mbox{and} \  I'_{\e_n}(u_{l,m_,{\e_n}})\neq 0,
\end{array}
$$
namely, considering the equality (\ref{CPS5.22})  and setting
$$ 
F_n:=\max\{|\lambda_{ x_{i,\e_n}}|,\ |\lambda_{ y_{j,\e_n}}|\
:\ i\in\{1,\ldots,l\},\ j\in\{1,\ldots,m_n\}\}\qquad\forall n\in\N
$$
where $\lambda_{ x_{i,\e_n}}$, $\lambda_{ y_{j,\e_n}}$
are the Lagrange multipliers  appearing in that equality, the relation $F_n > 0$ holds.

Up to a subsequence, we can suppose that either $m_n\equiv \bar m$ or $(m_n)_n$ is an
increasing sequence. 
Moreover, without any loss of generality, we can assume 
that one between $\lambda_{
  x_{1,\e_n}}$ and $\lambda_{ y_{1,\e_n}}$ is equal to $F_n,$ for all $n\in \N$. 
In what follows we consider $\lambda_{
  y_{1,\e_n}}=F_n$;   assuming that the equality $\lambda_{
  x_{1,\e_n}}=F_n$ holds, the argument can be carried out  analogously.
  
Therefore, we assume $|\lambda_{y_{1,\e_n}}|\neq 0$, for all $n$, and up
to a subsequence 
$$
\lim_{n\to \infty} {\lambda_{y_{1,\e_n}}\over
|\lambda_{y_{1,\e_n}}|}=\lambda.
$$

Let $(\sigma_n)_n$ be a sequence of real positive
numbers, such that
\beq
\label{37.3} 
\lim_{n\to \infty}\sigma_n=0,\qquad \lim_{n\to\infty} \sigma_nm_n=0,\qquad
\lim_{n\to \infty} \frac{\sigma_n}{|\lambda_{y_{1,\e_n}}|} = 0, 
\qquad
\lim_{n\to \infty} \frac{\sigma_n m_n}{|\lambda_{y_{1,\e_n}}|} = 0, 
\eeq
we set, for all $n\in \N$, 
$$
z_n=y_{1,\e_n}+\sigma_n\lambda
$$ 
and, for all $i\in\{1,\ldots,l+m_n\}$ 
$$
z_{i,n}=  x_{i,\e_n}\ \mbox{ if
}i \leq l\ \mbox{  and }
 z_{i,n}=  y_{i-l,\e_n}\ \mbox{ if }i \geq l+1. 
$$

Let be $v_n \in {M}^{\e_n}_{ x_{1,\e_n},\ldots, 
  x_{l,\e_n} , z_n , y_{2,\e_n},\ldots,  y_{m_n,\e_n} }$.
By definition of $\mu^{\e_n}_{l,m_n}$,
\beq
\label{37.1}
I_{\e_n} (v_n) \leq\mu^{\e_n}_{l,m_n} = I_{\e_n}(u_{l,m_n,\e_n}).
 \eeq
On the other hand, a Taylor expansion gives
\begin{eqnarray}
I_{\e_n} (v_n) - I_{\e_n}(u_{l,m_n,\e_n}) 
&=  & I'_{\e_n}(u_{l,m_n,\e_n})[v_n - u_{l,m_n,\e_n}] 
\nonumber\\
& & + \frac{1}{2} \left[
\int_{\R^N} \hspace{-3mm}|\nabla(v_n - u_{l,m_n,\e_n})|^2 dx 
+\int_{\R^N} \hspace{-3mm}  a_{\e_n}(x)(v_n - u_{l,m_n,\e_n})^2 dx \right]
\nonumber \\
& & - \frac{p}{2} \int_{\R^N} |\omega_n|^{p-1}(v_n -u_{l,m_n,\e_n})^2 dx
\label{37.2} 
\end{eqnarray}
with $\omega_n (x) = u_{l,m_n,\e_n}(x) + \tilde{\omega}_n (x)
(v_n -u_{l,m_n,\e_n})(x)$, for a suitable $\tilde{\omega}_n (x) \in [0,1]$.

Our aim is to obtain a contradiction with (\ref{37.1}) proving, by a
careful estimate of the terms of the above expansion, that, for large $n$, 
\beq 
\label{38.1} 
I_{\e_n} (v_n)- I_{\e_n}(u_{l,m_n,\e_n}) > 0 . 
\eeq 
First step is proving that 
\beq
\label{38.2} 
\lim_{n\to \infty}|v_n- u_{l,m_n,\e_n}|_\infty  = 0.
\eeq 
Applying Proposition \ref{CPS5.4} to $v_n$ and $u_{l,m_n,\e_n}$ and considering that  $|z_n-  y_{1,\e_n}|\tton
0,$ we obtain 
\beq 
\label{38.3} 
\lim_{n\to \infty} |v_n -u_{l,m_n,\e_n}|_{L^\infty(\cup_{i=1}^{l+m_n} 
B_{\bar{\rho}}(z_{i,n}))} = 0 \qquad \forall
\bar{\rho} > 0.
\eeq

On the other hand, since $\d < (\bar a /p)^{1/(p-1)}$, when $ -\delta<  u_{l,m_n,\e_n}(x) <
v_n(x)< \d,$  we have
$$
\Delta(v_n - u_{l,m_n,\e_n}) =  (v_n - u_{l,m_n,\e_n}) \left[ a_{\e_n}(x) - 
\frac{|v_n|^{p-1}v_n -|u_{l,m_n,\e_n}|^{p-1}u_{l,m_n,\e_n}}
{v_n - u_{l,m_n,\e_n}} \right]  \geq 0.
$$  
Analogously, if $ -\delta<\ v_n(x) < u_{l,m_n,\e_n} (x) < \d $, we deduce
$$ 
\Delta(v_n -u_{l,m_n,\e_n}) \leq 0. 
$$
Hence, for $\bar{\rho}>\rho$, taking into account the asymptotic
behaviour of $v_n$ and  $u_{l,m_n,\e_n}$ we can say that  the maximum of
$|v_n - u_{l,m_n,\e_n}|$   on 
$\R^N \setminus \cup_{i=1}^{l+m_n} 
B_{\bar{\rho}}(z_{i,n})$ is attained on the boundary $\cup_{i=1}^{l+m_n} 
\partial B_{\bar{\rho}}( z_{i,n})$.
Hence, by (\ref{38.3}),
\beq
\label{39.2}
\lim_{n\to\infty} |v_n -u_{l,m_n,\e_n}|_{L^\infty(\R^N\setminus
  \cup_{i=1}^{l+m_n}  B_{\bar{\rho}}(z_{i,n}))} = 0
\eeq
that, together with (\ref{38.3}), gives (\ref{38.2}).

Set $ s_n = |v_n - u_{l,m_n,\e_n}|_\infty$ and remark that $s_n>0$, for
all $n\in\N$, because $z_n\neq  y_{1,\e_n}$.
Then, we can define
$$ 
\phi_n(x) = \frac{v_n(x) - u_{l,m_n,\e_n}(x)}{s_n}.
$$

Now, let us denote by $\cI$ the set $\{1,\ldots,l+\bar m\}$, if $m_n\equiv\bar m$, and the set $\N$,
if $(m_n)_n$ is unbounded, and put, for all $t  \in \cI,$ 
$n(t )=\min\{n\in\N$ : $m_n\ge t-l \}$.

{\em  {\bf Claim}   A subsequence of $(\phi_n)_n$, still
denoted by $(\phi_n)_n$, exists such that for all  $t  \in \cI$ the
sequence $(\phi_n (x+ z_{t ,n}))_{n\ge n(t )}$ converges in 
$H^1_{\loc}(\R^N)$, and  uniformly on the compact subsets of
$\R^N$, to a solution $\phi$ of the equation
\beq
\label{*39} 
- \Delta \phi + a_\infty\phi = p\ w^{p-1} \phi \qquad \mbox{ in }  \R^N.
\eeq
Moreover, the convergence is uniform with respect to $t\in\cI$.
}

We postpone the proof of this claim to the end of the argument.

Let us, now, observe that Proposition \ref{CPS5.4}
and the choice of $\rho$ allow us to assert that, for large $n$, 
$\supp \ (v_n^-)^\d_{|_{B_\rho(z_n)}}\subset B_\rho(  y_{1,\e_n})$,
hence we deduce
$$
\int_{B_\rho(z_{i,n})} [((|v_n|^+)^\d )^2 - ((|u_{l,m_n,\e_n}|^+)^\d)^2 ]
(x-z_{i,n}) dx 
$$
\beq
\label{39.1}
 = \left\{
\begin{array}{cl}
0 
&  \mbox{ if }i\in\{1,\ldots,l,l+2,\ldots,l+m_n\}
\\
\sigma_n \l \int_{B_\rho(  y_{1,\e_n})} ((v_n^-)^\d)^2 dx   
& \mbox{
  if }i=l+1
\end{array}
\right.
\eeq
because, when $i=l+1$:
\begin{eqnarray*}
\int_{B_\rho(  y_{1,\e_n})} ((v_n^-(x))^\d )^2 (x- 
y_{1,\e_n}) dx & = & 
\int_{B_\rho(  y_{1,\e_n})} ((v_n^-(x))^\d )^2  (x-z_n)dx
\\ & & 
+ \int_{B_\rho(  y_{1,\e_n})}((v_n^-(x))^\d )^2  ( z_n- 
y_{1,\e_n} ) dx\\
& = & \sigma _n \l \int_{B_\rho( y_{1,\e_n})} ((v_n^-(x))^\d )^2  dx.
\end{eqnarray*}
On the other hand, for all $i\in\{1,\ldots,l+m_n\}$ let us evaluate
$$
\int_{B_\rho(z_{i,n})} [((|v_n|^+)^\d )^2 - ((|u_{l,m_n,\e_n}|^+)^\d)^2 ]
(x-z_{i,n}) dx 
$$
\beq 
\label{40.1}
= 2  \int_{B_\rho(z_{i,n})}  (|u_{l,m_n,\e_n}|^+)^\d
[|v_n|-|u_{l,m_n,\e_n}|]
 (x-z_{i,n})dx
 + \int_{B_\rho(z_{i,n})} \cR _{v_n,u_{l,m_n,\e_n}}(x)\,   (x-z_{i,n})dx,
\eeq
where
$$
\cR_{v_n,u_{l,m_n,\e_n}}  = ((|v_n|^+)^\d)^2
-((|u_{l,m_n,\e_n}|^+)^\d)^2 - 
2 (|u_{l,m_n,\e_n}|^+)^\d [|v_n | -|u_{l,m_n,\e_n}| ].
$$
By a direct computation, using the convexity of
the real map 
$t \mapsto g(t)=((t-\d)^+)^2$ and the fact that $\min_{t\in\R}[(t-y)^2 -
\cR(t,y)]=0$  for all fixed $y \in \R$  
(being $\cR(t,y) = g(t) - g(y) - g'(y)(t-y)$), we deduce that
\beq
\label{40*}
0 \leq \cR_{v_n,u_{l,m_n,\e_n}}  \leq (|v_n | -
|u_{l,m_n,\e_n}|)^2 \leq s_n^2.
\eeq
Comparing (\ref{39.1}) and (\ref{40.1}) and using (\ref{40*}), 
we infer that when $i\neq l+1$ 
\beq
\label{40.2}
\left|2 \int_{B_\rho(0)} (|u_{l,m_n,\e_n}(x+z_{i,n})|^+)^\delta 
 [|v_n (x + z_{i,n})| - |u_{l,m_n,\e_n}(x+z_{i,n})|]\, x \,
 dx\right|\le c\,s_n^2, 
\eeq
with $c$ independent of $i$, while, when we consider $i=l+1$,
$$
2 \int_{B_\rho(0)} (u_{l,m_n,\e_n}^-)^\delta (x+  y_{1,\e_n})
 [v_n (x +   y_{1,\e_n}) - u_{l,m_n,\e_n}(x+  y_{1,\e_n})]\, 
x \, dx + O(s_n^2) 
$$
\beq
\label{40.3}
 =\sigma_n \l \int_{B_\rho(0)} ((v_n^-)^\delta (x +   y_{1,\e_n}))^2 dx.
\eeq
Therefore, since $s_n \neq 0$ and Lemma \ref{N2'} and (\ref{*39}) give
the existence of vectors $\tau_i \in \R^N$, $i \in \cI$, such that
\beq
\label{40.4}
\lim_{n\to \infty}\phi_n(x+z_{i,n}) = (\nabla w(x) \cdot
\tau_i)\qquad\forall i\in\cI,
\eeq
by using Proposition \ref{CPS5.4}, the equality (\ref{40.2}) and the
choice of $\rho$, we deduce, for all $i\in \cI\setminus\{l+1\}$, the equality 
$$
0  = 2 \lim_ {n\to \infty}\int_{B_\rho(0)}
(|u_{l,m_n,\e_n}(x+z_{i,n})|^+)^\delta  \phi_n (x + z_{i,n})\, x \, dx 
= 2 \int_{B_\rho(0)} (|w|^+)^\d (\nabla w \cdot \tau_i) \, x\, dx 
$$
$$
= \int_{B_\rho(0)} (\nabla ((w^+) ^\d)^2 \cdot \tau_i)\, x \, dx = 
-  \tau_i \int_{B_\rho(0)} ((w^+)^\d)^2  dx
$$
that implies, for all  $i\in \cI\setminus\{l+1\},$ $\tau_i = 0.$ 
Analogously, when we consider  $y_{1,\e_n}$,
 using (\ref{40.3}), we obtain
$$
\l \left(\lim_{n\to \infty}\frac{\sigma_n }{s_n}
\right)\int_{B_\rho(0)} 
((w^+)^\d )^2  dx = - \tau_{l+1} \int_{B_\rho(0)} ((w^+)^\d
)^2 dx
$$
which gives
\beq
\label{41.1}
\tau_{l+1} = 
- \l \left(\lim_{n\to \infty}\frac{\sigma_n }{s_n}\right).
\eeq
Now, we observe that $\tau_{l+1} \neq 0$. 
In fact, otherwise, from (\ref{40.4}), we would deduce
that
$\lim_{n\to \infty} |\phi_n|_{L^\infty(\cup_{i = 1}^{l+m_n}
  B_{\bar{\rho}}(z_{i,n}))} = 
0$, 
 for all $\bar{\rho} > 0$; moreover, by the argument
 used to obtain (\ref{39.2}), 
  for $\bar{\rho}>\rho$  we  deduce the 
relation $\lim_{n\to \infty}$ $|\phi_n|_{L^\infty(\R^N
  \setminus\cup_{i = 1}^{l+m_n}
  B_{\bar{\rho}}(z_{i,n}))}
= 0$, and, 
hence, $\lim_{n\to \infty} |\phi_n|_\infty = 0$,  contradicting
$|\phi_n|_\infty = 1$. 
Thus, by (\ref{41.1}), $\sigma_n$ and $ s_n$ have the same order, namely
\beq
\label{42.2}
\lim_{n\to \infty} \frac{\sigma_n}{s_n} = \gamma \in \R^+ \setminus
\left\{0\right\}\quad\mbox{and} \quad \tau_{l+1} = -\gamma \l. 
\eeq

Now, let $\tilde{\rho}\in (\frac{3}{4}\rho , \rho)$ be fixed and
consider a cut-off decreasing function $\bar \chi \in C^\infty (\R^+,
[0,1])$ such that $\bar \chi(t) = 1$ if $t \leq \tilde{\rho}$, $\bar
\chi (t) = 0$ if $t \geq \rho$. 
Thus, putting $\bar \chi_{z_{i,n}} (x) = \bar \chi (|x - z_{i,n}|)$
for $i\in\cI$ and $n\ge n(i)$,  we have 
$$
\bar \chi_{z_{i,n}} (x) [|v_n (x)| - |u_{l,m_n,\e_n}(x)|] \in H^1_0
(B_\rho (z_{i,n})), 
$$ 
moreover, for large $n$, by Proposition \ref{CPS5.4},
\beq
\label{42.1}
\left\{
\begin{array}{l}
\supp (|u_{l,m_n,\e_n}|^+)^\delta  \subset  
\cup_{i = 1}^{l+m_n}
  B_{\tilde{\rho}}(z_{i,n})
\\ 
  \supp\ (|v_n|^+)^\d  \subset  
\cup_{i = 1}^{l+m_n}
  B_{\tilde{\rho}}(z_{i,n})
 \\
\bar\chi_{z_{i,n}} (x) = 1  \quad \forall x \in B_{\tilde{\rho}}(z_{i,n}) \ : \  
(|u_{l,m_n,\e_n}(x)|^+)^\d \neq 0.
\end{array}
\right.
\eeq

We are, now, in position of estimating the terms of the expansion
(\ref{37.2}). 
Indeed, considering (\ref{42.1}), (\ref{CPSe3.9}), and (\ref{1224}),
we can write, for large $n$: 
$$
\hspace{-8cm} I'_{\e_n}(u_{l,m_n,\e_n})[v_n -u_{l,m_n,\e_n} ] 
$$
\begin{eqnarray*}
& = & 
I'_{\e_n}(u_{l,m_n,\e_n})
[
(1 - \sum_{i = 1}^{l+m_n}\bar\chi_{z_{i,n}}) 
(v_n - u_{l,m_n,\e_n})
]+ I'_{\e_n}(u_{l,m_n,\e_n}) 
[\sum_{i = 1}^{l+m_n}\bar\chi_{z_{i,n}}(v_n - u_{l,m_n,\e_n}) ]
\\ 
& = &
\sum_{i = 1}^{l+m_n}I'_{\e_n}(u_{l,m_n,\e_n})
[\bar\chi_{z_{i,n}}(v_n - u_{l,m_n,\e_n}) ]
\\
& = &
\sum_{i = 1}^{l+m_n}
\int_{B_{\rho}(z_{i,n})}(|u_{l,m_n,\e_n}|^+)^\delta 
\bar\chi_{z_{i,n}}(v_n - u_{l,m_n,\e_n})
 (\l_{z_{i,n}} \cdot (x- z_{i,n})) dx
\\
& = &
\sum_{i = 1}^{l+m_n}
\int_{B_{\rho}(z_{i,n})}(|u_{l,m_n,\e_n}|^+)^\delta 
(v_n - u_{l,m_n,\e_n})
 (\l_{z_{i,n}} \cdot (x- z_{i,n})) dx.
\end{eqnarray*}

Thus, considering (\ref{40.2}) and (\ref{42.2}), we obtain
$$
\lim_{n\to \infty}  \frac{1}{s_n |\l_{  y_{1,\e_n}}|}
 I'_{\e_n}( u_{l,m_n,\e_n})\left[v_n - u_{l,m_n,\e_n} \right] 
$$
\begin{eqnarray}
& = &\lim_{n\to \infty}\left[ \left( \frac{\l_{ 
      y_{1,\e_n}}}{|\l_{  y_{1,\e_n}}|} \cdot  
\int_{B_{\rho}(0)} (u^-_{l,m_n,\e_n}(x+  y_{1,\e_n}))^\delta \phi_n
(x +  y_{1,\e_n})  \, x \, dx \right) 
\right.
\nonumber \\
& & \left.
\hspace{12mm}+
 \sum_{i \neq l+1} \frac{1}{s_n} \frac{\l_{z_{i,n}}}{|\l_{  y_{1,\e_n}}|} 
O(s_n^2)
\right]
\nonumber \\
 & = & 
-\gamma \int_{B_{\rho}(0)} (w^+)^\d(x) (\nabla
w(x)\cdot \l) (\l \cdot x)\, dx.
\label{44.1}
\end{eqnarray}

Moreover, taking into account (\ref{37.3}) and (\ref{42.2}) and
choosing $R> \rho$ so large that $p \,w^{p-1}(x)<\bar a$ in $\R^N
\setminus B_R(0)$, we can write, up to a subsequence: 
\begin{eqnarray}
&\liminf\limits_{n\to \infty}& \frac{1}{s_n |\l_{  y_{1,\e_n}| }}
  \int_{\R^N} \left[ |\nabla (v_n - u_{l,m_n,\e_n} )|^2 + a_{\e_n}(x)(v_n
    - u_{l,m_n,\e_n})^2 +\right.
\nonumber \\
&  & \left. -p \ |\omega_n|^{p-1} (v_n - u_{l,m_n,\e_n} )^2 \right] dx
\nonumber \\
&\geq & \liminf_{n\to \infty}\frac{1}{s_n |\l_{  y_{1,\e_n}}|} 
\left[\sum_{i = 1}^{l+m_n} \int_{B_R(z_{i,n})} - p\
  |\omega_n|^{p-1} (v_n - u_{l,m_n,\e_n} )^2 dx \right.
 \nonumber \\
& & \hspace{1cm}\left. +\int_{\R^N\setminus
    \cup_{i=1}^{l+m_n}B_R(z_{i,n})}
(a_{\e_n}(x)-p|\omega_n|^{p-1}) (v_n - u_{l,m_n,\e_n} )^2 dx\right]
\nonumber \\
&\geq & \lim_{n\to \infty} -\  \frac{C (l+m_n) s_n}{|\l_{  y_{1,\e_n}}|} = 0.
\label{44.2}
\end{eqnarray}
So, finally, combining (\ref{37.2}), (\ref{44.1}) and (\ref{44.2}), we get
$$
 \liminf_{n\to \infty} \frac{I_{\e_n} (v_n) - I_{\e_n}(u_{l,m_n,\e_n})} 
{s_n |\l_{  y_{1,\e_n}|}} 
 \geq 
 -\gamma \int_{B_{\rho}(0)} (w^+)^\d (\nabla w \cdot \l)
(\l \cdot x)dx > 0 
$$
and, as a consequence, (\ref{38.1}), as desired.

\vspace{ 2mm}

To complete the proof, let us now \textbf{prove the claim}.
It is clear that, if $\cI$ is finite, it is sufficient to prove the
claim for a fixed index $t\in\cI$, with  estimates independent of $t$.
But this is also sufficient in the nontrivial case $\cI=\N$. 
Indeed, in such a case, we can extract from every subsequence
$(\phi_{{(t-1)}_n}(x+z_{t-1,{(t-1)}_{n}}))_n$ a new subsequence
$(\phi_{{t}_n}(x+z_{t,{t}_{n}}))_n$ verifying the claim for
the index $t$. 
Finally, the subsequence $(\phi_{n_n}(x+z_{t,n_{n}}))_n$ is the
subsequence we are looking for, $\forall t\in\N$.
Then, from now on we can consider fixed $t\in\cI$. 
Moreover, to simplify the notation, we indicate with the same symbols
subsequences of a given sequence. 

 Being 
$u_{l,m_n,\e_n} \in {M}^{\e_n}_{  x_{1,\e_n}, 
  \ldots ,  x_{l,\e_n},  y_{1,\e_n},y_{2,\e_n},\ldots ,  y_{m_n,\e_n} 
}$ 
and 
$v_n \in {M}^{\e_n}_{
  x_{1,\e_n},\ldots ,  x_{l,\e_n},z_n,   y_{2,\e_n},\ldots ,  y_{m_n,\e_n} 
}$, by Corollary \ref{CPS5.5}, they verify, respectively,  the
Euler-Lagrange equations  
$$
-\Delta u_{l,m_n,\e_n} (x)+a_{\e_n}(x) u_{l,m_n,\e_n} (x)= |u_{l,m_n,\e_n} (x)
|^{p-1} u_{l,m_n,\e_n} (x) 
$$
\beq
\label{delta1}
+ \sum_{i=1}^l (u_{l,m_n,\e_n}^+)_i^\delta (x)(\lambda_{  x_{i,\e_n}}
\cdot (x-  x_{i,\e_n}))
+ \sum_{j=1}^{m_n} (u_{l,m_n,\e_n}^-)_j^\delta (x)(\lambda_{  y_{j,\e_n}}
\cdot (x-  y_{j,\e_n})),
\eeq
$$
-\Delta v_n (x)+a_{\e_n}(x)v_n (x)= |v_n (x)
|^{p-1} v_n (x) + \sum_{i=1}^l (v_n^+)_i^\delta (x)(\bar\lambda_{  x_{i,\e_n}}
\cdot (x-  x_{i,\e_n}))
$$
\beq
\label{delta2}
+ (v_n^-)_1^\delta(x)(\bar\lambda_{z_n}\cdot(x-z_n))
+ \sum_{j=2}^{m_n} (v_n^-)_j^\delta (x)(\bar\lambda_{  y_{j,\e_n}}
\cdot (x-  y_{j,\e_n})),
\eeq
where we have denoted by 
$\bar\lambda_{  x_{i,\e}}$, $\bar\lambda_{z_n}$, and 
$\bar\lambda_{  y_{j,\e}}$ the Lagrange multipliers related to $v_n$. 
Thus  we have
$$
-\Delta(u_{l,m_n,\e_n}-v_n)+ a_{\e_n}(u_{l,m_n,\e_n}-v_n)=(|u_{l,m_n,\e_n}|^{p-1}u_{l,m_n,\e_n}-|v_n|^{p-1}v_n)
$$
$$
+\sum_{i=1}^l
(u_{l,m_n,\e_n}^+)_i^\delta((\lambda_{  x_{i,\e_n}}-
\bar \lambda_{  x_{i,\e_n}})\cdot(x- 
x_{i,\e_n}))
+\sum_{i=1}^l
[(u_{l,m_n,\e_n}^+)_i^\delta-(v_n^+)_i^\delta](\bar\lambda_{ 
  x_{i,\e_n}}\cdot (x-  x_{i,\e_n}))
$$
$$
+(v_n^-)_1^\delta(\bar \lambda_{z_n}\cdot(z_n-  y_{1,\e_n})) +(u_{l,m_n,\e_n}^-)_1^\delta((\lambda_{  y_{1,\e_n}}-
\bar \lambda_{z_n})\cdot(x-  y_{1,\e_n}))
$$
$$
+
[(u_{l,m_n,\e_n}^-)_1^\delta-(v_n^-)_1^\delta]
(\bar\lambda_{z_n}\cdot (x- y_{1,\e_n}))
$$
\beq
\label{46.1}
+\sum_{j=2}^{m_n}(u_{l,m_n,\e_n}^-)_j^\delta((\lambda_{  y_{j,\e_n}}-
\bar \lambda_{  y_{j,\e_n}})\cdot(x- 
y_{j,\e_n}))+\sum_{j=2}^{m_n}
[(u_{l,m_n,\e_n}^-)_j^\delta-(v_n^-)_j^\delta](\bar\lambda_{ 
  y_{j,\e_n}}\cdot (x-  y_{j,\e_n})).
\eeq

Let us fix $t \in \cI$ and, for $n\ge n(t)$, set 
$$
\hat s_n\ =\ \left\{
\begin{array}{ll}
\max\{s_n,|\lambda_{  x_{1,\e_n}}-\bar \lambda_{z_n}|\}
&\mbox{ if } t =l+1
\\
\max\{s_n,|\lambda_{z_{t ,n}}-\bar\lambda_{z_{t ,n}}|\}
&\mbox{ if }t \in \cI\setminus\{l+1\}
\end{array}
\right.
$$
and
$$
\hat\phi_n(x)= \frac{u_{l,m_n,\e_n}(x +z_{t ,n} )- v_n(x + z_{t ,n})}
{\hat{s}_n }. 
$$
Remark that 
\beq
\label{a}
|\hat{\phi}_n|_\infty \leq 1.
\eeq

Dividing by $\hat s_n$, we deduce from  (\ref{46.1})
$$
-\Delta\hat\phi_n(x)+a_{\e_n}(x+z_{t ,n})\hat\phi_n(x)=
b_n(x+z_{t ,n})\hat\phi_n(x)
$$
$$
+\sum_{i=1}^l
(u_{l,m_n,\e_n}^+)_i^\delta(x+z_{t ,n})\left(
{\lambda_{  x_{i,\e_n}}-
\bar \lambda_{  x_{i,\e_n}}\over \hat s_n}\cdot(x+z_{t ,n}- 
x_{i,\e_n})\right)
$$
$$+
\sum_{i=1}^l
{(u_{l,m_n,\e_n}^+)_i^\delta(x+z_{t ,n})
  -(v_n^+)^\delta_i(x+z_{t ,n})
\over   \hat s_n}
\,
(\bar\lambda_{   x_{i,\e_n}}\cdot (x+z_{t ,n}-  x_{i,\e_n}))
$$
$$
+(v_n^-)^\delta_1(x+z_{t ,n})
\left(\hspace{-1mm} \bar\lambda_{z_n}\cdot {z_n-  y_{1,\e_n}\over
    \hat s_n}\hspace{-1mm}\right)+(u_{l,m_n,\e_n}^-)_1^\delta(x+z_{t ,n})
\left(\hspace{-1mm}
{\lambda_{  y_{1,\e_n}}- \bar \lambda_{z_n}\over \hat s_n}
\cdot(x+z_{t ,n}-  y_{1,\e_n})
\hspace{-1mm}\right)
$$
$$
+{(u_{l,m_n,\e_n}^-)_1^\delta(x+z_{t ,n})
  -(v_n^-)^\delta_1(x+z_{t ,n} )\over \hat s_n}
\, 
(\bar\lambda_{z_n}\cdot (x+z_{t ,n}-  y_{1,\e_n}))
$$
$$
+\sum_{j=2}^{m_n}(u_{l,m_n,\e_n}^-)_j^\delta(x+z_{t ,n})\left(
{\lambda_{  y_{j,\e_n}}-
\bar \lambda_{  y_{j,\e_n}}\over \hat s_n}\cdot(x+z_{t ,n}- 
y_{j,\e_n})\right)$$
\beq
\label{f}
+\sum_{j=2}^{m_n}
{(u_{l,m_n,\e_n}^-)_j^\delta(x+z_{t ,n})
  -(v_n^-)^\delta_j(x+z_{t ,n})
\over   \hat s_n}
\,
(\bar\lambda_{   y_{j,\e_n}}\cdot (x+z_{t ,n}-  y_{j,\e_n}))
\eeq
where
$b_n(x)={|u_{l,m_n,\e_n}(x)|^{p-1}u_{l,m_n,\e_n}(x)-|v_n(x)|^{p-1}v_n(x)\over 
u_{l,m_n,\e_n}(x)-v_n(x)}=p|\hat\omega_n|^{p-1}(x)$,
$\hat\omega_n$ being a function taking its values between the values
of $u_{l,m_n,\e_n}$ and $v_n$, while $b_n(x)=0$ if $u_{l,m_n,\e_n}(x)=v_n(x)$.
We also remark that the relations:
\beq
\label{b}
\frac{|\l_{z_{t ,n}} -\bar{\l}_{z_{t ,n}}
  |}{\hat{s}_n} \leq 1 \quad\mbox{ if }t \in \cI\setminus\{l+1\}, 
\quad\quad  
\frac{|\l_{  y_{1,\e_n}} -\bar{\l}_{z_n}|}{\hat{s}_n} \leq
1\quad \mbox{ if }t =l+1
\eeq
clearly hold true, while the relations
\beq
\label{b2}
\l_{z_{t ,n}}\tton 0, \quad 
\bar \l_{z_{t ,n}}\tton 0\quad \mbox{ if }t \in \cI\setminus\{l+1\},\qquad 
\bar{\l}_{z_n}\tton 0,\quad  \mbox{ if }t = l+1,
\eeq
hold uniformly with respect to the choice of $t$ and can be obtained
arguing as for proving (\ref{r.2}) and  by
using (\ref{delta1}) and (\ref{delta2}) respectively. 
Moreover, since for large $n$
$$ 
z_n-  y_{1,\e_n} = 
\frac{1}{|(v_n^-)_1^\d|_2^2} \int_{B_\rho(  y_{1,\e_n})}\hspace{-2mm} x\ [(v_n^-)_1^\d (x)]^2 dx - 
\frac{1}{|(u_{l,m_n,\e_n}^-)_1^\d|_2^2} \int_{B_\rho(  y_{1,\e_n})} \hspace{-2mm}x \
[(u_{l,m_n,\e_n}^-)_1^\d (x)]^2 dx,
$$
 considering that 
$|(v_n^-)_1^\d (x)-(u_{l,m_n,\e_n}^-)_1^\d (x)| \leq s_n$, 
 $(v_n^-)_1^\d(\cdot+z_n) \tton (w^+)^\d$ and 
$(u_{l,m_n,\e_n}^-)_1^\d $ $(\cdot +   y_{1,\e_n})\tton (w^+)^\d$,
 by a direct computation, we deduce that, for large $n,$
\beq
\label{b'}
|z_n -   y_{1,\e_n}| \leq  c   s_n. 
\eeq

Now, let us observe that by (\ref{f})
$$ 
- \Delta \hat{\phi}_n + a_{\e_n}(x +z_{t ,n}) \hat{\phi}_n(x) 
= b_n(x + z_{t ,n}) \hat{\phi}_n (x ) \qquad \forall
x\in\R^N\setminus\cup_{i=1}^{l+m_n}B_{\rho+\sigma_n|\lambda|}(z_{i,n}-z_{t ,n}), 
$$
where $|b_n | \leq p \d^{p-1}$, so $|\hat{\phi}_n |$ satisfies
$-\Delta|\hat{\phi}_n (x)| + (\bar a - p\, \d^{p-1})|\hat{\phi}_n|\leq 0$ in $
\R^N\setminus\cup_{i=1}^{l+m_n}$ $B_{\rho+\sigma_n|\lambda|}(z_{i,n}-z_{t ,n})$. 
Then, by Lemma \ref{N*}, fixing $\hat
r>\rho+|  y_{1,\e_n}-z_n|$ and taking into account Corollary
\ref{CPS5.2}, we have that, for large $n$, $|\hat\phi_n(x)|  
< c e^{-b|x|}$, $x\in B_{2\hat r}(0)\setminus  B_{\rho + \sigma_n |\l|}
(0)$, for a fixed $b \in (0, \sqrt{\bar a - p \d^{p-1}})$ and $c$
independent on $\hat r$.
Therefore, taking into account (\ref{a}), we infer that there exists
$C>0$ such that $\forall \hat r>0$ and for large $n$
\beq 
\label{g} 
\int_{B_{2 \hat{r}} (0)}|\hat{\phi}_n| dx \leq C. 
\eeq 
Thus,
using (\ref{a}),(\ref{g}),(\ref{b}),(\ref{b2}),(\ref{b'}) and
observing that, by Proposition \ref{CPS5.4}, 
$|u_{l,m_n,\e_n}(x+z_{t ,n})|\tton w (x)$, 
$|v_n(x+z_{t ,n})|\tton w (x)$, 
and 
$|\hat{\omega}_n (x+z_{t ,n})|\tton w(x)$,
 from (\ref{f}) we deduce that there exists $\bar n(\hat r)$ such that 
\beq 
\label{h} 
\int_{B_{2 \hat{r}} (0)} | \Delta \hat{\phi}_n||\hat{\phi }_n| dx \leq
C_1\qquad \forall n>\bar n(\hat r),
\eeq 
where $C_1$ is a real positive constant depending neither on $n$
nor on $\hat{r}$. 
Fix now a mollifying positive function
 $\tilde{\xi} \in\mathfrak{D} (B_2(0))$,
$\tilde{\xi} = 1$ on $B_1(0)$, $\tilde{\xi} \le 1$; setting
$\tilde{\xi}_{\hat{r}}(x) = \tilde{\xi} (\frac{x}{\hat{r}})$,
remark that $\int_{B_{2\hat{r}}(0)}|\Delta \tilde{\xi}_{\hat{r}}|
= C_2 \hat{r}^{N-2}$.
Then, using (\ref{h}), we get the relation 
\begin{eqnarray*}
\int_{B_{\hat{r}}(0)}|\nabla \hat{\phi}_n |^2 
& \leq &
\int_{\R^N} |\nabla \hat{\phi}_n |^2 \tilde{\xi}_{\hat{r}}dx 
 =
- \int_{\R^N} 
\hat{\phi}_n \nabla ( \tilde{\xi}_{\hat{r}}\nabla \hat{\phi}_n )
dx 
\nonumber \\
& = &
- \int_{\R^N} \hat{\phi}_n \tilde{\xi}_{\hat{r}} \Delta\hat{\phi}_n dx 
- \frac{1}{2}\int_{\R^N}  
\nabla \hat{\phi}_n^2 \nabla \tilde{\xi}_{\hat{r}}
dx 
\nonumber \\
& \leq &
\int_{B_{2 \hat{r}} (0)} | \Delta \hat{\phi}_n| |\hat{\phi}_n|dx 
+ \frac{1}{2} 
\int_{B_{2 \hat{r}} (0)} \hat{\phi}_n^2 \Delta \tilde{\xi}_{\hat{r}}
dx
\nonumber\\
& \leq & C_1+{C_3\over 2}\, e^{- 2 b \hat{r}} \hat{r}^{N-2} 
\end{eqnarray*}
showing that $(\hat{\phi}_n)_n$ is bounded in $H^1_{loc}(\R^N)$. 

Hence, we can assume that a function $\phi$ exists such that, 
up to a subsequence, $\hat{\phi}_n\rightharpoonup \phi $ in 
$H^1_{loc}(\R^N)$ and, in view of (\ref{a}), 
$|\hat{\phi}_n - \phi|_{L^q(K)} \to 0$ for all $q < \infty$ and
for all compact sets $K\subset\R^N$. 
Furthermore, we can pass to the limit in (\ref{f}) obtaining
\beq
\label{l}
-\Delta \phi +a_\infty \phi - p w^{p-1} \phi = (w^+)^\d (\lambda' \cdot x),
\eeq
where  $\lambda' = \lim_{n\to \infty }\frac{\lambda_{z_{t ,n}} -\bar{\lambda}_{z_{t ,n}}}{\hat{s}_n}$ if
  $t \in \cI\setminus\{ 1+1\}$ and $\lambda' = \lim_{n\to \infty }\frac{\lambda_{ 
    x_{1,\e_n}} -\bar{\lambda}_{z_n}}{\hat{s}_n}$ if $t =l+1$.

Now, to complete the argument, we only need to show that $\l' = 0$. 
Indeed, in this case, for large $n$, $\hat{s}_n = s_n$,  $\hat{\phi}_n
= \phi_n$ and what above proved for $\hat{\phi}_n$ is just what
asserted in the claim. 
From (\ref{l}), using Lemma \ref{N2'} and Fredholm alternative
theorem, we deduce that $(w^+)^\d (\l'\cdot x)$ must be orthogonal to
$\frac{\partial w}{\partial \nu}$, for all $\nu\in\R^N\setminus\{0\}$,  so,
choosing $\nu = \l'$, we obtain 
\begin{eqnarray*}
0& =& \int_{\R^N} \frac{\partial w}{\partial \l'}\,  (w^+)^\d(\l'\cdot
x)dx = \frac{1}{2}\intr \frac{\partial ((w^+)^\d)^2}{\partial\l'}(\l'\cdot x)dx 
 \\ 
& = & -\frac{1}{2}\int_{\R^N} ((w^+)^\d)^2 \frac{\partial
}{\partial\l'} (\l'\cdot x)dx = -\frac{1}{2} |\l'|^2\intr ((w^+)^\d)^2 dx, 
\end{eqnarray*}
that gives $\l' = 0,$ as desired.

\qed

\vspace{1cm}

{\underline {\sf Proof of Theorem \ref{T3}} \hspace{2mm}}
To prove the theorem, it is sufficient to follow step by step the proof of
Theorem \ref{T1} setting $l=0$  and $A=\emptyset$. 
So, for all suitably small $\e$, one can prove,  for
all $m\in\N\setminus\{0\}$ the existence of a non positive solution of $(\tilde
P_\e)$ with exactly $m$ negative bumps in $\R^N\setminus B/\e$.
The  strict negativity of the solutions follows from Lemma \ref{L1207}.

\qed



\end{document}